\documentclass[12pt]{article}
\usepackage{ct}

\usepackage{paralist}
\usepackage{lmodern}
\usepackage{microtype}
\usepackage{verbatim}
\usepackage[inline]{enumitem}
\usepackage[labelformat=simple,labelsep=none,labelfont=up]{subcaption}

\usepackage{mathtools,thmtools}
\usepackage{tikz}

\usetikzlibrary{positioning}

\usepackage{mathdots}
\usepackage{mleftright} \mleftright
\DeclarePairedDelimiter{\abs}{\lvert}{\rvert}
\DeclarePairedDelimiter{\brac}{[}{]}
\DeclarePairedDelimiter{\set}{\{}{\}}
\DeclarePairedDelimiter{\paren}{(}{)}

\DeclareMathOperator{\arm}{arm}

\newcommand{\EP}{EP}

\DeclareMathOperator{\SVT}{SVT}
\newcommand{\REEtoFT}{\overline{\varphi}}
\DeclareMathOperator{\R}{\mathfrak{R}}
\DeclareMathOperator{\B}{\mathfrak{B}}

\newcommand{\LP}{\mathcal{NIP}}
\newcommand{\tiling}{\mathsf{t}}
\newcommand{\OOTtoFT}{\Phi}
\newcommand{\cOOT}{\mathcal{OOT}}
\newcommand{\cOOE}{\mathcal{OOE}}
\newcommand{\cSF}{\mathcal{SF}}
\newcommand{\cRE}{\mathcal{RE}}
\newcommand{\lmstar}{\lambda^*/\mu^*}

\newcommand{\setz}{\mathbb{Z}}
\newcommand{\vocab}{\emph}

\makeatletter
\newcommand{\ycelldraw}[1]{\tikz[x=\macro@boxdim@YT,y=\macro@boxdim@YT] #1}
\makeatother
\newcommand{\bd}{\ycelldraw{\draw[red, thick] (0,1)--(1,0);}}
\newcommand{\peak}{\ycelldraw{\draw[fill=gray] (0,0)--(0,1)--(1,1)--cycle;}}
\newcommand{\del}{\ycelldraw{\draw[red] (0,0)--(1,1) (1,0)--(0,1);}}
\newcommand{\dg}{*(cyan)}

\makeatletter
\usepackage[boxsize=1.1em,aligntableaux=bottom]{ytableau}
\newskip\@ytabspace \@ytabspace=0ex plus 0.3ex
\def\ytab{\@ifstar\@@ytab\@ytab}
\def\ytabset#1{\globaldefs=-1\ytableausetup{#1}\globaldefs=0}
\newcommand{\@ytab}[1][]{\@@ytab[nobaseline,#1]}
\newcommand{\@@ytab}[1][]{\ytabset{#1}\hspace{\@ytabspace}\@ytab@i}
\newcommand{\@ytab@i}[2][]{\ytableaushort[#1]{#2}\hspace{\@ytabspace}}
\newcommand{\ydiag}[1][]{\ytabset{#1}\hspace{\@ytabspace}\@ydiag@i}
\newcommand{\@ydiag@i}[1]{\ydiagram{#1}\hspace{\@ytabspace}}
\makeatother

\newcommand{\lm}{\ensuremath{{\lambda/\mu}}}

\DeclareMathOperator{\SSYT}{SSYT}
\DeclareMathOperator{\OOT}{OOT}
\DeclareMathOperator{\RPP}{RPP}
\DeclareMathOperator{\rpp}{rpp}
\newcommand{\rppq}[1]{\ensuremath{{\rpp_{#1}(q)}}}
\newcommand{\rpplq}{\rppq{\lambda}}
\newcommand{\rpplmq}{\rppq{\lm}}
\newcommand{\rppratio}{\ensuremath{{\frac{\rpplmq}{\rpplq}}}}

\renewcommand{\vec}{\mathbf}
\newcommand{\vx}{\vec{x}}
\newcommand{\vy}{\vec{y}}
\newcommand{\vz}{\vec{z}}
\newcommand{\vd}{\vec{d}}
\newcommand{\va}{\vec{a}}
\DeclarePairedDelimiterX\mpar[2]{(}{)}{#1\;\delimsize\vert\;#2}
\newcommand{\eval}[1]{\left. #1\right\rvert}

\newcommand{\cale}{\mathcal{E}}
\newcommand{\calr}{\mathcal{R}}

\newcommand{\ssf}{\mathsf{f}}

\newcommand{\hash}{\mathbin{\#}}
\makeatother

\newtheorem*{cor1}{Corollary}

\specs{1}{2 (1)}{2022}{}

\dateline{Sep 22, 2020}{Sep 12, 2021}{March 31, 2022}



\title{On the Okounkov--Olshanski formula \\ for standard tableaux of skew shapes}

\author[1]{Alejandro H. Morales\thanks{Partially supported by the NSF Grant DMS-1855536.}}
\author[2]{Daniel G. Zhu}

\affil[1]{Department of Mathematics and Statistics, UMass, Amherst, MA~01003, U.S.A.\newline
\email{ahmorales@math.umass.edu}%
}

\affil[2]{Massachusetts Institute of Technology, Cambridge, MA~02139, U.S.A.\newline
\email{zhd@mit.edu}%
}

\begin{document}

\maketitle


\begin{abstract}
 The classical hook length formula counts the number of standard tableaux of straight shapes. In 1996, Okounkov and Olshanski found a positive formula for the number of standard Young tableaux of a skew shape. We prove various properties of this formula, including three determinantal formulas for the number of nonzero terms, an equivalence between the Okounkov--Olshanski formula and another skew tableaux formula involving Knutson--Tao puzzles, and two $q$-analogues for reverse plane partitions, which complement work by Stanley and Chen for semistandard tableaux. We also give several reformulations of the formula, including two in terms of the excited diagrams appearing in a more recent skew tableaux formula by Naruse. Lastly, for thick zigzag shapes we show that the number of nonzero terms is given by a determinant of the Genocchi numbers and improve on known upper bounds by Morales-Pak-Panova on the number of standard tableaux of these shapes.
\end{abstract}



\section{Introduction}
Young tableaux are fundamental objects in algebraic and enumerative combinatorics and discrete probability (see \cite{Romik}). While standard Young tableaux, developed by Frobenius in 1903, were first used to study representations of the symmetric group (see \cite{saganbook}), extensions of the concept have been found numerous applications elsewhere: semistandard tableaux, an extension of standard tableaux, are inherent in the representation theory of general linear groups, and they define Schur functions, which are one of the key bases of the ring of symmetric functions (see \cite[Ch. 7]{EC2}). A further extension of these ideas to skew shapes yields a rich theory involving counting permutations by descents, the jeu de taquin relation, and the Littlewood-Richardson coefficients (see \cite[Ch. 2]{EC1}, \cite[Ch. 7]{EC2}).

In 1954, Frame, Robinson, and Thrall \cite{frthlf} discovered the hook length formula, a deceptively simple expression that counts the number of standard Young tableaux $f^\lambda$ of a certain shape $\lambda$:
\[f^\lambda = \frac{\abs{\lambda}!}{\prod_{u \in [\lambda]} h(u)},\]
where $h(u) = \lambda_i - i + \lambda'_j -j +1$ is the {\em hook length} of the cell $u=(i,j)$.
The structural simplicity of the formula leads to a wide variety of proofs and applications; for instance probabilistic \cite{GNW} and bijective \cite{NPS} proofs and yielding shapes $\lambda$ for which $f^\lambda$ is maximized when $\abs{\lambda}$ is kept fixed \cite{vershikkerov1985}. In 1971, Stanley \cite{stanley1971theory} found a $q$-analogue of the hook length formula for the generating function of semistandard tableaux:
\[s_\lambda(1,q,q^2,\ldots) = q^{b(\lambda)} \prod_{u \in [\lambda]} \frac{1}{1-q^{h(u)}},\]
where $b(\lambda)=\sum_{i} \binom{\lambda'_i}{2}$. Up to the power of $q^{b(\lambda)}$, the RHS above also gives a $q$-analogue for the generating function of reverse plane partitions of shape $\lambda$ that we denote by $\rpplq$.

Considering skew shapes \lm, there is no known product formula that gives the number  $f^{\lm}$ of standard Young tableaux of skew shape. However, there are recent formulas for $f^{\lm}$ as nonnegative sums of products, indexed by combinatorial objects, that come from rules for \vocab{equivariant Littlewood--Richardson coefficients}. In particular, Okounkov and Olshanski \cite{okounkovolshanski98} discovered the following formula, which will be our focus.
\begin{theorem}[Okounkov--Olshanski \cite{okounkovolshanski98}] \label{thm:oof}
\begin{equation} \label{eq:oof} \tag{OOF}
f^\lm = \frac{\abs{\lm}!}{\prod_{u\in[\lambda]} h(u)} \sum_{T \in \SSYT(\mu, d)} \prod_{u \in [\mu]} (\lambda_{d + 1 - T(u)} - c(u)),
\end{equation}
where $c(u) = j-i$ is the {\em content} of the cell $u=(i,j)$, $d=\ell(\lambda)$, and $\SSYT(\mu,d)$ is the set of SSYT of shape $\mu$ with entries $\leq d$.
\end{theorem}

Here, similar to Morales-Pak-Panova's study \cite{mpp1,mpp2,mpp3,mpp4} of the Naruse hook length formula and its \vocab{excited diagrams} (see \eqref{eq:nhlf} and Section~\ref{sec:finalnaruse}), we prove various properties of the Okounkov--Olshanski formula.

\subsection{Number of nonzero terms}
We examine properties of nonzero terms in \eqref{eq:oof}, allowing their number, denoted $\OOT(\lm)$, to be counted by determinants.
\begin{theorem} \label{thm:ooterms}
The number of nonzero terms of the Okounkov--Olshanski formula of the shape $\lm$ is
\[\OOT(\lm) = \det\brac*{\binom{\lambda_i - \mu_j + j - 1}{i-1}}_{i,j=1}^d = \det\brac*{\binom{\lambda'_i}{\mu'_j + i - j}}_{i,j=1}^{\mu_1}.\]
\end{theorem}

We also show in Theorem~\ref{cor:oolascouxpragacz} that the number $\OOT(\lm)$ satisfies a \vocab{Lascoux--Pragacz} \cite{lascoux1988ribbon} determinantal identity by using the recent \vocab{Hamel--Goulden} \cite{hamelgoulden1995} determinantal identities for Macdonald's \vocab{ninth variation Schur functions} \cite{Mac,bc20,fk20}. 

\subsection{Number of nonzero terms and bounds for standard tableaux of thick zigzags}

Theorem~\ref{thm:ooterms} and Theorem~\ref{cor:oolascouxpragacz} allow  $\OOT(\lm)$ to be evaluated in certain special cases. Most notably, in the case of (thick) zigzag skew shapes $\lm =  \delta_{n+2k}/\delta_n$ where $\delta_n=(n-1,n-2,\ldots,1)$ denotes the staircase shape. For these shapes, $f^{\delta_{n+2}/\delta_n}=E_{2n+1}$ where $E_n$ denotes the $n$th Euler number \cite[\href{http://oeis.org/A000111}{A000111}]{oeis} and $f^{\delta_{n+2k}/\delta_n}$ is given by a determinant of Euler numbers. We show that $\OOT(\delta_{n+2k}/\delta_n)$ is given by (determinants) of 
 \vocab{Genocchi numbers} $G_{2n}$ \cite[\href{http://oeis.org/A110501}{A110501}]{oeis}.
 
\begin{theorem} \label{thm:ootthickzigzag}
For positive integers $k$ and nonnegative integers $n$ we have that 
\[\OOT(\delta_{n+2k}/\delta_n) = \prod_{i=1}^k \frac{(2(n+i+k-1))!}{(2i-1)!} \det \brac*{\hat{G}_{2(n+i+j-1)}}_{i,j=1}^{k},\]
where $\hat G_{2n} = G_{2n}/(2n)!$. 
In particular, for the zigzag $\delta_{n+2}/\delta_n$ we have that $\OOT(\delta_{n+2}/\delta_n)=G_{2n+2}$.
\end{theorem}

Since Euler numbers and Genocchi numbers are proportional ($E_{2n+1}= G_{2n+2} \cdot 2^{2n+1}/(2n+2)$), as a corollary we obtain that  $f^{\delta_{n+2k}/\delta_n}$ and $\OOT(\delta_{n+2k}/\delta_n)$ are also proportional (see Corollary~\ref{cor: rel SYTzigzag OOTzigzag}):
\begin{equation} \label{eq:curiuos identity}
f^{\delta_{n+2k}/\delta_n} \,=\, 2^{k(2n+2k-1)}\,(k(2n+2k-1))!\, \prod_{i=1}^k \frac{(2i-1)!}{(2n+2i+2k-2)!}\,\OOT(\delta_{n+2k}/\delta_n).
\end{equation}
We use this curious proportionality identity to give bounds for  $f^{\delta_{2k}/\delta_{k}}$ \cite[\href{http://oeis.org/A278289}{A278289}]{oeis}, an asymptotic problem  first studied in \cite{mpp_asymptotics}.

\begin{theorem} \label{thm: bounds thick zigzags}
Let $k$ be an even nonnegative integer and $n = \abs{\delta_{2k}/\delta_{k}}$, then
\begin{equation} \label{eq:bounds thick zigzags}
\frac{1}{2} - \frac{9}{2} \log 2 + 2 \log 3 + o(1) \,\leq\, \frac{1}{n}\log f^{\delta_{2k}/\delta_{k}} - \frac{1}{2} \log n \,\leq\, \frac{1}{2} - \frac{13}{2} \log 2 + \frac{7}{2} \log 3 + o(1).
\end{equation}
\end{theorem}

The RHS above is $\approx -0.4219$ and the LHS is $\approx -0.1603$. The lower bound is weaker but the upper bound is sharper than previously known bounds of $\approx -0.3237$ and $\approx -0.0621$ obtained in \cite{mpp_asymptotics} using \eqref{eq:nhlf}. The best known bounds are $\approx -0.2368$ and $\approx -0.1648,$ computed by Pak in \cite{P}. The existence of such a constant was proved in \cite{MPT} and its conjectured value is $\approx -0.1842$, \cite{mpp_asymptotics}.

\subsection{Reformulations of the Okounkov--Olshanski formula}
We give bijections between the tableaux in $\SSYT(\mu,d)$ contributing to the RHS of \eqref{eq:oof} and a zoo of objects including \vocab{flagged skew tableaux} that appeared in \cite{LP} in a different context, \vocab{lozenge tilings}, and variations of excited diagrams (see Figures~\ref{fig:oopaths},\ref{fig:schematic}). The latter include a variation very similar to the excited diagrams of Naruse--Ikeda \cite{ikedanaruse2009excited} and Knutson--Miller--Yong \cite{KMY}: we start with cells of $[\mu]$ and apply excited moves but possibly beyond $[\lambda]$. The other new variation is called \vocab{reverse excited diagrams} that we denote by $\calr\cale(\lm)$: we start with cells of $[\lm]$, viewed as a skew shifted shape, and apply reverse excited moves. As an application of these bijections we give four reformulations of \eqref{eq:oof} (Corollaries~\ref{cor:edform}--\ref{cor:flag-form}). One of them (Corollary~\ref{cor:redform}) is similar to the Naruse hook length formula since it is in terms of reverse excited diagrams. 

\begin{cor1}[Okounkov--Olshanski - reverse excited diagram formulation]
\[
f^{\lm} \,=\, \frac{|\lm|!}{\prod_{u\in [\lambda]} h(u)} \sum_{D\in \calr\cale(\lm)}  \prod_{(i,j) \in B(D)} (\lambda_i+d-j),
\]
where $d=\ell(\lambda)$, $B(D)$ are certain cells of $[\lm]$ (viewed as a shifted skew shape) associated to $D$ and $\lambda_i+d-j$ equals the arm-length of the cell $(i,j)$.
\end{cor1}

Another of the recent formulas for $f^\lm$ coming from geometry is in terms of Knutson--Tao puzzles (see Section~\ref{sec:KT-puzzles} for notation). 

\begin{theorem}[{Knutson--Tao \cite{knutsontao2003puzzles}, see also \cite[\S 9.4]{mpp1}}]
\[f^\lm = \frac{\abs{\lm}!}{\prod_{u\in[\lambda]} h(u)} \sum_{P \in \Delta^{\lambda\mu}_\lambda} \prod_{p \in \lozenge(P)}\operatorname{ht}(p). \]
\end{theorem}

By extending the bijections mentioned earlier, we show that the Knutson--Tao formula and \eqref{eq:oof} are, despite their appearance, essentially the same.

\begin{theorem}
\label{thm:ooiskt}
The Okounkov--Olshanski and Knutson--Tao formulas for $f^\lm$ are term-by-term equivalent.
\end{theorem}

\subsection{\texorpdfstring{$q$}{q}-analogues for skew reverse plane partitions}
There has also been work to find $q$-analogues of the Okounkov--Olshanski formula. Chen and Stanley \cite{chenstanley2016qoo} proved the following result for the generating function of skew semistandard tableaux. For a skew shape $\lm$ and $u$ in $[\mu]$, let $w(u,k)=\lambda_{d+1-k}-c(u)$.

\begin{theorem}[Chen-Stanley \cite{chenstanley2016qoo}] \label{thm: Chen-Stanley q-analogue}
\begin{equation} \label{eq:qoof-ssyt}
\frac{s_\lm(1,q,q^2,\ldots)}{s_\lambda(1,q,q^2,\ldots)} = \sum_{T \in \SSYT(\mu, d)} \prod_{u \in [\mu]} q^{T(u)-d}(1 - q^{w(u,T(u))}).
\end{equation}
\end{theorem}

We give two $q$-analogues of the Okounkov--Olshanski formula for the generating function of reverse plane partitions of skew shape $\rpplmq \coloneqq\sum_{T \in \RPP(\lm)} q^{|T|}$. In contrast with straight shapes, this generating function is nontrivially different from $s_{\lambda/\mu}(1,q,q^2,\ldots)$. We prove our analogues both using factorial Grothendieck polynomials or using a determinantal identity by Krattenthaler \cite{krattenthaler89}. 

\begin{theorem} \label{thm:oorpptableaux}
\begin{equation} \label{eq:qoof-rpp}
\rppratio =  \sum_{T \in \SSYT(\mu, d)} q^{p(T)} \prod_{u \in [\mu]} (1 - q^{w(u,T(u))}),
\end{equation}
where  $p(T)=\sum_{u\in[\mu],m_T(u)\leq k<T(u)} w(u,k)$, and $m_T(u)$ is the minimum $k\leq T(u)$ such that replacing $T(u)$ by $k$ still gives a semistandard tableaux.
\end{theorem}

\begin{theorem} \label{thm:oorpptalt}
\[\rppratio = \sum_{T \in \SSYT(\mu, d)} q^{p^*(T)} \prod_{u \in [\mu]} (1 - q^{w(u, T(u))}),\]
where $p^*(T) = \sum_{u \in [\mu], T(u) < k \leq \min(d, M_T(u))} w(u, k)$, and $M_T(u)$ is the maximum $k\geq T(u)$ such that replacing $T(u)$ by $k$ still gives a semistandard tableaux.
\end{theorem}

These $q$-analogues can be stated using variants of excited diagrams (Theorems \ref{thm:oorppexcited} and \ref{thm: 2nd q-analogue rev excited} respectively) highlighting similarities between our $q$-analogues and a reverse plane partition $q$-analogue of the Naruse hook length formula by Morales, Pak and Panova (see \cite[Cor. 6.17]{mpp1}).

\subsection*{Outline}
 In Section \ref{sec:prelim} we introduce notation and definitions used in the paper. In Section \ref{sec:background} we introduce the Okounkov--Olshanski formula in a context that highlights its connection to other rules for $f^\lm$. Section \ref{sec:posterms} contains bijections between nonzero terms in the Okounkov--Olshanski formula and other combinatorial objects. Section~\ref{sec:reformulations} uses these bijections to give reformulations of \eqref{eq:oof}, including one in terms of Knutson--Tao puzzles, and four others in terms of (reverse) excited diagrams, lozenge tilings and flagged skew tableaux. Section~\ref{sec: det formulas for oot} has the proof of the determinantal identities for $\OOT(\lm)$ in Theorem~\ref{thm:ooterms} and \ref{cor:oolascouxpragacz}. Section~\ref{sec:special} evaluates $\OOT(\lm)$ for special classes of skew shapes including \vocab{slim shapes}, rectangles, and thick zigzags $\delta_{n+2k}/\delta_n$. For the latter, this section also includes asymptotic bounds to $f^{\delta_{2k}/\delta_{k}}$. In  Section \ref{sec:qoo} describes the reverse plane partition $q$-analogues of the Okounkov--Olshanski formula. We conclude with final remarks and open questions in Section~\ref{sec: final remarks}.

\section{Preliminaries} \label{sec:prelim}

\subsection{Young diagrams and tableaux} \label{subsec:prelim_tableaux}
Let $\lambda = (\lambda_1, \lambda_2, \ldots, \lambda_d)$ be an integer partition. We denote the \vocab{size}, \vocab{length}, \vocab{conjugate paritition}, and \vocab{Young diagram} of $\lambda$ by $\abs{\lambda}$, $\ell(\lambda)$, $\lambda'$, and $[\lambda]$, respectively. Given a cell $u = (i,j) \in [\lambda]$, define the \vocab{content} $c(u) = j - i$, the \vocab{arm length} $\arm(u)=\lambda_i-i+1$, and the \vocab{hook length} $h(u) = \lambda_i + \lambda'_j - i - j + 1$. A \vocab{skew partition} is denoted by $\lm$ for $[\mu] \subseteq [\lambda]$ and its \vocab{skew Young diagram} by $[\lm]$. In what follows, $\lambda,\mu,\rho,\xi$ denote partitions, with $[\mu] \subseteq [\lambda]$, and $\lambda$ always having $d$ parts. We will use $\theta$ to denote an arbitrary skew partition. 

For a strict partition $\lambda^*=(\lambda_1,\lambda_2,\ldots,\lambda_d)$ with $\lambda_1 > \lambda_2>\cdots>\lambda_d$, its {\em shifted Young diagram} $[\lambda^*]$ is obtained from the Young diagram by shifting row $i$ to start at position $(i,i)$. We can similarly define shifted skew shapes. Given an ordinary skew shape $\lm$ of length $d$, we denote by $\lmstar$ the shifted skew shape $(\lambda_1+d-1,\lambda_2+d-2,\ldots,\lambda_d)/(\mu_1+d-1,\mu_2+d-2,\ldots,\mu_d)$.

Given a skew partition $\theta$, a \vocab{reverse plane partition (RPP)} of shape $\theta$ is a filling $T$ of the boxes of $[\theta]$ with nonnegative integers so that all rows and columns of $[\theta]$ have weakly increasing entries going left to right or top to bottom. The set of all reverse plane partitions of shape $\theta$ is denoted by $\RPP(\theta)$. We denote the generating function of RPP of shape $\theta$ by
\[
\rppq{\theta} \coloneqq \sum_{T \in \RPP(\theta)} q^{|T|},
\]
where $\abs{T}$ denotes the sum of the entries in $T$.

A \vocab{semistandard Young tableau (SSYT)}  is a reverse plane partition of positive entries such that all columns are strictly increasing. The set of all semistandard Young tableaux of shape $\theta$ is denoted by $\SSYT(\theta)$. Let $\SSYT(\theta, L)$ be the set of semistandard Young tableaux with all entries at most $L$. For a straight shape $\theta=\lambda$, the size $\abs{\SSYT(\lambda,L)}$ is given by Stanley's hook content formula (e.g. see \cite[Cor. 7.21.4]{EC2}). A \vocab{standard Young tableau (SYT)} is a semistandard Young tableau of shape $\theta$ with entries exactly $\{1, 2, \ldots, \abs{\theta}\}$. The number of standard Young tableaux of shape $\theta$ is $f^\theta$.

Given a skew shape $\theta$ with $d$ rows and a sequence ${\bf f}=(\ssf_1,\ldots,\ssf_{d})$ of weakly increasing nonnegative integers, a \vocab{flagged tableau} is a semistandard Young tableau of shape $\theta$ such that every entry in row $i$ is at most $\ssf_i$. These tableaux were first studied by Lascoux and Sch\"utzenberger \cite{LS82}, Wachs \cite{wachs85}, and Gessel--Viennot \cite{gv89}; see also \cite{proctor2017row}. Let $\SSYT(\theta,{\bf f})$ denote the set of such tableaux.

A \vocab{set-valued semistandard Young tableau} is a filling  of $[\theta]$ with nonempty sets of positive integers, such that for every way to choose an element from the entry of each cell, the chosen elements form a valid semistandard tableau. Set-valued tableaux were introduced by Buch in \cite{AB}. We denote the set of set-valued tableaux of shape $\mu$ by $\SVT(\mu)$. Let $\SVT(\mu,L)$ be the set of set-valued semistandard tableaux of shape $\mu$ with entries at most $L$.

\subsection{Schur functions and generalizations} \label{sec: generalizations Schur functions}
Given a partition $\lambda=(\lambda_1, \lambda_2, \ldots, \lambda_d)$ and an infinite sequence of variables $\vy = (y_1, y_2, \ldots)$, define $\vy_\lambda = (y_{\lambda_1 + d}, y_{\lambda_2 + d - 1}, \ldots, y_{\lambda_d + 1})$.

If $\vx = (x_1, x_2, \ldots)$ is an infinite sequence of variables, define the \vocab{Schur function}
\[s_\theta(\vx) \,\coloneqq\, \sum_{T\in\SSYT(\theta)} \vx^T,\]
where $\vx^T = x_1^{\text{\#1s in $T$}}x_2^{\text{\#2s in $T$}}\cdots$. If $\vx = (x_1, x_2, \ldots, x_k)$ is a finite sequence of variables and $\va = (a_1, a_2, \ldots)$ is an infinite sequence of variables, define the \vocab{factorial Schur function}
\[s_\theta\mpar{\vx}{\va} \coloneqq \sum_{T\in\SSYT(\theta)} \prod_{u\in[\theta]} (x_{T(u)} - a_{T(u) + c(u)}).\]
It is known \cite{molevsagan1999} that for ordinary partitions $\lambda$
\begin{equation}
s_\lambda\mpar{\vx}{\va} = \frac{\det \brac*{(x_i - a_1)(x_i - a_2) \cdots (x_i - a_{\lambda_j + k - j})}_{i,j=1}^k}{\prod_{i<j}(x_i - x_j)}, 
\label{eq:factschurdet}
\end{equation}
where we take $\lambda_j = 0$ if $j > d$. Note that from this analogue of the bialternant definition of Schur functions we conclude that $s_\lambda\mpar{\vx}{\va}$ is symmetric in $\vx$.

\begin{remark}
The factorial Schur function $s_\lambda\mpar{\vx}{\va}$ equals a \vocab{double Schubert polynomial} of a \vocab{Grassmannian permutation} associated to $\lambda$ \cite{LS82} as shown in \cite{Lascoux_lectures}.
\end{remark}

Given a skew shape $\theta$ and a two-dimensional array of variables $\vz = (z_{i,j})_{i \geq 1, j \in \setz}$, define the \vocab{ninth variation Schur function}

\begin{equation} \label{eq:def ninth variation}
s^\dagger_{\theta} (\vz) \coloneqq \sum_{T \in \SSYT(\theta)} \prod_{u \in [\theta]} z_{T(u), c(u)},
\end{equation}
 defined by Macdonald \cite{Mac} and recently studied in \cite{NNSY,bc20,fk20}.

\begin{remark}
If $\vz$ is of the form
\[z_{i,j} = \begin{cases} x_i - a_{i+j} & i \leq d \\ 0 & \text{else}\end{cases},\] then $s^\dagger_\theta(\vz) = s_\theta\mpar{x_1, x_2, \ldots, x_d}{\va}$. Thus the $s^\dagger$ reduce to factorial Schur functions which in turn reduce to Schur functions. 
\end{remark}

Let $a \oplus b = a + b - ab$ and $\ominus a = \frac{a}{a-1}$ be the unique value so that $a \oplus (\ominus a) = 0$. If $\vx = (x_1, x_2, \ldots, x_k)$ is a finite sequence of variables and $\va = (a_i)_{i\in \mathbb{Z}}$ is an infinite sequence of variables, define the \vocab{factorial Grothendieck polynomial} \cite{KMY, mcnamara2006} to be
\[G_\mu \mpar{\vx}{\va} \,\coloneqq\, \sum_{T\in \SVT(\mu,d)} (-1)^{\abs{T} - \abs{\mu}} \prod_{\substack{u\in [\mu]\\r\in T(u)}} (x_r \oplus a_{r+c(u)}).\]

\begin{remark}
Just as factorial Schur functions, $G_\mu \mpar{\vx}{\va}$ equals a \vocab{double Grothendieck polynomial} of a \vocab{Grassmannian permutation} associated to $\mu$ \cite{mcnamara_addendum}. The latter were defined more generally for all \vocab{vexillary permutations} by Knutson--Miller--Yong \cite{KMY} before \cite{mcnamara2006}.
\end{remark}

\subsection{Excited diagrams}
Fix a skew shape $\lambda/\mu$. Given a subset $D$ of $[\lambda]$, consider a subset $D'$ obtained from $D$ by applying the following move to an element of $D$ (represented in blue):
\[\ytab[centertableaux]{{\dg}{},{}{}} \; \longrightarrow\; \ytab[centertableaux]{{}{},{}{\dg}} \label{eq:augmentation} \tag{$\searrow$} \]
(This is only allowed if the white cells on the left side of \eqref{eq:augmentation} are not in $D$ and exist in $[\lambda]$.) We call this process an \vocab{excited move}. Then, we define an \vocab{excited diagram} of \lm{} to be any set of $\abs{\mu}$ cells obtained by starting with the cells of $[\mu] \subseteq [\lambda]$ and applying any number of excited moves. We let $\cale(\lm)$ be the set of excited diagrams of \lm.

Given an excited diagram $D \in \cale(\lm)$ define $\varphi(D)$ to be the tableau $T$ of shape $\mu$ such that $T(u)$ is the row number of the cell that the original cell $u$ of $\mu$ gets excited to in $D$. This yields the following characterization of excited diagrams as flagged tableaux already known by Kreiman \cite[Sec. 6]{Kre1} and Knuton--Miller--Yong \cite[Sec. 5]{KMY} in the context of Schubert calculus and studied in \cite[Sec. 3]{mpp1} in the context of \eqref{eq:nhlf}.

\begin{theorem}[\cite{Kre1,KMY,mpp1}] \label{thm:excitedflag}
Let $\ssf_i$ be the maximum row number to which last cell in row $i$ in $[\mu]$ can be excited. Then the function $\varphi$ is a bijection between $\cale(\lm)$ and the set $\SSYT(\mu,(\ssf_1,\ldots,\ssf_{\ell(\mu)}))$. 
\end{theorem}

\subsection{The triangular lattice and lozenge tilings}
Throughout this paper we will often work in the \vocab{triangular lattice}, composed from tiling the plane with equilateral triangles of side length $1$, oriented so that each triangle has a horizontal edge, either on the top or bottom. To describe the non-horizontal directions in this lattice, we use the terms NE, NW, SE, and SW. A \vocab{lozenge} is a rhombus composed of two adjacent triangles, which can arise in three different orientations, which we will describe by the orientation of their longer diagonal, so that there are vertical lozenges, NE-SW lozenges, and NW-SE lozenges. We concede that in the above definitions a term such as ``NE'' can possibly refer to two different directions, but this ambiguity should be resolved by context.

Given a region $R$ in the triangular lattice, we let $L(R)$ be the set of tilings of $R$ by lozenges, or \vocab{lozenge tilings} for short. It is well-known that these are in bijection with families of non-intersecting paths on $R$:
\begin{theorem}[see e.g.\ \cite{adc2,propp,gorinlectures}] \label{thm:lozengepath}
Let $H(R)$ be the set of horizontal edges in $R$, and define a graph with vertices $H(R)$, connecting two horizontal edges if there exists a lozenge within $R$ with those two edges as opposite sides. Then, the NE-SW and NW-SE rhombi induce a system of vertex-disjoint paths on this graph, which begin at the horizontal edges which border the bottom of $R$ and end at the horizontal edges which border the top of $R$. Moreover, this correspondence from lozenge tilings to path systems is a bijection.
\end{theorem}

Crucially, similar results hold for the other two possible orientations. For an illustration, see Figures \ref{fig:oopathstiling}, \ref{fig:oopathsrgb}, and \ref{fig:oopathsbtiling}.

\subsection{Equivariant Knutson--Tao puzzles} \label{sec:KT-puzzles}
Given a partition $\lambda$ that fits inside an $a \times b$ rectangle, consider the path in the rectangle which delimits the boundary of $[\lambda]$. The following is an example for $\lambda= (6,4,4,1)$:
\begin{center}
\begin{tikzpicture}[scale=0.45]
\draw (0,0) grid (6,4);
\draw[ultra thick] (0,0)--(1,0)--(1,1)--(4,1)--(4,3)--(5,3)--(6,3)--(6,4);
\end{tikzpicture}
\end{center}

Define $w_\lambda$ to be the string with $a$ zeros and $b$ ones that represents the moves of this path from the bottom-left to the top-right, with a zero representing a vertical move that a one representing a horizontal move. Here, $w_\lambda = 1011100101$.

Given three partitions $\mu,\nu,\lambda$ that fit inside an $a \times b$ box, let an \vocab{(equivariant) Knutson--Tao puzzle} \cite{knutsontao2003puzzles} be the set of fillings of an equilateral triangle of side length $a+b$ with the following pieces (of side length $1$) where the northwest edge, northeast edge, and bottom edge read $w_\mu, w_\nu, w_\lambda$ from left to right, respectively:
\[\tikz[scale=0.32,yscale=1.732] \draw (0,0)--node {0}(2,0)--node {0}(1,1)--node {0}cycle; \qquad \tikz[scale=0.32,yscale=1.732] \draw (0,0)--node {1}(2,0)--node {1}(1,1)--node {1}cycle; \qquad \tikz[scale=0.32,yscale=1.732] \draw (0,0)--node {0}(1,-1)--node {1}(2,0)--node {0}(1,1)--node {1}cycle; \qquad \tikz[scale=0.32,yscale=1.732] \filldraw[fill=gray!50] (0,0)--node {1}(1,-1)--node {0}(2,0)--node {1}(1,1)--node {0}cycle;\]
Here only the first three pieces can be rotated and none can be reflected. As with lozenge tilings, the orientiations of the third piece are called vertical, NW-SE, and SW-NE. We call the last piece the \vocab{equivariant piece}. Call the set of all such puzzles $\Delta^{\mu\nu}_{\lambda}$.

Given a puzzle $P = \Delta^{\mu\nu}_{\lambda}$, define $\lozenge(P)$ to be the set of equivariant pieces in $P$. Given a piece $p$ in $\lozenge(P)$, let $\operatorname{ht}(p)$ be the height of $P$, i.e. the positive integer from $1$ to $a+b-1$ describing the distance from the center of $p$ to the bottom edge.

\subsection{Outside decompositions}
Let a \vocab{strip} be a connected skew shape with no $2 \times 2$ square. Given a strip $\theta$, note that there exist integers $m, n$ such that the contents of the cells of $\theta$ occupy every integer from $m$ to $n$ (inclusive) exactly once. Let $\theta[a, b]$ be the substrip of $\theta$ which starts with the cell of content $a$ and ends with the content of $b$. Let $\theta[a, b]$ be undefined if $b < a$.

Consider a skew shape \lm{} and let the minimum and maximum contents in $[\lm]$ be $m, n$. If $\theta$ is a strip with minimum and maximum contents $m, n$, then let the \vocab{outside decomposition} of \lm{} with \vocab{cutting strip} $\theta$ be the set of strips created by moving $\theta$ diagonally in a northwest-southeast direction and recording the connected components of the intersections of the shifted versions of $\theta$ and $\lm$ \cite{cyy2005}. Say this results in strips $(\theta_1, \theta_2, \ldots, \theta_k)$ and let the contents of $\theta_i$ range from $m_i$ to $n_i$. Then define $\theta_i \hash \theta_j = \theta[m_i, n_j]$.

The \vocab{Lascoux--Pragacz decomposition} of a skew shape \lm{} is the outside decomposition where the cutting strip is the outer border strip of $\lambda$.

\subsection{Euler and Genocchi numbers}
The $n$-th \vocab{Euler numbers} $E_{n}$ count the number of alternating permutations of size $n$ 
\cite[\href{http://oeis.org/A000111}{A000111}]{oeis}. The unsigned \vocab{Genocchi number} is defined as $G_{2n} = 2(-1)^{n-1}(4^n-1) B_{2n}$, where $B_{2n}$ is the $2n$th Bernoulli number. Note that $E_{2n-1} = G_{2n} \cdot 2^{2n-1}/2n$.

A \vocab{pistol} is a sequence of positive integers $a_1, a_2, \ldots, a_n$ so that $a_k \leq \frac{k+1}{2}$ for all $1 \leq k \leq n$. A pistol is \vocab{strictly alternating} if $a_k \geq a_{k+1}$ if $k$ is odd and $a_k < a_{k+1}$ if $k$ is even, for all $1 \leq k < n$. It is known that $G_{2n}$ counts the number of strictly alternating pistols of length $2n-1$ \cite{dumont80}. Also, define the \vocab{median Genocchi numbers} by letting $H_{2n+1}$ be the number of strictly alternating pistols of length $2n$ \cite{zeng06}. The first few terms for each sequence are $G_2 = G_4 = 1, G_6 = 3, G_8 = 17$ 
\cite[\href{http://oeis.org/A110501}{A110501}]{oeis} and $H_3 = 1, H_5 = 2, H_7 = 8, H_9 = 56$ \cite[\href{http://oeis.org/A005439}{A005439}]{oeis}.

\subsection{Superfactorials} \label{sec:superfactorials}
We will use the following notation: the \emph{superfactorial} $\Phi(n) = 1! \cdot 2! \cdots n!$ (the integer values of the Barnes $G$-function \cite{barnesG}), the \emph{double superfactorial} $\Psi(n) = 1! \cdot 3! \cdots (2n-1)!$, and the \emph{super doublefactorial} $\Lambda(n) = 1!! \cdot 3!! \cdots (2n-1)!!$ where $(2n-1)!!=(2n-1)\cdot (2n-3)\cdots 3\cdot 1$ is a \emph{double factorial}. Known asymptotics for these functions can be found in \cite[\S 2.10]{mpp3}.

\section{Basic Properties} \label{sec:background}
In \cite{okounkovolshanski98}, the Okounkov--Olshanski formula is proven using representation-theoretic methods. Here, we give a quick proof independent of these concepts and discuss positivity.

\subsection{A quick proof of the Okounkov--Olshanski formula}
\begin{definition}
If \lm{} is a skew shape, define $T_\lm \coloneqq \eval{s_\mu\mpar{\vy_\lambda}{\vy}}_{y_i = i}$.
\end{definition}

The original form of the Okounkov--Olshanski formula \cite[Thm.\ 8.1]{okounkovolshanski98} is the following.
\begin{theorem}[{Okounkov--Olshanski \cite[Thm. 8.1]{okounkovolshanski98}}] \label{thm:flmintermsofs}
\[f^\lm = \frac{f^\lambda \cdot \abs{\lm}!}{\abs{\lambda}!} T_\lm\]
\end{theorem}

\begin{proof}
Let $\ell_i = \lambda_i + d - i$ and $m_i = \mu_i + d - i$. By the determinantal form of the factorial Schur function \eqref{eq:factschurdet}, 
\begin{align*}
T_\lm &= s_\mu\mpar{\ell_1+1, \ell_2 + 1, \ldots}{1, 2, \ldots} \\
&= \frac{\prod_i \ell_i!}{\prod_{i<j} (\ell_i - \ell_j)} \cdot \det\brac*{\frac{1}{(\ell_j - m_i)!}}_{i,j=1}^d.
\end{align*}
By using the ordinary product formula for $f^\lambda$ \cite[Thm.\ 14.5.1]{adinroichman} and the Jacobi-Trudi formula for $f^\lm$ \cite[Thm.\ 14.5.6]{adinroichman}, this yields
\[
T_\lm = \frac{\abs{\lambda}!}{f^\lambda} \cdot \frac{f^\lm}{\abs{\lm}!}.
\]
Rearranging gives the desired result.
\end{proof}

\begin{proof}[Proof of Theorem \ref{thm:oof}]
By the hook length formula and Theorem \ref{thm:flmintermsofs}, we have
\[f^\lm = \frac{\abs{\lm}!}{\prod_{u \in [\lambda]} h(u)} T_\lm. \]
Since factorial Schur functions $s_\mu\mpar{\vx}{\va}$ are symmetric in $\vx$, we can evaluate
\begin{align*}
T_\lm &= s_\mu\mpar{\lambda_d + 1, \lambda_{d-1} + 2, \ldots, \lambda_1 + d}{1,2,\ldots} \\
&= \sum_{T \in \SSYT(\mu, d)}\prod_{u\in[\mu]} (\lambda_{d+1-T(u)} + T(u) - (T(u)+c(u))) \\
&= \sum_{T \in \SSYT(\mu, d)}\prod_{u\in[\mu]} (\lambda_{d+1-T(u)} - c(u)).
\end{align*}

Combining these two expressions yields the Okounkov--Olshanski formula.
\end{proof}

\begin{example} \label{ex:oof}
For the shape $\lm=2221/11$, there are six tableaux in $\SSYT(\lambda,4)$:
\begin{equation} \label{eq:oo-ex-ssyt}
\ytab[centertableaux]{3,4},\quad  \ytab[centertableaux]{2,4},\quad \ytab[centertableaux]{2,3},\quad \ytab[centertableaux]{1,4},\quad \ytab[centertableaux]{1,3},\quad \ytab[centertableaux]{1,2},
\end{equation}
and the Okounkov--Olshanski formula \eqref{eq:oof} gives 
\begin{equation} \label{eq:ex-oof} 
f^{\lm} \,=\, \frac{5!}{2\cdot 3\cdot 3 \cdot 4 \cdot 5}\left(2\cdot 3 + 2\cdot 3 + 2\cdot 3 + 1\cdot 3 + 1\cdot 3 + 1 \cdot 3 \right) \,=\, 9.
\end{equation}
\end{example}

\subsection{Nonnegativity of the formula}
A quick look at the formula \eqref{eq:oof} suggests that there could be negative terms. However, we next show that every term in the formula is nonnegative.

\begin{proposition} \label{prop:nonnegativity}
Every term in the Okounkov--Olshanski formula is nonnegative. Moreover, every positive term has all weights positive.
\end{proposition}
\begin{proof}
Suppose that there exist $i,j$ so that $\lambda_{d+1-T(i,j)} - c(i,j) < 0$. It suffices to show that there exists some $j'$ so that $\lambda_{d+1-T(i,j')} - c(i,j') = 0$.

Since $c(i,j) > 0$ then $j > i$. For $i \leq k \leq j$, consider the quantity $a_k = \lambda_{d+1-T(i,k)} - c(i,k)$. Note that $a_{k+1} - a_k = -1 + \lambda_{d+1-T(i,k+1)} - \lambda_{d+1-T(i,k)} \geq -1$. Thus, since $a_i = \lambda_{d+1-T(i,i)} \geq 0$ and $a_j<0$, there must be some $k$ with $a_k = 0$, as desired.
\end{proof}

For most of the remainder of the paper, we will focus on the nonzero terms in \eqref{eq:oof}. These are indexed by the following objects:

\begin{definition}
For a skew shape $\lm$ of length $d$, an {\em Okounkov--Olshanski tableau} is a SSYT $T$ in $\SSYT(\mu, d)$ such that $c(u) < \lambda_{d + 1 - T(u)}$ for all $u \in [\mu]$. Let $\cOOT(\lm)$ be the set of such tableaux and let $\OOT(\lm)$ be the size of $\cOOT(\lm)$.
\end{definition}

\begin{corollary}
The nonzero terms in the Okounkov--Olshanski formula for the shape $\lm$ correspond to SSYT $T$ in $\cOOT(\lm)$.
\end{corollary}

\begin{proof}
This characterization follows from  Proposition \ref{prop:nonnegativity}.   
\end{proof}

\section{Correspondences on Okounkov--Olshanski Tableaux} \label{sec:posterms}
In this section we will demonstrate bijections which operate on the set $\cOOT(\lm)$, converting its members into more useful forms for enumeration or other analysis. Central to this idea are two well-known concepts: Theorem \ref{thm:lozengepath}, which relate lattice paths on a shape to lozenge tilings, and a correspondence between semistandard tableaux to lattice paths:

\begin{theorem}[\cite{gv89}, see also {\cite[\S 4.5]{saganbook}}] \label{thm:ssytlattice}
The set $\SSYT(\lm, m)$ is in bijection with systems of vertex-disjoint paths in $\setz^2$ which travel rightwards and upwards from $(\mu_i - i, 1)$ to $(\lambda_i - i, m)$, for $1 \leq i \leq m$, such that the $i$th path has rightward steps $(-i+j-1, T(i,j)) \to (-i+j, T(i,j))$ for all $\mu_i < j \leq \lambda_i$.
\end{theorem}

In what follows, we will first leverage these tools to convert nonzero Okounkov--Olshanski terms, which correspond to certain special semistandard tableaux, into lozenge tilings of a certain shape. Then, we will apply a similar bijection backwards to yield a different class of semistandard tableaux. We will then extend this to excited diagrams and reverse excited diagrams. Figure~\ref{fig:schematic} contains a ``commutative diagram'' of all the bijections.

\subsection{From Okounkov--Olshanski tableaux  to lozenge tilings} \label{subsec:OOtoLT}
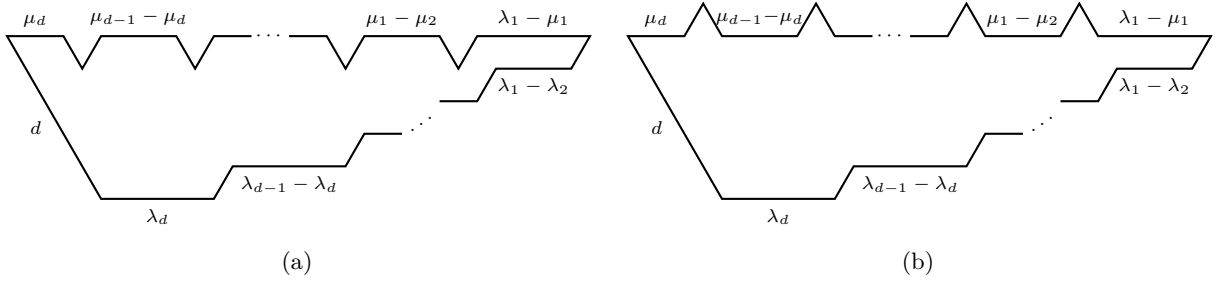
\begin{figure}[tb]
\centering
\begin{subfigure}{0.5\linewidth} \centering
\begin{tikzpicture}[x=0.5cm,y=0.5cm,yscale=0.866,xscale=0.5,thick,font=\scriptsize]
\draw (8,5)--++(-2,0)--++(-1,-1)--++(-1,1)--node[anchor=south]{$\mu_{d-1}-\mu_d$}++(-4,0)--++(-1,-1)--++(-1,1)--node[anchor=south]{$\mu_d$}
(-5,5)--node[shift=(30:-8pt)]{$d$}
(0,0)--node[anchor=north]{$\lambda_d$}++(6,0)--++(1,1)--node[anchor=north]{$\lambda_{d-1}-\lambda_d$}++(6,0)--++(1,1)--++(2,0);
\node at (17,2.5) {$\iddots$};
\node at (9,5) {$\cdots$};
\draw (10,5)--++(2,0)--++(1,-1)--++(1,1)--node[anchor=south]{$\mu_1-\mu_2$}++(4,0)--++(1,-1)--++(1,1)--node[anchor=south]{$\lambda_1-\mu_1$}++(6,0)--++(-1,-1)--node[anchor=north]{$\lambda_1 - \lambda_2$}++(-4,0)--++(-1,-1)--++(-2,0);
\end{tikzpicture}
\caption{}
\label{fig:oolozenge1}
\end{subfigure}%
\begin{subfigure}{0.5\linewidth} \centering
\begin{tikzpicture}[x=0.5cm,y=0.5cm,yscale=0.866,xscale=0.5,thick,font=\scriptsize]
\draw (8,5)--++(-2,0)--++(-1,1)--++(-1,-1)--node[anchor=south]{$\mu_{d-1}{-}\mu_d$}++(-4,0)--++(-1,1)--++(-1,-1)--node[anchor=south]{$\mu_d$}
(-5,5)--node[shift=(30:-8pt)]{$d$}
(0,0)--node[anchor=north]{$\lambda_d$}++(6,0)--++(1,1)--node[anchor=north]{$\lambda_{d-1}-\lambda_d$}++(6,0)--++(1,1)--++(2,0);
\node at (17,2.5) {$\iddots$};
\node at (9,5) {$\cdots$};
\draw (10,5)--++(2,0)--++(1,1)--++(1,-1)--node[anchor=south]{$\mu_1-\mu_2$}++(4,0)--++(1,1)--++(1,-1)--node[anchor=south]{$\lambda_1-\mu_1$}++(6,0)--++(-1,-1)--node[anchor=north]{$\lambda_1 - \lambda_2$}++(-4,0)--++(-1,-1)--++(-2,0);
\end{tikzpicture}
\caption{}
\label{fig:oolozenge2}
\end{subfigure}
\caption{The regions $\nabla^*_\lm$ and $\nabla_\lm$ both correspond to nonzero Okounkov-Olshanski terms. Figure \ref{fig:oolozenge1} is derived from the lattice paths corresponding to the SSYT that appear in \eqref{eq:oof}. Figure \ref{fig:oolozenge2} is a slight modification which is easier to work with.}
\label{fig:oolozenge}
\end{figure}
In this section we are concerned with the following two shapes:
\begin{definition}
Given a skew shape \lm, we let $\nabla^*_\lm$ be the region shown in Figure \ref{fig:oolozenge1}. Similarly, let $\nabla_\lm$ be the region shown in Figure \ref{fig:oolozenge2}.
\end{definition}

\begin{remark} \label{prop:ooterms3}
Observe that $L(\nabla^*_\lm)$ and $L(\nabla_\lm)$ are ``canonically isomorphic'' since each tiling in $L(\nabla^*_\lm)$ can be augmented with vertical lozenges in the triangular gaps along the top boundary to yield a tiling in $L(\nabla_\lm)$, and vice versa. For the remainder of this paper, we will assume this correspondence and frequently interchange the two sets.
\end{remark}

\begin{figure}[tb]
\centering
\begin{subfigure}[b]{0.3\linewidth} \centering
\begin{tikzpicture}[x=0.45cm,y=0.45cm]
\begin{scope}[gray, thin]
\foreach \x in {-6,...,0} \draw(\x,1)--(\x,5);
\draw (1,2)--(1,5);
\draw (2,3)--(2,5);
\draw (3,4)--(3,5);
\draw(-6.5,1)--(0,1);
\draw(-6.5,2)--(1,2);
\draw(-6.5,3)--(2,3);
\draw(-6.5,4)--(3,4);
\draw(-6.5,5)--(4,5);

\end{scope}
\filldraw (-1,1) circle[radius=2pt] 
(-2,1) circle[radius=2pt]
(-3,1) circle[radius=2pt]
(-4,1) circle[radius=2pt]
(-5,1) circle[radius=2pt]
(2,5) circle[radius=2pt]
(0,5) circle[radius=2pt]
(-2,5) circle[radius=2pt]
(-4,5) circle[radius=2pt]
(-5,5) circle[radius=2pt];
\draw[thick] (-2,2) circle[radius=2pt]
(-3,3) circle[radius=2pt]
(-4,4) circle[radius=2pt];
\draw[very thick] (-1,1)--(0,1)--(0,2)--(1,2)--(1,3)--(2,3)--(2,5);
\draw[very thick] (-2,1)--(-2,2)--(-1,2)--(-1,5)--(0,5);
\draw[very thick] (-3,1)--(-3,5)--(-2,5);
\draw[very thick] (-4,1)--(-4,5);
\draw[very thick] (-5,1)--(-5,5);
\end{tikzpicture}
\caption{}
\label{fig:oopathstrunc}
\end{subfigure}%
\begin{subfigure}[b]{0.4\linewidth} \centering
\begin{tikzpicture}[x=0.45cm,y=0.45cm,yscale=0.866,xscale=0.5]
\begin{scope}[fill=red!25,semithick]
\filldraw (-1,1)--++(2,0)--++(1,1)--++(2,0)--++(1,1)--++(2,0)--++(2,2)--++(1,-1)--++(-2,-2)--++(-2,0)--++(-1,-1)--++(-2,0)--++(-1,-1)--++(-2,0)--cycle;
\filldraw (-2,2)--++(2,0)--++(3,3)--++(2,0)--++(1,-1)--++(-2,0)--++(-3,-3)--++(-2,0)--cycle;
\filldraw (-2,2)--++(2,2)--++(2,0)--++(-1,1)--++(-2,0)--++(-2,-2)--cycle;
\filldraw (-4,4)--(-3,5)--(-2,4)--(-3,3)--cycle;
\end{scope}
\begin{scope}[very thin]
\draw(2,0)--(0,2) (3,1)--(1,3)--(-1,3)--(-2,4) (5,1)--(4,2) (6,2)--(5,3)--(3,3)--(2,4)--(0,4)--(-1,5) (8,2)--(7,3) (9,3)--(8,4)--(6,4)--(5,3) (4,4)--(3,5) (10,4)--(12,4)--(13,5);
\draw [gray] (-9,1)--(-1,1)--(-2,0)--(-6,4)(-8,2)--(-2,2)--(-4,0)--(-7,3)(-7,3)--(-3,3)--(-6,0)--(-8,2)(-6,4)--(-4,4)--(-8,0)--(-9,1);
\end{scope}
\draw[thick] (0,0)--++(2,0)--++(1,1)--++(2,0)--++(1,1)--++(2,0)--++(1,1)--++(2,0)--++(1,1)--++(2,0)--++(1,1)--++(-4,0)--++(-1,-1)--++(-1,1)--++(-2,0)--++(-1,-1)--++(-1,1)--++(-2,0)--++(-1,-1)--++(-1,1)--++(-2,0)--++(-1,-1)--++(-1,1)--++(-1,-1)--++(-1,1)--cycle;
\draw[thick, black!20] (0,0)--(-10,0)--(-5,5);
\begin{scope}[y={(1,1)},x={(2,0)},shift={(0.5,-0.5)}]
\draw[red!25, very thick] (-2,1)--(-2,2) (-3,1)--(-3,3) (-4,1)--(-4,4) (-5,1)--(-5,5);
\begin{scope}[every node/.style={fill=red,circle,inner sep=0pt,minimum size=4pt}]
\begin{scope}[every node/.append style={fill=red!25}]
\node at (-2,1) {};
\node at (-3,1) {};
\node at (-4,1) {};
\node at (-5,1) {};
\end{scope}
\node at (-1,1) {};
\node at (-2,2) {};
\node at (-3,3) {};
\node at (-4,4) {};
\node at (2,5) {};
\node at (0,5) {};
\node at (-2,5) {};
\node at (-4,5) {};
\node at (-5,5) {};
\end{scope}
\draw[red, very thick] (-1,1)--(0,1)--(0,2)--(1,2)--(1,3)--(2,3)--(2,5)
 (-2,2)--(-2,2)--(-1,2)--(-1,5)--(0,5)
 (-3,3)--(-3,5)--(-2,5)
 (-4,4)--(-4,5);
\end{scope}
\end{tikzpicture}
\caption{}
\label{fig:oopathstiling}
\end{subfigure}%
\begin{subfigure}[b]{0.3\linewidth} \centering
\begin{tikzpicture}[x=0.45cm,y=0.45cm,yscale=0.866,xscale=0.5]
\draw (-1,1)--++(2,0)--++(1,1)--++(2,0)--++(1,1)--++(2,0)--++(2,2)--++(1,-1)--++(-2,-2)--++(-2,0)--++(-1,-1)--++(-2,0)--++(-1,-1)--++(-2,0)--cycle;
\draw (-2,2)--++(2,0)--++(3,3)--++(2,0)--++(1,-1)--++(-2,0)--++(-3,-3)--++(-2,0)--cycle;
\draw (-2,2)--++(2,2)--++(2,0)--++(-1,1)--++(-2,0)--++(-2,-2)--cycle;
\draw (-4,4)--(-3,5)--(-2,4)--(-3,3)--cycle;
\draw(2,0)--(0,2) (3,1)--(1,3)--(-1,3)--(-2,4) (5,1)--(4,2) (6,2)--(5,3)--(3,3)--(2,4)--(0,4)--(-1,5) (8,2)--(7,3) (9,3)--(8,4)--(6,4)--(5,3) (4,4)--(3,5) (6,4)--(7,5);
\draw (11,5)--(10,4)--(12,4)--(13,5);
\begin{scope}[ultra thick]
\draw[red] (-0.5,0.5)--++(2,0)--++(1,1)--++(2,0)--++(1,1)--++(2,0)--++(3,3);
\draw[red] (-1.5,1.5)--++(2,0)--++(3,3)--++(2,0)--++(1,1);
\draw[red] (-2.5,2.5)--++(2,2)--++(2,0)--++(1,1);
\draw[red] (-3.5,3.5)--++(2,2);
\draw[red] (-4.5,4.5)--++(1,1);

\draw[blue] (-4.5,5.5)--++(3,-3)--++(2,0)--++(2,-2);
\draw[blue] (-2.5,5.5)--++(2,-2)--++(2,0)--++(1,-1)--++(2,0)--++(1,-1);
\draw[blue] (1.5,5.5)--++(2,-2)--++(4,0)--++(1,-1);
\draw[blue] (5.5,5.5)--++(1,-1)--++(2,0)--++(1,-1)--++(2,0);
\draw[blue] (9.5,5.5)--++(1,-1)--++(4,0);

\draw[green!80!black] (1,0)--++(-2,2)--++(2,2)--++(-1,1);
\draw[green!80!black] (4,1)--++(-1,1)--++(2,2)--++(-1,1);
\draw[green!80!black] (7,2)--++(-1,1)--++(2,2);
\draw[green!80!black] (10,3)--++(2,2);
\draw[green!80!black] (13,4)--++(1,1);
\end{scope}
\draw[thick] (0,0)--++(2,0)--++(1,1)--++(2,0)--++(1,1)--++(2,0)--++(1,1)--++(2,0)--++(1,1)--++(2,0)--++(1,1)--++(-4,0)--++(-1,1)--++(-1,-1)--++(-2,0)--++(-1,1)--++(-1,-1)--++(-2,0)--++(-1,1)--++(-1,-1)--++(-2,0)--++(-1,1)--++(-1,-1)--++(-1,1)--++(-1,-1)--cycle;
\end{tikzpicture}
\caption{}
\label{fig:oopathsrgb}
\end{subfigure}\\[0.2cm]%
\begin{subfigure}[b]{0.3\linewidth} \centering
\begin{tikzpicture}[x=0.45cm,y=0.45cm,yscale=0.866,xscale=0.5]
\begin{scope}[fill=blue!15,semithick]
\filldraw (2,0)--++(1,1)--++(-2,2)--++(-2,0)--++(-2,2)--++(-1,-1)--++(2,-2)--++(2,0)--cycle;
\filldraw (5,1)--++(1,1)--++(-1,1)--++(-2,0)--++(-1,1)--++(-2,0)--++(-1,1)--++(-1,-1)--++(1,-1)--++(2,0)--++(1,-1)--++(2,0)--++(1,-1)--cycle;
\filldraw (8,2)--++(1,1)--++(-1,1)--++(-4,0)--++(-1,1)--++(-1,-1)--++(1,-1)--++(4,0)--cycle;
\filldraw (11,3)--++(1,1)--++(-2,0)--++(-1,1)--++(-2,0)--++(-1,-1)--++(2,0)--++(1,-1)--cycle;
\filldraw (14,4)--++(1,1)--++(-4,0)--++(-1,-1)--cycle;
\end{scope}
\begin{scope}[very thin]
\draw (-1,1)--(1,1)--(4,4) (-3,3)--(-2,4) (-2,2)--(0,4) (0,2)--(2,4) (4,2)--(6,4) (7,3)--(9,5) (9,3)--(10,4) (12,4)--(13,5);
\end{scope}
\draw[thick] (0,0)--++(2,0)--++(1,1)--++(2,0)--++(1,1)--++(2,0)--++(1,1)--++(2,0)--++(1,1)--++(2,0)--++(1,1)--++(-4,0)--++(-1,-1)--++(-1,1)--++(-2,0)--++(-1,-1)--++(-1,1)--++(-2,0)--++(-1,-1)--++(-1,1)--++(-2,0)--++(-1,-1)--++(-1,1)--++(-1,-1)--++(-1,1)--cycle;
\begin{scope}[y={(-1,1)},x={(2,0)},shift={(0.5,-0.5)}]
\begin{scope}[every node/.style={fill=blue,circle,inner sep=0pt,minimum size=4pt}]
\node at (1,1) {};
\node at (3,2) {};
\node at (5,3) {};
\node at (7,4) {};
\node at (9,5) {};
\node at (0,5) {};
\node at (1,5) {};
\node at (3,5) {};
\node at (5,5) {};
\node at (7,5) {};
\end{scope}
\draw[blue, very thick] (0,5)--(0,3)--(1,3)--(1,1)
(1,5)--(1,4)--(2,4)--(2,3)--(3,3)--(3,2)
(3,5)--(3,4)--(5,4)--(5,3)
(5,5)--(6,5)--(6,4)--(7,4)
(7,5)--(9,5);
\end{scope}
\end{tikzpicture}
\caption{}
\label{fig:oopathsbtiling}
\end{subfigure}%
\begin{subfigure}[b]{0.4\linewidth} \centering
\begin{tikzpicture}[x=0.45cm,y=0.45cm]
\draw[gray, thin] (-6.5,0.5) grid[step=0.45cm] (5.5,5);
\filldraw (2,5) circle[radius=2pt] 
(0,5) circle[radius=2pt]
(-2,5) circle[radius=2pt]
(-4,5) circle[radius=2pt]
(-5,5) circle[radius=2pt];
\draw [thick]
(4,5) circle[radius=2pt]
(2,4) circle[radius=2pt]
(0,3) circle[radius=2pt]
(-2,2) circle[radius=2pt]
(-4,1) circle[radius=2pt];
\draw[very thick] (-5,5)--(-5,3)--(-4,3)--(-4,1) (-4,5)--(-4,4)--(-3,4)--(-3,3)--(-2,3)--(-2,2) (-2,5)--(-2,4)--(0,4)--(0,3) (0,5)--(1,5)--(1,4)--(2,4)--(2,3) (2,5)--(4,5)--(4,4);
\draw[very thick,gray,->] ((-4,1)--(-4,0.5);
\draw[very thick,gray,->] (-2,2)--(-2,0.5);
\draw[very thick,gray,->] (0,3)--(0,0.5);
\draw[very thick,gray,->] (2,4)--(2,0.5);
\draw[very thick,gray,->] (4,5)--(4,0.5);
\end{tikzpicture}
\caption{}
\label{fig:oopathspreflag}
\end{subfigure}%
\begin{subfigure}[b]{0.3\linewidth} \centering
\begin{tikzpicture}[x=0.45cm,y=0.45cm]
\draw[red,thick](-2,5)--(-3,6);
\draw[red,thick](0,5)--(-1,6);
\draw[red,thick](1,2)--(0,3);
\draw[red,thick](2,3)--(1,4);
\draw[red,thick](0,1)--(-1,2);
\draw[red,thick](-1,2)--(-2,3);
\fill[cyan!70] (4,5) rectangle (2,6) rectangle (3,4) (0,6) rectangle (1,4) rectangle (-1,5) rectangle (0,3) (-3,5) rectangle (-2,3);
\begin{scope}[very thin]
\draw(-5,5)--(-5,6);
\draw(-4,4)--(-4,6);
\draw(-3,3)--(-3,6);
\draw(-2,2)--(-2,6);
\draw(-1,1)--(-1,6);
\draw(0,1)--(0,6);
\draw (1,2)--(1,6);
\draw (2,3)--(2,6);
\draw (3,4)--(3,6);
\draw(-1,1)--(0,1);
\draw(-2,2)--(1,2);
\draw(-3,3)--(2,3);
\draw(-4,4)--(3,4);
\draw(-5,5)--(4,5);
\draw(-5,6)--(4,6);
\end{scope}
\draw[very thick] (-1,1)--++(1,0)--++(0,1)--++(1,0)--++(0,1)--++(1,0)--++(0,1)--++(1,0)--++(0,1)--++(1,0)--++(0,1)
--++(-2,0)--++(0,-1)--++(-1,0)--++(0,-1)--++(-1,0)--++(0,-1)--++(-1,0)--cycle;
\end{tikzpicture}
\caption{}
\label{fig:oorevexcite}
\end{subfigure}\\[0.2cm]%
\begin{subfigure}[b]{0.3\linewidth} \centering
\ytableaushort{
123,25,5}
\caption{}
\label{fig:oosorigtab}
\end{subfigure}%
\begin{subfigure}[b]{0.3\linewidth} \centering
\globaldefs=-1
\ytableausetup{smalltableaux}
\globaldefs=0
\ytableaushort{
\none\none\none11,\none\none12,\none22,23,3}
\caption{}
\label{fig:ooskewtab}
\end{subfigure}%
\begin{subfigure}[b]{0.4\linewidth} \centering
\begin{tikzpicture}[x=0.45cm,y=0.45cm]
\fill[cyan!70] (0,5) rectangle (1,3) (2,0) rectangle (3,1) (4,0) rectangle (5,1) (2,3) rectangle (3,4) (4,2) rectangle (5,3);
\draw[step=1.0, very thin] (0,5) grid (9.5, 0);
\draw [very thick] (0,5)--++(0,-3)--++(1,0)--++(0,1)--++(1,0)--++(0,1)--++(1,0)--++(0,1)--cycle;
\draw [very thick, gray, dashed] (0,5)--++(1,0)--++(0,-1)--++(2,0)--++(0,-1)--++(2,0)--++(0,-1)--++(2,0)--++(0,-1)--++(2,0)--++(0,-1)--++(-9,0)--cycle;
\end{tikzpicture}
\caption{}
\label{fig:ooe}
\end{subfigure}%
\caption{Constructions related to the bijection described in Section \ref{sec:posterms}.}
\label{fig:oopaths}
\end{figure}
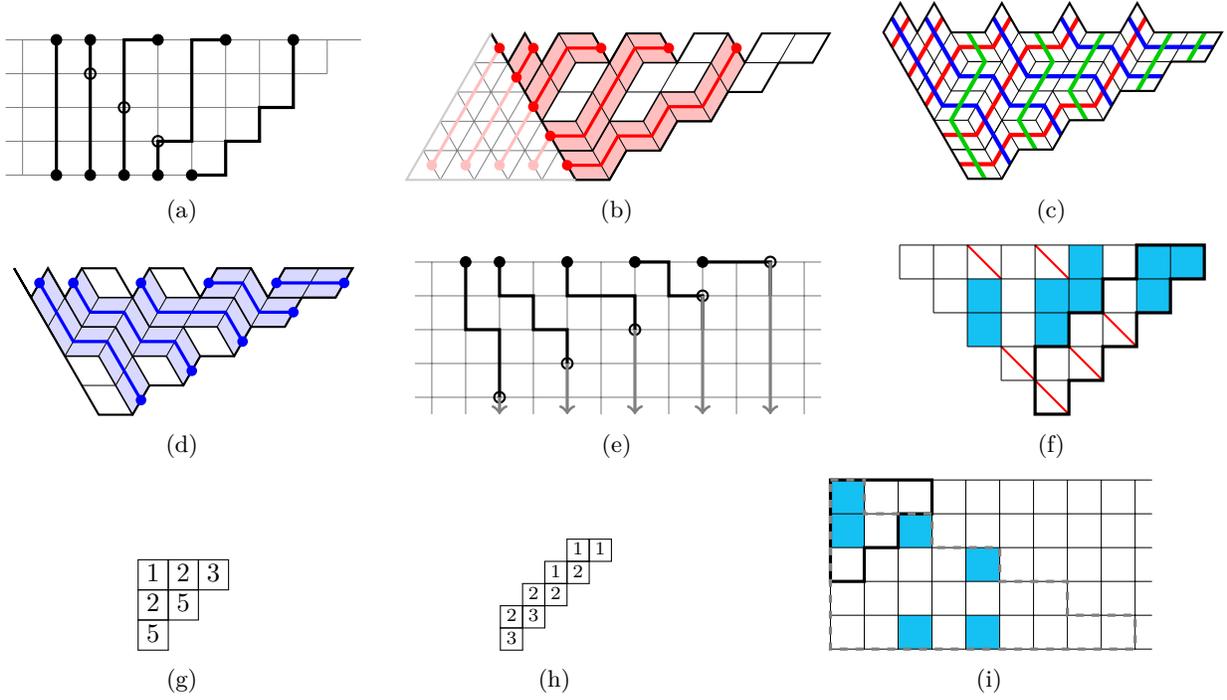
Now, we detail a correspondence between $\cOOT(\lm)$ and $L(\nabla_\lm)$ using Theorems \ref{thm:ssytlattice} and \ref{thm:lozengepath}.

\begin{definition}
Given a skew shape $\lm$, let $\LP(\lm)$ be the collection of systems of non-intersecting lattice paths $(\rho_1,\ldots,\rho_d)$ such that 
\begin{itemize}
    \item[(i)] The $i$th path $\rho_i$ goes from $(-i,i)$ to $(\mu_i-i, d)$ for $1 \leq i \leq d$.
\item[(ii)] The paths only move upwards and rightwards.
\item[(iii)] The support of the  paths is contained in the set $\set{(x,y)\mid x \leq \lambda_{d + 1 - y}}$.
\end{itemize}
\end{definition}

\begin{proposition} \label{prop:ooterms1}
The SSYT in $\cOOT(\lm)$ are in bijection with systems of non-intersecting lattice paths in $\LP(\lm)$.
\end{proposition}

\begin{proof}
This is straightforward upon applying the bijection in Theorem \ref{thm:ssytlattice} to the tableaux in $\cOOT(\lm)$.
\end{proof}

\begin{proposition} \label{prop:ooterms2}
Non-intersecting systems of lattice paths in $\LP(\lm)$ are in bijection with lozenge tilings in $L(\nabla^*_\lm)$.
\end{proposition}
\begin{proof}
By Theorem \ref{thm:lozengepath}, $L(\nabla_\lm^*)$ is in bijection with NE-right paths within $\nabla^*_\lm$ that begin at the NW-SE edges on the left side of $\nabla_\lm^*$ and end at the the NW-SE edges on the top boundary of the shape.

It is straightforward to show that shearing path systems in $\LP(\lm)$ so that an upwards move becomes a NE move makes it so that the endpoints of the paths can be placed on the boundary NW-SE edges. Moreover, the bound $x \leq \lambda_{d + 1- y}$ transforms to the right boundary of $\nabla^*_\lm$, meaning that this shear is indeed a bijection between path systems. This concludes the proof.
\end{proof}

\begin{definition}
Let $\R: L(\nabla^*_\lm) \to \cOOT(\lm)$ be the map $\mathsf{t}\mapsto T$ obtained by composing the inverses of the bijections in the proofs of Propositions~\ref{prop:ooterms1} and \ref{prop:ooterms2}, respectively.
\end{definition}

\begin{lemma} \label{lem:SSYT2lozengetilings}
For a skew shape $\lm$, the map $\R$ is a bijection between  $L(\nabla^*_\lm) \cong L(\nabla_{\lm})$ and $\cOOT(\lm)$. 
\end{lemma}

\begin{proof}
This follows by combining Propositions~\ref{prop:ooterms1} and  \ref{prop:ooterms2} (see Remark~\ref{prop:ooterms3}).
\end{proof}

\begin{example} \label{exmp:ooterms}
Consider the  skew shape $\lm=$ where $\lambda = 54321$, $\mu = 32100$, and the following SSYT of shape $\mu$ in $\cOOT(\lm)$:
\begin{equation} \label{ex:runningexSSYT}
T = \ytab[centertableaux]{123,25,5}.
\end{equation}
The associated objects that are in correspondence with this SSYT can be seen in Figure \ref{fig:oopaths}. In particular, the bijection detailed in Proposition~\ref{prop:ooterms1} yields Figure~\ref{fig:oopathstrunc}. Applying the bijections described in Propositions~\ref{prop:ooterms2} and Remark~\ref{prop:ooterms3} to this example yield Figures~\ref{fig:oopathstiling} and \ref{fig:oopathsrgb}, respectively.
\end{example}

\subsection{From lozenge tilings to flagged tableaux}\label{subsec:LTtoFT}
In this subsection we apply the previous techniques to the lozenge tiling, but in a different orientation. This will yield a  flagged tableau.

\begin{definition}
Let $\cSF(\lm)=\SSYT(\lm,(1,2,\ldots,\ell(\lambda))$ be the set of SSYT $T$ of shape \lm\ where $T(i,j) \leq i$ for $(i,j)\in [\lm]$.
\end{definition}

\begin{remark}
The flagged skew tableaux in $\cSF(\lm)$ appeared in the context of \emph{dual stable Grothendieck polynomials} in \cite{LP}, see Section~\ref{sec:dual stable}. 
\end{remark}

\begin{lemma} \label{prop:flaggedssyt}
For a skew shape $\lm$, there is a bijection $\B$ between  $L(\nabla^*_\lm)$ and $\cSF(\lm)$.
\end{lemma}

\begin{proof}
The desired bijection $\B$ between $L(\nabla^*_\lm)$ and $\cSF(\lm)$ is obtained as follows:

By Theorem~\ref{thm:lozengepath}, there is a bijection from tilings in $L(\nabla^*_\lm)$ to non-intersecting right-SE path systems connecting the $d$ SW-NE edges on the top of $\nabla^*_\lm$ to the $d$ SW-NE edges along the bottom boundary. Furthermore, notice that the condition that the paths stay inside $\nabla^*_\lm$ is unnecessary, since the only obstruction to the paths staying inside is of them colliding with the triangular indentations on the top boundary. However, the nonexistence of these collisions is already implied by the path system being non-intersecting.

Establishing a coordinate system in which $(1,0)$ is horizontal while $(0,1)$ points in the NW  direction, while the SW corner is $(1/2,1/2)$, then the $i$th path $\beta_i$ from the left goes from $(\mu_{d+1-i} + i, d)$ to $(\lambda_{d+1-i} + i, i)$. Therefore we have established a bijection from the right-SE path systems  to down-right lattice path systems  in $\setz^2$ from $(\mu_{d+1-i} + i, d)$ to $(\lambda_{d+1-i} + i, i)$.

By reversing the indices $i \leftrightarrow d+1-i$ and applying a translation, this is equivalent to the lattice paths from $(\mu_i - i, d)$ to $(\lambda_i - i, d+1-i)$. Furthermore, applying a vertical reflection makes the lattice paths go in the up-right direction from $(\mu_i - i, 1)$ to $(\lambda_i - i, i)$. These correspond precisely to $\cSF(\lm)$ under the canonical bijection from semistandard tableaux to lattice paths.
\end{proof}

\begin{example} \label{exmp:ooterms-flagged tableaux}
Continuing with Example~\ref{exmp:ooterms}, the bijection described in Lemma~\ref{prop:flaggedssyt} uses the blue paths in Figure~\ref{fig:oopathsrgb}, denoted in Figure~\ref{fig:oopathsbtiling}, to obtain the skew flagged tableaux in  Figure~\ref{fig:ooskewtab}.
\end{example}

At this point, we have now investigated two of the three path systems that can be obtained from a lozenge tiling of $\nabla^*_\lm$ (the red and blue paths in Figure~\ref{fig:oopathsrgb}). The third is created by connecting the midpoints of opposite horizontal edges of lozenges, which travel from the $\lambda_1$ horizontal edges on the bottom to the $\lambda_1$ horizontal edges on the top (the green path in Figure~\ref{fig:oopathsrgb}) and we denote the system by $(\gamma'_1,\ldots,\gamma'_d)$. Defining the \textit{depth} of a segment on one of these paths to be a variable that ranges from $1$ at the top to $d$ at the bottom, we can now state the following description of these paths.

\begin{proposition} \label{prop:greenpaths}
Let $\tiling$ be a tiling in $L(\nabla^*_\lm)$ with corresponding path system $(\gamma'_1,\ldots,\gamma'_d)$ such that $T\coloneqq\R(\tiling) \in  \cOOT(\lm)$ and $U\coloneqq\B(\tiling) \in \cSF(\lm)$. Then each $\gamma'_i$:
\begin{itemize}
    \item[(i)] has segments with depths from $1$ to $\lambda'_i$,
    \item[(ii)] has a SE step at depth $d+1-k$  if and only if $T$ has a $k$ entry in the $i$th column, 
    \item[(iii)] has a SW-NE step at depth $k$ if and only if $U$ has an entry of $k$ in the $i$th column.
\end{itemize}
\end{proposition}
\begin{proof}
Note that the lozenge in $\tiling$ directly above the $i$th horizontal edge on the bottom boundary has depth $\lambda'_i$, so the depths of the $i$th path range from $1$ to $\lambda'_i$.

Given a system $(\rho'_1,\ldots,\rho'_d)$ of right-NE paths on $\nabla^*_\lm$, one can construct NE-NW paths $(\gamma'_1,\ldots,\gamma'_d)$ by starting from the bottom edge, and iteratively moving to the NW when there is a right-NE path immediately above the current location, and moving to the NE otherwise. By translating this to $\setz^2$, the set of horizontal segments that the $i$th path $\gamma_i'$ intersects can be found by first taking the lowest edge that goes from $x = i-2$ to $x = i-1$, and then iteratively finding the lowest edge above the current one in the column one unit to the left, until no others exist. Since the $x$-coordinate of a horizontal segment is described by the content, the entries these horizontal edges correspond to first taking the smallest entry of content $i - 1$, and then iteratively taking the smallest larger entry with a content that is one lower. This is easily seen to simply be the $i$th column of $T$.

Since an entry of $k$ in $T$ corresponds to a NW-SE lozenge of depth $d + 1 - k$, then the NW segments in the $i$th path $\gamma_i$, which each cross a right-NE path, occur at depths $d + 1 - k$, where $k$ ranges over all the entries in column $i$ of $T$.

The relation to $U$ is similar. Here it is more useful to view the paths $\gamma_i'$ in reverse, taking an SW step whenever there is a right-SE path to be crossed and an SE step otherwise. Then, since the (reverse) $i$th path $\gamma'_i$ starts $d - \mu'_i + i - 1$ units from the left, in the perspective with lattice paths from $(\mu_i - i, 1)$ to $(\lambda_i - i, i)$, the SW edges correspond to taking the lowest edge that goes from $x = -\mu'_i + i - 2$ to $x = -\mu'_i + i - 1$ and then iteratively taking the lowest higher edge one column to the left. These corresponds to precisely the $i$th column of $U$. Therefore, since a SW edge at depth $k$ corresponds to an entry of $k$, we are done.
\end{proof}

As a corollary of the above result, we can give a direct description of the correspondence between $\cOOT(\lm)$ and $\cSF(\lm)$.

\begin{definition}
Given $T \in \cOOT(\lm)$, let $\OOTtoFT(T)\coloneqq U$ be the semistandard tableaux of shape $\lm$ defined as follows: for each $i=1,\ldots,\lambda_1$ let $U$ be the tableau of shape $\lm$ whose entries in the $i$th column are the  integers $k$ in $\{1,\ldots,\lambda'_i\}$ such that $d+1-k$ is not an entry in the $i$th column of $T$.
\end{definition}

\begin{corollary} \label{cor:directbij}
The map $\OOTtoFT(T)$ is a bijection between $\cOOT(\lm)$ and $\cSF(\lm)$.
\end{corollary}

\begin{proof}
For each column $i$, by Proposition \ref{prop:greenpaths}~(iii) we have that $k$ is in the $i$th column of $U$ precisely when the $i$th path $\gamma'_i$ has a SW step at depth $k$. Therefore we can use Proposition \ref{prop:greenpaths}~(ii) again to finish. 
\end{proof}

\begin{example} \label{exmp:oot-to-flagged tableaux}
Continuing with Example~\ref{exmp:ooterms}, the bijection $\OOTtoFT$ maps the tableaux in Figure~\ref{fig:oosorigtab} to Figure~\ref{fig:ooskewtab}.
\end{example}

\begin{remark}
See Section~\ref{sec: symmetries lozenge tilings} for more details on tableaux that can be obtained from the path systems of lozenge tilings in $L(\nabla^*_\lm)$.
\end{remark}

\subsection{Reverse excited diagrams} \label{subsec: reverse excited diagrams}
An alternate formulation of the flagged tableaux comes in the form of reverse excited diagrams.

\begin{definition}
Given a skew shape $\lm$, its {\em reverse excited diagrams} are the diagrams obtained by starting with the cells of the shifted skew shape $\lmstar$ and applying {\em reverse excited moves}:
\begin{equation}
\includegraphics[scale=0.8]{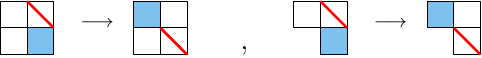} \label{eq:rev_excited_move} \tag{$\nwarrow$}
\end{equation}
(we ignore momentarily the red diagonals on certain cells). We let $\calr\cale(\lm)$ denote the set of reverse excited diagrams of $\lm$.
\end{definition}

\begin{example}
The skew shape $\lm=2221/11$ has two excited diagrams and six reverse excited diagrams as illustrated in Figure~\ref{fig:ex_reverse_excited}.
\begin{figure}
\centering
\begin{subfigure}{0.25\linewidth}
\includegraphics[scale=0.6]{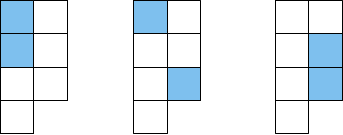}
\caption{}
\label{fig:ex_excited}
\end{subfigure}
\begin{subfigure}{0.7\linewidth}
\includegraphics[scale=0.6]{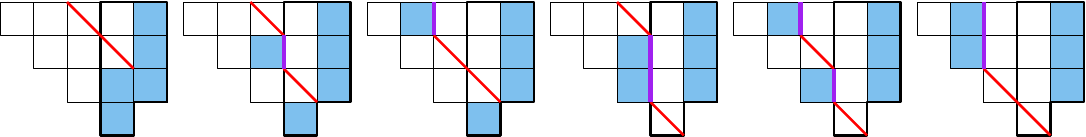}
\caption{}
\label{fig:ex_reverse_excited}
\end{subfigure}
\begin{subfigure}{0.25\linewidth}
\phantom{\includegraphics[scale=0.6]{resources/ex_excited.pdf}}
\end{subfigure}
\begin{subfigure}{0.7\linewidth}
\includegraphics[scale=0.6]{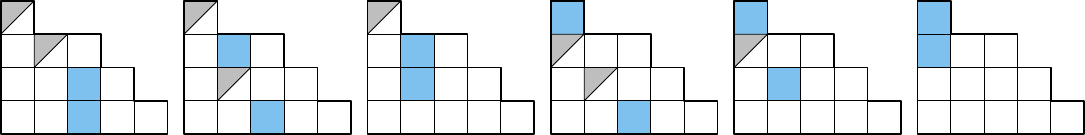}
\caption{}
\label{fig:ex_oo_excited}
\end{subfigure}

\caption{Excited diagrams, reverse excited diagrams, and Okounkov--Olshanski excited diagrams of the shape $\lm = 2221/11$.}
\label{fig:ex_exicted_okounkov_olshanski}
\end{figure}
\end{example}

Just as excited diagrams in $\cale(\lm)$ are in correspondence with flagged tableaux of shape $\mu$ (see Theorem~\ref{thm:excitedflag}), reverse excited in $\cRE(\lm)$ are in correspondence with flagged tableaux of shape $\lm$.  
\begin{definition}
Given a reverse excited diagram $D \in \cRE(\lm)$ define $\REEtoFT(D)$ to be the tableau of shape $\lm$ such that $T(u)$ is the row number of the cell that the original cell $u$ of $[\lm]$ gets excited to in $D$.
\end{definition}

\begin{lemma} \label{lem:redexcite}
The map $\REEtoFT$ is a bijection between $\cRE(\lm)$ and $\cSF(\lm)$.
\end{lemma}

\begin{proof}
Consider three partitions $\lambda, \mu, \nu$ with $[\nu]\subseteq[\mu]\subseteq[\lambda]$. Letting $\cale_\nu(\lambda/\mu)$ be the subset of $\cale(\lm)$ so that no cells in $\nu$ are excited, applying Theorem \ref{thm:excitedflag} yields that $\cale_\nu(\lambda/\mu)$ is in bijection with semistandard tableaux of shape $\mu$ with entries bounded by $\ssf_i$ so that every cell in $\nu$ is filled with its row number. These are equivalent to semistandard tableaux of shape $\mu/\nu$ such that the entries in row $i$ are at least $i$ but at most $\ssf_i$.

Suppose $\lambda$ and $\mu$ have $d$ parts and fix a large integer $N$ at least $\lambda_1$. Let $\lambda_r = (N - \lambda_d, N - \lambda_{d-1}, \ldots, N-\lambda_1)$ and define $\mu_r$ similarly. Then, since the allowed moves in a reverse excited diagram are exactly those for an excited diagram, but rotated 180 degrees, $\cRE(\lambda/\mu)$ corresponds to elements in $\cale_{\lambda_r}(\nu/\mu_r)$, where $\nu = (N + d - 1)^d$. The aforementioned correspondence means that these excited diagrams biject to semistandard tableaux of shape $\mu_r/\lambda_r$ so that the entries in row $i$ are at least $i$ but at most $d$. Rotating back and replacing each entry $r$ with $d + 1 - r$ yields tableaux of shape $\lm$ so that the entries in row $i$ are at most $i$, exactly equal to $\cSF(\lm)$. This tableaux contains the row data of the original reverse excited diagram since if a cell is excited to row $i$, the rotated version has the corresponding cell excited to row $d + 1 - i$.
\end{proof}

\begin{example} \label{exmp:ree-to-flagged tableaux}
Continuing with Example~\ref{exmp:ooterms}, the bijection $\REEtoFT$ maps reverse excited diagram in Figure~\ref{fig:oorevexcite} to the flagged tableau in Figure~\ref{fig:ooskewtab}.
\end{example}

To continue, we define broken diagonals associated with each reverse excited diagram:
\begin{figure}[tbp]
    \centering
    \includegraphics{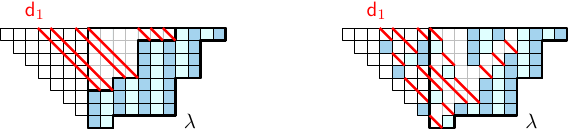}
    \caption{Broken diagonals of the reverse excited diagrams $[\lmstar]$ and $D$.}
    \label{fig:reverse_diagonals}
\end{figure}
\begin{definition}[Broken diagonals] \label{def:brokendiag}
For each reverse excited diagram $D$ in $\calr\cale(\lm)$ we define \vocab{broken diagonals} $\mathsf{d}_1,\ldots,\mathsf{d}_{\ell(\mu')}$, subsets of the complement of $D$, as follows:
\begin{itemize}
\item For the initial diagram $[\lmstar]$, let $\mathsf{d}_i$ be the cells in $[\mu^*]$ with contents $i-\mu'_i+d-1$. 
\item If $D'$ is obtained from $D$ by applying the reverse excited move $(i,j)\to (i-1,j-1)$, then define $\mathsf{d}_r(D')=\mathsf{d}_r(D)$ for $r \neq t$ and $\mathsf{d}_t(D') = (\mathsf{d}_t(D) \setminus (i-1,j)) \cup \set{(i,j)}$, where $t$ is the unique index such that $(i-1,j) \in \mathsf{d}_t(D)$.
\end{itemize}
\end{definition}

\begin{definition}
Define $B(D) = \bigcup_i \mathsf{d}_i(D)$, i.e.\ the cells of all the broken diagonals of $D$. 
\end{definition}

\begin{remark}
For a visualization of how broken diagonals change with an excited move, see the red diagonals in \eqref{eq:rev_excited_move}. For an example see Figure~\ref{fig:ex_reverse_excited} and Figure~\ref{fig:reverse_diagonals}.
\end{remark}

\begin{remark}
The notion of broken diagonals of regular excited diagrams in $\cale(\lm)$ already appeared in \cite[\S 7]{mpp3}. 
\end{remark}

Definition \ref{def:brokendiag} may appear to have ambiguities; in particular, it is not clear that the broken diagonals of a reverse excited diagram are well-defined, as opposed to depending on the order of the excited moves. Moreover, the existence and uniqueness of $t$ doesn't follow immediately either from the local description of excited diagrams. To remedy these issues, one can construct a global characterization of excited diagonals by considering the locations of the cells excited from each column in $\lm$.

\begin{definition}
For an excited diagram $D$ in $\cRE(\lm)$, consider the cells excited from column $j$ in $\lm$ originally located at $(\mu'_j + 1, j + d - 1), \ldots, (\lambda'_j, j + d - 1)$. Suppose that cell $(i, j+d-1)$ is excited $a_i$ times to end up at $(i-a_i, j+d-a_i-1)$. By Lemma \ref{lem:redexcite}, $a_i$ is a strictly decreasing sequence in $i$, so for $\mu'_j < i < \lambda'_j$ there is a natural set of $a_i-a_{i+1}$ diagonal cells denoted by $\mathsf{c}_{i, j}(D)$ from $(i-a_i+1, j+d-a_i)$ from $(i-a_{i+1}, j+d-a_{i+1}-1)$, connecting the $i$th and $(i+1)$th excited cells. Extending this pattern, there are also diagonals, which we denote $\mathsf{c}_{\lambda_j', j}(D)$ and $\mathsf{c}_{\mu_j', j}(D)$, from $(\lambda_j' - a_{\lambda_j'} + 1, j + d - a_{\lambda_j'})$ to $(\lambda_j', j + d - 1)$ and from $(1, j + d - \mu'_j)$ to $(\mu'_j - a_{\mu'_j + 1}, j + d - a_{\mu'_j + 1} - 1)$, connecting the cells at each end of the column to the limits of the shape; these consist of $a_{\lambda_j'}$ and $\mu_j' - a_{\mu_j'+1}$ cells, respectively. Let $\mathsf{C}_j(D) = \bigcup_i \mathsf{c}_{i,j}(D)$.
\end{definition}

\begin{proposition} \label{prop:diagglobal}
For a reverse excited diagram $D$ in $\cRE(\lm)$  we have that $\mathsf{C}_j(D) = \mathsf{d}_j(D)$. Also, $\mathsf{d}_j(D)$ are all disjoint from each other and from $D$.
\end{proposition}
\begin{proof}
We will use induction on the number of excited moves of $D$.

If $D = [\lmstar]$, then all the $a_i$ are zero, so $\mathsf{C}_j(D) = \mathsf{c}_{\mu'_j,j}(D)$, which extends from from $(1, j + d - \mu'_j)$ to $(\mu'_j, j + d - 1)$. All these cells have contents $j - \mu_j' + d - 1$, matching $\mathsf{d}_j(D)$. It is also apparent that extending this diagonal further would cause it to leave $[\mu^*]$, since the bottom edge of the cell $(\mu'_j, j + d - 1)$ is on the boundary of $[\mu^*]$. Thus $\mathsf{C}_j(D) = \mathsf{d}_j(D)$. Also, since $j - \mu_j' + d - 1$ is strictly increasing in $j$, all the $\mathsf{d}_j(D)$ are disjoint; since they are additionally all contained in $[\mu^*]$, they are also disjoint from $D$.

Now suppose $D'$ is obtained from $D$ by applying the move $(i,j+d-1) \to (i-1,j+d-2)$. If the original location of this cell was $(i+k,j+k+d-1)$, then $\mathsf{c}_{i+k-1,j+k}(D') = \mathsf{c}_{i+k-1,j+k}(D) \setminus \set{(i-1,j+d-1)}$ and $\mathsf{c}_{i+k,j+k}(D) = \mathsf{c}_{i+k,j+k}(D) \cup \set{(i,j+d-1)}$. Now, since $(i-1,j+d-1) \in \mathsf{c}_{i+k-1,j+k}(D) \subseteq \mathsf{d}_{j+k}(D)$ and all the $\mathsf{d}_i(D)$ are disjoint, $t = j+k$, proving that $\mathsf{C}_r(D') = \mathsf{d}_r(D')$ for all $r$.

To finish, we need to prove that all the $\mathsf{d}_r(D')$ are disjoint from $D'$ and each other. First, the only cell added to any of the $\mathsf{d}_r(D')$ was $(i, j+d-1)$ added to $\mathsf{d}_{j+k}(D')$. Since this was in $D$, it cannot be in any of the $\mathsf{d}_r(D') = \mathsf{d}_r(D)$ for any $r \neq j+k$, meaning that all the $\mathsf{d}_r(D')$ are disjoint. To prove that none of them overlap with $D'$, it suffices to show that $(i, j+d-1) \notin D'$ and $(i-1, j+d-2) \notin B(D)$. The first is obvious by the definition of an excited move. To prove the second statement, observe that if $(i-1, j+d-2)$ were part of a broken diagonal, it would be the case that $(i, j+d-2) \in D$, $(i, j+d-1) \in B(D)$, or $(i, j+d-2) \notin [\lambda^*]$. Each of these is clearly impossible.
\end{proof}

This alternate definition of the $\mathsf{d}_j(D)$ yields an important connection with lozenge tilings.

\begin{lemma} \label{lem:redaffine}
Consider a lozenge tiling $\tiling$ in $L(\nabla_{\lm})$. Upon applying a shear preserving the horizontal direction and setting SW-NE lozenges to squares, the tiling can be naturally overlaid on $[\lambda^*]$. Under this association, the SW-NE lozenges of $\tiling$ are exactly the elements of the corresponding reverse excited diagram $D$ and the right edges of the NW-SE lozenges lie exactly on the cells of $B(D)$.
\end{lemma}
\begin{proof}
Under this shear, the paths of horizontal edges in Proposition \ref{prop:greenpaths} become paths of parallelograms, which direct the path in a NW-SE direction, and squares, which direct that path vertically. Observe that the $i$th path from the left starts at the bottom of the $i$th column.

Consider all cells such that their bottom edge is in the $i$th path; the bottommost of these cells has content $i - \lambda'_i + d - 1$. Since moving through a parallelogram doesn't change the content, while moving upwards through a square increases the content by $1$, the content of the $k$th square from the bottom in that path has content $i + k - \lambda'_i + d - 1$.

By Proposition \ref{prop:greenpaths}, there are $\lambda'_i - \mu'_i$ squares, occurring in the row numbers equal to the entries of the $i$th column in the corresponding element of $\cSF(\lm)$. Since the contents match those of the $i$th column of $[\lmstar]$, the squares in the $i$th path are exactly the elements of $D$ excited from the $i$th column.

By Proposition \ref{prop:diagglobal}, the elements of $\mathsf{d}_i(D)$ come in ``strings'' connecting the cells excited from the $i$th column of $[\lmstar]$. A close examination of the placement of these strings yields that they are located at the right edges of the parallelograms in the $i$th path. Since every SW-NE and NW-SE lozenge lies in a path of horizontal edges, we are done.
\end{proof}

\begin{example} \label{exmp:oot-to-reverse excited}
The shear described in Lemma~\ref{lem:redaffine} maps the tiling in Figure~\ref{fig:oopathsrgb} to the reverse excited diagram in  Figure~\ref{fig:oorevexcite}.
\end{example}

\subsection{Loose ends and summary} \label{subsec:ooesum}
If a partition $\tilde\lambda$ with $d$ parts satisfies $\tilde \lambda_d \geq \mu_1 + d - 1$, then one can check that $\mathsf{f}_i = d$. Therefore, by Theorem \ref{thm:excitedflag}, $\SSYT(\mu, d)$ is in bijection with $\cale(\tilde\lambda/\mu)$ for all such shapes $\tilde \lambda$. Since the image of the excited diagram doesn't depend on the exact value of $\tilde\lambda$, we will instead denote this set by $\cale(\Lambda^d/\mu)$ where $\Lambda$ denotes an unbounded part. Observe that $\cOOT(\lm)$ is in bijection with a subset of $\cale(\Lambda^d/\mu)$.

\begin{definition} \label{def:ooe}
Let $\cOOE(\lm) \subseteq \cale(\Lambda^d/\mu)$ be the image of the map $\varphi$ defined in Theorem \ref{thm:excitedflag} on $\cOOT(\lm)$.
\end{definition}

Analogously to Lemma \ref{lem:redaffine}, we have the following result. 

\begin{proposition} \label{prop:ooetolozenge}
Given a lozenge tiling $\tiling$ in $L(\nabla_{\lm})$, apply a shear preserving the horizontal direction sending NW-SE lozenges to squares. Then, applying a vertical flip and aligning the left edge with the left edge of $\Lambda^d$, the squares are precisely those in the associated diagram in $\cOOE(\lm)$.
\end{proposition}

\begin{example} \label{exmp:oot-to-ooexcited}
The shear described in Lemma~\ref{lem:redaffine} maps the tiling in Figure~\ref{fig:oopathsrgb} to the  excited diagram in  Figure~\ref{fig:ooe}.
\end{example}

\begin{remark} \label{rmk: description of OOE}
From this characterization, it follows that the diagrams in $\cOOE(\lm)$ can be described as diagrams obtained by starting with the cells of $[\mu]$ and applying excited moves within an upside-down $[\lambda^*]$, where $\lambda^*$ is interpreted as a regular partition. See Figure~\ref{fig:ex_oo_excited}. 
\end{remark}

\begin{figure}[tbp] \centering
\begin{tikzpicture}
\node (lozenge) {$L(\nabla_\lm) \simeq L(\nabla_\lm^*)$};
\node (green) [below=of lozenge] {\color{green!50!black} path systems};
\node (red) [left=of green] {\color{red} path systems};
\node (oot) [left=of red] {$\cOOT(\lm)$};
\node (ooe) [above=of oot] {$\cOOE(\lm)$};
\node (blue) [right=of green] {\color{blue} path systems};
\node (flag) [right=of blue]{$\cSF(\lm)$};
\node (rev) [above=of flag] {$\cRE(\lm)$};
\draw[<->]
(lozenge) to node[auto] {Prop\ \ref{prop:greenpaths}} (green);
\draw[<->](lozenge) to  node[auto, swap] {Prop.\ \ref{prop:ooterms2}} (red);
\draw[<->] (lozenge) to node[auto] {Lem.\ \ref{prop:flaggedssyt}} (blue);
\draw[<->] (oot) to node[auto] {Prop.\ \ref{prop:ooterms1}} (red);
\draw[<->, bend right=15, looseness=0.5] (oot) to node[auto, swap] {Prop. \ref{prop:greenpaths}} (green);
\draw[<->] (flag) to  node[auto, swap] {Lem.\ \ref{prop:flaggedssyt}}(blue);
\draw[<->, bend left=15, looseness=0.5] (flag)  to node[auto] {Prop.\ \ref{prop:greenpaths}} (green);
\draw[<->] (rev) to  node[auto,swap] {Lem.\ \ref{lem:redexcite}} (flag);
\draw[<->] (rev) to node[auto, swap] {Lem.\ \ref{lem:redaffine}} (lozenge);
\draw[<->, bend right=25, looseness=0.8] (oot) to node[auto, swap] {Cor.\ \ref{cor:directbij}} (flag);
\draw[<->] (ooe) to node[auto] {Prop.\ \ref{prop:ooetolozenge}} (lozenge);
\draw[<->] (ooe) to node[auto] {Def.\ \ref{def:ooe}} (oot);
\end{tikzpicture}
\caption{Bijections constructed in Section \ref{sec:posterms}, labeled by their respective result. Path systems are colored based on their color in Figure \ref{fig:oopathsrgb}.} \label{fig:schematic}
\end{figure}
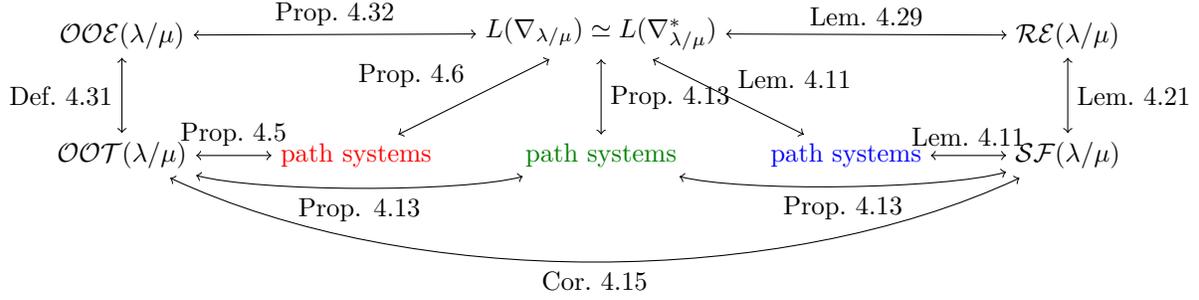

To summarize, we have found that a nonzero Okounkov--Olshanski term corresponds to two families of Young tableaux, $\cOOT(\lm)$, and $\cSF(\lm)$. Each of these is associated with four types of ``graphical data'': two types of path systems regarding the rows and columns, an excited diagram, and a lozenge tilings. The correspondence between these families is derived from such graphical data: the column path systems of each are similar, which translates to identical lozenge tilings. A visual representation of the relations we have found is shown in Figure \ref{fig:schematic}.

The astute reader may notice that these bijections come in two main families: using path systems or excited diagrams. To establish all the related objects, we worked in the context of path systems, which broke the combinatorial objects into rows or columns and operated on them individually. Later, through analyzing these operations, we were able to construct global bijections that did not rely on this decomposition, passing through the language of excited diagrams. Therefore, while in a sense path systems were integral to defining these bijections, excited diagrams make it easier to work with them.

\section{Reformulations of the Okounkov--Olshanski Formula} \label{sec:reformulations}
\subsection{Weight transformation rules} \label{subsec:weighttrans}
In this subsection we apply the bijections described in Section \ref{sec:posterms} to reformulate the Okounkov--Olshanski formula. We will operate along the top of Figure \ref{fig:schematic}, excluding path systems.

In order to phrase these reformulations, we first make a few definitions:
\begin{definition}
Given $\tiling \in L(\nabla_\lm)$, let $R_\leftarrow(\tiling)$ be the set of NW-SE rhombi in $\tiling$ and define $w \colon R_\leftarrow(\tiling) \to \setz$ so that $w(r)$ is the distance from the center of $r$ to the rightmost point in $\nabla_\lm$ at the same height.
\end{definition}

\begin{definition}
Given $T \in \cSF(\lm)$, let $M(T)$ be the set of $(j, t) \in \setz^2$, with $1 \leq j \leq \ell(\lambda')$ and $1 \leq t \leq \lambda'_j$, so that there is no $i$ with $T(i,j) = t$. Given $m \in M(T)$ with $m = (j,t)$, let $i(m)$ be the unique integer $\mu'_j \leq i \leq \lambda'_j$ so that $T(i',j) < t$ for all $i' \leq i$ and $T(i', j) > t$ for all $i' > i$.
\end{definition}
It is often useful to visualize $M(T)$ as a set labels on the horizontal edges of $T$; in particular, representing $(j, t) \in M(T)$ as a label of $t$ between the cells $(i(j,t), j)$ and $(i(j,t)+1, j)$ makes it so that column $j$ contains every single number from $1$ to $\lambda'_j$, in increasing order from top to bottom:
\begin{center}
\[\vcenter{\hbox{\begin{tikzpicture}[x=0.45cm,y=0.45cm,scale=0.5,font=\footnotesize]
\draw (0, 0)-- (0, 4)-- (2, 4)-- (2, 6)-- (4, 6)-- (4, 8)-- (6, 8)-- (6, 10)-- (10, 10)-- (10, 8)-- (8, 8)-- (8, 6)-- (6, 6)-- (6, 4)-- (4, 4)-- (4, 2)-- (2, 2)-- (2, 0)-- (0, 0)-- (0, 0);
\draw
 (0, 2)-- (2, 2)-- (2, 4)-- (4, 4)-- (4, 6)-- (6, 6)-- (6, 8)-- (8, 8)-- (8, 10);
\node at (9, 9) {1};
\node at (7, 9) {1};
\node at (7, 7) {2};
\node at (1, 3) {2};
\node at (5, 7) {1};
\node at (3, 5) {2};
\node at (3, 3) {3};
\node at (1, 1) {3};
\node at (5, 5) {2};
\phantom{\node at (1, 0) {$4, 5$};}
\end{tikzpicture}}}
\quad
\longleftrightarrow
\quad
\vcenter{\hbox{\begin{tikzpicture}[x=0.45cm,y=0.45cm,scale=0.5,font=\footnotesize]
\draw (0, 0)-- (0, 4)-- (2, 4)-- (2, 6)-- (4, 6)-- (4, 8)-- (6, 8)-- (6, 10)-- (10, 10)-- (10, 8)-- (8, 8)-- (8, 6)-- (6, 6)-- (6, 4)-- (4, 4)-- (4, 2)-- (2, 2)-- (2, 0)-- (0, 0)-- (0, 0);
\draw
 (0, 2)-- (2, 2)-- (2, 4)-- (4, 4)-- (4, 6)-- (6, 6)-- (6, 8)-- (8, 8)-- (8, 10);
\node at (9, 9) {1};
\node at (7, 9) {1};
\node at (7, 7) {2};
\node at (5, 4) {3};
\node at (1, 3) {2};
\node at (5, 7) {1};
\node at (3, 5) {2};
\node at (3, 3) {3};
\node at (1, 1) {3};
\node at (5, 5) {2};
\node at (3, 6) {1};
\node at (1, 4) {1};
\node at (3, 2) {4};
\node at (1, 0) {$4, 5$};
\end{tikzpicture}}}\]
\end{center}
These tableaux with labeled edges have appeared before in the literature, as the equivariant semistandard tableaux of Thomas and Yong \cite{thomasyong2012}. We defer further discussion of the connection to Section \ref{subsec:eqssyt}.

Now fix $T \in \cOOT(\lm)$ and say the defined maps
\[\cOOT(\lm) \to \cOOE(\lm) \to L(\nabla_{\lm}) \to \cRE(\lm) \to \cSF(\lm)\]
map a tableau $T \in \cOOT(\lm)$ as $T \mapsto D \mapsto \tiling \mapsto D' \mapsto T'$.

These maps induce natural bijections on the constituents of these objects, mapping
\begin{equation} \label{eq: mapping constituent objects}
[\mu] \to D \to R_\leftarrow(\Delta) \to B(D') \to M(T').
\end{equation}
These are defined by
\begin{itemize}
    \item mapping a cell $(i,j) \in [\mu]$ to the cell in $D$ it is excited to, as discussed in Theorem \ref{thm:excitedflag};
    \item mapping a cell in $D$ to the NW-SE rhombus obtained after flipping $D$ upside down and applying a shear to align its left boundary with the left boundary of $\nabla_\lm$, as discussed in Proposition \ref{prop:ooetolozenge};
    \item mapping a NW-SE rhombus to the cell centered on its right edge after applying a shear sending $\nabla_\lm$ to $[\lambda^*]$, as described in Lemma \ref{lem:redaffine};
    \item sending an element $(i,j) \in \mathsf{d}_k(D')$ to $(k, i)$, which is a bijection as $\mathsf{d}_k(D')$ contains a cell in row $i$ if and only if $i$ is not an entry in column $k$ of $T'$, by Proposition \ref{prop:diagglobal}.
\end{itemize}
These observations can be encoded into the following proposition:
\begin{proposition} \label{prop:wherecellsgo}
The bijections in \eqref{eq: mapping constituent objects} send a cell $(i,j) \in [\mu]$ with $T(i,j) = t$, to
\begin{itemize}
    \item $(t, j-i+t) \in D$
    \item a NW-SE rhombus $r$ in the $t$'th row from the bottom with $w(r) = \lambda_{d+1-t} + i - j$
    \item $(d+1-t, j-i + d) \in B(D')$
    \item $(j, d+1-t) \in M(T)$
\end{itemize}
Moreover, $i(j, d+1-t) = d+i-t$.
\end{proposition}
\begin{proof}
Most of this is a straightforward application of the definitions; the only wrinkle is proving that $k$ in the definition of the map $B(D') \to M(T')$ is equal to $j$, i.e.\ that $(d+1-t, \lambda_{d+1-t} + i - j) \in \mathsf{d}_j(D')$.

To prove this, observe that, in proving Proposition \ref{prop:ooetolozenge}, we show that the shear sending $D$ to $R_\leftarrow(\tiling)$ sends the cells excited from the $j$th column in $T$ to the NW-SE rhombi on the $j$th path on horizontal edges from the left. Therefore, $r$ lies on this $j$th path. Now, the proof of Lemma \ref{lem:redaffine} shows that the image of the NW-SE rhombi on the $j$th path under the shear $R_\leftarrow(\tiling) \to B(D')$ is precisely $\mathsf{d}_j(D')$, as desired.

Finally, observe that every cell in $\mathsf{c}_{i,j}(D')$ has content $j + d - 1 - i$, so $(d+1-t, j-i+d) \in \mathsf{c}_{d+i-t, j}(D')$. Therefore, $T'(d+i-t,j) < d+1-t < T'(d+i-t+1, j)$, so $i(j, d+1-t) = d+i-t$, as desired.
\end{proof}

Due to these correspondences, it is immediate to establish the following reformulations of \eqref{eq:oof}.

\begin{corollary}[Excited diagram formulation] \label{cor:edform}
\[f^{\lm} = \frac{\abs{\lambda/\mu}!}{\prod_{u \in [\lambda]} h(u)} \sum_{D \in \cOOE(\lm)} \prod_{(i,j) \in D} (\lambda_{d+1-i} + i - j).\]
\end{corollary}
\begin{corollary}[Lozenge tiling formulation] \label{cor:lozengeform}
\[f^{\lm} = \frac{\abs{\lambda/\mu}!}{\prod_{u \in [\lambda]} h(u)} \sum_{\tiling \in L(\nabla_\lm)} \prod_{r \in R_\leftarrow(\tiling)} w(r).\]
\end{corollary}
\begin{corollary}[Reverse excited diagram formulation] \label{cor:redform}
\[f^{\lm} = \frac{\abs{\lambda/\mu}!}{\prod_{u \in [\lambda]} h(u)} \sum_{D \in \cRE(\lm)} \prod_{(i,j) \in B(D)} (\lambda_i + d - j).\]
\end{corollary}
\begin{corollary}[Flagged tableaux formulation] \label{cor:flag-form}
\[f^{\lm} = \frac{\abs{\lambda/\mu}!}{\prod_{u \in [\lambda]} h(u)} \sum_{T \in \cSF(\lm)} \prod_{(j, t) \in M(T)} (\lambda_t-t+ i(j,t)-j+1).\]
\end{corollary}

\begin{remark} \label{rmk: explanation of weights of reformulations}
We give interpretations to some of the weights of some of the reformulations above:
\begin{itemize}
\item[(i)] With the description of $\cOOE(\lm)$ in Remark~\ref{rmk: description of OOE}, the quantity $\lambda_{d+1-i}+i-j$ in the excited diagram formulation of Corollary~\ref{cor:edform} is the arm-length of the cell $(i,j)$ in $D\in \cOOE(\lm)$. See Figure~\ref{fig:ex_oo_excited}.

\item[(ii)] Since the outer edge of $[\lambda^*]$ occurs at $(i, \lambda_i + d - 1)$, the quantity $\lambda_i + d - j$ in the reverse excited diagram formulation of Corollary \ref{cor:redform} is the arm-length of the cell $(i, j)$ in $D\in \cRE(\lm)$.  See Figure~\ref{fig:ex_reverse_excited}.

\item[(iii)] The quantity $\lambda_t -t + i(j,t)-j+1$ in the flagged tableaux formulation of Corollary~\ref{cor:flag-form} is the hook length of the cell $(t,j)$ to $(t,\lambda_t)$ and $(i(j,t),j)$ for $(j,t)$ in $M(T)$.
\end{itemize}

\end{remark}

\subsection{Equivalence to equivariant Knutson--Tao puzzles}
Starting from Corollary \ref{cor:lozengeform}, we can now prove Theorem \ref{thm:ooiskt}. The proof will require the following result from \cite{knutsontao2003puzzles}.

\begin{lemma}[{\cite[Corollary 3]{knutsontao2003puzzles}}] \label{lem:ktlemma}
If a Knutson--Tao puzzle has boundaries equal to $\lambda, \mu, \lambda$, it does not contain any SW-NE rhombi.
\end{lemma}
\input{resources/figurektcompressed.tex}
\begin{proof}[Proof of Theorem \ref{thm:ooiskt}]
We will describe a bijection between $L(\nabla_\lm)$ and $\Delta^{\lambda\mu}_{\lambda}$ that preserves weights.

Given an element of $L(\nabla_\lm)$, split each vertical rhombus into two equilateral triangles and remove the equilateral triangles sticking out of the top part of the region. Then label all horizontal edges $1$, with the exception of the horizontal edges used to split the vertical lozenges. Label all other edges $0$.

Now split the tiling into horizontal \vocab{bands} of height one and rotate each of the bands clockwise by 60 degrees. Fit the bands into an equilateral triangle of the appropriate size such that the ends of the bands are the edges of the triangle. Finally, between the shapes, add vertical non-equivariant rhombi and all-$1$ triangles. Observe that this process turns vertical rhombi to 0-triangles, NW-SE rhombi to equivariant pieces, and NE-SW rhombi to NW-SE non-equivariant rhombi. Thus the result is an equivariant Knutson--Tao puzzle (see Figure~\ref{fig:kt2}).

Due to the definition of the top boundary of $\nabla_\lm$, which had peaks where $w_\mu$ had zeroes, the northeast edge is labeled with $w_\mu$. Since the $i$th row from the top in $\nabla_\lm$ was decomposed into $2(\lambda_i + d - i) + 1$ triangles, and the ends of these rows were labeled $0$, the northwest and bottom edges have zeroes in the $(\lambda_i + d - i + 1)$th positions, exactly the locations of the zeroes in $w_\lambda$. Finally, since the bands composed of non-equivalent rhombi and 1-triangles are bounded by edges labeled 1, the northwest and bottom edges have $w_\lambda$. Therefore we have established a mapping from $L(\nabla_\lm)$ to $\Delta^{\lambda\mu}_{\lambda}$.

Now we will show that this correspondence is a bijection by constructing its inverse. Consider an element of $\Delta^{\lambda\mu}_\lambda$. By Lemma \ref{lem:ktlemma}, we can split up the triangle into bands that run northwest to southeast.

Consider the bands that have a horizontal or NE-SW edge $e$ labeled $1$. Since the only allowed pieces with such an edge are 1-triangles and non-equivariant vertical rhombi, any piece bordering $e$ must be one of these two pieces, each of which induces a second 1-edge further down the band. Thus, by an easy induction, all edges strictly inside the band are labeled $1$, and the band must be composed of all-$1$ triangles or non-equivariant vertical rhombi.

Therefore, except for a 1-triangle at the bottom of the band, the band is composed of vertical non-equivariant rhombi and pairs of 1-triangles attached along a horizontal edge; in particular, since the edges at the ends of the band are labeled $1$, they must border ones on the boundary. Also, the southeast edge of the band has the same labels as the northwest edge, except that the northwest edge has an extra $1$ at the bottom. This allows us to contract all the bands with edges labeled $1$ to join all the bands bordering zeroes along their northwest edges. Rotate everything 60 degrees counterclockwise to get a tiling of $\nabla_{\lm}$, minus the triangular points at the top.

Since there are no 1-triangles or non-equivariant vertical rhombi left, the tiling must be composed of 0-triangles or rhombi, either NW-SE or NE-SW, with horizontal edges labeled $1$ and diagonal edges labeled $0$. If $w_\lambda$ has $k$ trailing ones, then the top edge is $w_\mu$ minus the last $k$ characters, which must be all $1$s as $[\mu] \subseteq [\lambda]$.

Since 1-edges can only be bordered by NW-SE or NE-SW rhombi, all 1-edges are vertical, and they form paths from the 1-edges at the top to the 1-edges at the bottom. Furthermore, the number of ones at the top, $\lambda_1$, is equal to the number of total horizontal edges at the bottom of $\nabla_\lm$, so each horizontal edge on the bottom is labeled $1$.

As a consequence, the horizontal 0-edges must either be on the top boundary, or border two 0-triangles. By merging these pairs of 0-triangles and adding peaks on the top, we have recovered the lozenge tiling upon forgetting all the numbers.

It remains to show that the weights in the two formulas are the same. In the lozenge tiling, the weighted objects are the NW-SE rhombi, with a weight equal to the distance to the right edge of the shape. Under the correspondence, those rhombi become equivariant pieces, and since each row of the lozenge tiling becomes a band of the Knutson--Tao puzzle, the distance to the right edge is equal the height of the piece. This concludes the proof.
\end{proof}

\begin{example}
For our running example (Example~\ref{exmp:ooterms}), the lozenge tiling first gets transformed into the configuration shown in Figure \ref{fig:kt1}. Then, the horizontal strips get formed into a Knutson--Tao puzzle as shown in Figure \ref{fig:kt2}. The areas shaded in gray represent the final locations of the horizontal strips of the lozenge tiling. The weight of the original SSYT in Figure~\ref{ex:runningexSSYT} is $1\cdot 1 \cdot 1 \cdot 3 \cdot 5 \cdot 7$ which is the product of the heights of the equivariant pieces of the corresponding puzzle (shaded in darker gray).
\end{example}

\subsection{Is the Okounkov--Olshanski formula ``geometric''?} \label{sec: oof geometric?}
While a skew partition $\lm$ is, strictly speaking, essentially a pair of partitions $(\lambda, \mu)$, it is common to conflate this concept with the Young diagram $[\lm]$. While this is harmless for ordinary partitions, the map $(\lambda, \mu) \mapsto [\lm]$ is far from injective.

Nonetheless, many properties of a skew partition depend only on its Young diagram, such as $f^{\lm}$, or more generally, the functions $s_\lm$ or $\rpplmq$. We will call these properties \vocab{geometric}. To extend this notion further, some properties are \vocab{strongly geometric}, being invariant under symmetries of $\setz^2$, such as translation or taking the conjugate skew partition. For example, $f^{\lm}$ is strongly geometric, invariant under all these symmetries. However, $s_{\lm}$ is invariant under translations, but not under taking the conjugate partition. Also, the Naruse hook length formula is strongly geometric (see Section~\ref{sec:finalnaruse}). 

A natural question to ask is whether the Okounkov--Olshanski formula is geometric under these notions. While the original formulation \eqref{eq:oof} seems to involve $\mu$ and $\lambda$, this dependence could be superficial. The next result shows that the number $\OOT(\lm)$ of terms is geometric but not strongly geometric since this number is not necessarily invariant under vertical translation (see Example~\ref{ex: oof geometric?}).

\begin{proposition} \label{prop:ootgeometric}
The number of terms $\OOT(\lm)$ depends only on the shape of $[\lm]$ and the row numbers of the cells. In particular the number is invariant under horizontal translations.
\end{proposition}

\begin{proof}
By Corollary~\ref{cor:directbij} we have that $\OOT(\lm)=|\cSF(\lm)|$. It is clear that $|\cSF(\lm)|$ depends only on the shape of $[\lm]$ and the row numbers of each of the cells.
\end{proof}

Next, we look at the geometricity of the Okounkov--Olshanski formula itself. Proposition \ref{prop:ootgeometric} implies that the set of combinatorial objects indexing the terms of the formula is only dependent on $[\lm]$. However, the formula \eqref{eq:oof} is not geometric, with the ratios between terms varying depending on $\lambda$ and $\mu$ even as $[\lm]$ is fixed (see Example~\ref{ex: oof geometric?}). However, we can salvage an element of geometricity via the following observation.

\begin{proposition} \label{prop:emptycolumns}
For a positive integer $1 \leq k \leq \lambda_1$ such that $\lambda'_k = \mu'_k$, the weights corresponding to column $k$ for every term of the Okounkov--Olshanski formula (i.e.\ $k$th column in $[\mu]$, $\mathsf{d}_k(D)$ in the reverse excited diagram formulation, etc.), are equal to $\prod_{(i,k)\in [\lambda]} h(i, k)$.
\end{proposition}
\begin{proof}
The proof of this is probably most elegant in the flagged tableaux formulation, Corollary~\ref{cor:flag-form}. There, the weights corresponding to the $k$th column are associated with to $(k, t) \in M(T)$, with $1 \leq t \leq \lambda'_k$, with $i(k, t) = \lambda'_k$. Thus the weights are $\lambda_t + \lambda'_k - t - k + 1 = h(t, k)$,
as desired. 
\end{proof}

Therefore, we can essentially ignore columns in which there is no cell of $[\lm]$. A consequence of this concerns horizontal translation. If a column is inserted at the left end of the shape, then, in the reverse excited diagram formulation, the reverse excited diagram corresponding to a flagged tableau shifts to the right by one unit, along with its broken diagonals. Since $[\lambda]$ shifts to the right as well, the only change to the weights is the broken diagonals associated with the extra column, which cancel with the extra hook lengths in the $1/\prod_{u \in [\lambda]} h(u)$ prefactor.

On the other hand, adding columns at the right edge of the Young diagram, which would fix $[\lm]$, would fix the reverse excited diagrams, but change the weights. If one ignores weights associated with empty columns, the broken diagonals and prefactor remain the same, expressing the the Okounkov--Olshanski formula as a polynomial on the variable partition parts.

\begin{example}
Let $\lm = x21/x$, where $x \geq 2$ is a variable positive integer. Then, the Okounkov--Olshanski formula, in the reverse excited diagram form, states that
{\ytabset{boxsize=0.6em}
\[2 = \frac{3!}{(x+2)\cdot x\cdot 3 \cdot 1 \cdot 1}\paren{
\underbrace{x \cdot (x-1)}_{\ytab{{}{}{\bd}{\bd},\none{}{\dg}{\dg},\none\none{\dg}}} +
\underbrace{2 \cdot (x-1)}_{\ytab{{}{\dg}{}{\bd},\none{}{\bd}{\dg},\none\none{\dg}}} +
\underbrace{1 \cdot (x-1)}_{\ytab{{}{\dg}{}{\bd},\none{\dg}{}{\dg},\none\none{\bd}}} +
\underbrace{2 \cdot 1}_{\ytab{{}{\dg}{\dg}{},\none{}{\bd}{\bd},\none\none{\dg}}} +
\underbrace{1 \cdot 1}_{\ytab{{}{\dg}{\dg}{},\none{\dg}{}{\bd},\none\none{\bd}}}
},\]
}
which is true as a polynomial identity in $x$. (Empty columns have been omitted from both the calculation and the illustrations.)
\end{example}

We finish with an example of all the principles in this subsection. 

\begin{example} \label{ex: oof geometric?}
We consider the following five shapes: $22/2$, $32/3$, $33/31$, $2/\varnothing$, and $22/11$, shown below in cyan. All these shapes are geometrically the same, two  squares sharing an edge, and have one standard tableau.
\[ \ytabset{aligntableaux=top,textmode} 
\ytab{{}{},{\dg}{\dg}} \qquad
\ytab{{}{}{},{\dg}{\dg}} \qquad
\ytab{{}{}{},{}{\dg}{\dg}} \qquad
\ytab{{\dg}{\dg}} \qquad
\ytab{{}{\dg},{}{\dg}} \]

We will then proceed to evaluating the Okounkov--Olshanski formula in the reverse excited diagram form, we obtain the following:
{ \ytabset{boxsize=0.6em, centertableaux}
\begin{align*}
\ytab{{}{},{\dg}{\dg}}:\qquad 1 &= \frac{2!}{3 \cdot 2 \cdot 2 \cdot 1} \paren{\underbrace{2 \cdot 1}_{\ytab{{}{\bd}{\bd},\none{\dg}{\dg}}} + \underbrace{2 \cdot 1}_{\ytab{{\dg}{}{\bd},\none{\bd}{\dg}}}+\underbrace{2 \cdot 1}_{\ytab{{\dg}{\dg}{},\none{\bd}{\bd}}}} \\
\ytab{{}{}{},{\dg}{\dg}}:\qquad 1 &= \frac{2!}{4 \cdot 3 \cdot 1 \cdot 2 \cdot 1} \paren{\underbrace{3 \cdot 2 \cdot 1}_{\ytab{{}{\bd}{\bd}{\bd},\none{\dg}{\dg}}} + \underbrace{2 \cdot 2 \cdot 1}_{\ytab{{\dg}{}{\bd}{\bd},\none{\bd}{\dg}}}+\underbrace{2 \cdot 1 \cdot 1}_{\ytab{{\dg}{\dg}{}{\bd},\none{\bd}{\bd}}}} \\
\ytab{{}{}{},{}{\dg}{\dg}}:\qquad 1 &= \frac{2!}{4 \cdot 3 \cdot 2 \cdot 3 \cdot 2 \cdot 1} \paren{\underbrace{4 \cdot 3 \cdot 2 \cdot  1}_{\ytab{{\bd}{}{\bd}{\bd},\none{\bd}{\dg}{\dg}}} + \underbrace{4 \cdot 3 \cdot 2 \cdot  1}_{\ytab{{\bd}{\dg}{}{\bd},\none{\bd}{\bd}{\dg}}}+\underbrace{4 \cdot 3 \cdot 2 \cdot  1}_{\ytab{{\bd}{\dg}{\dg}{},\none{\bd}{\bd}{\bd}}}} \\
\ytab{{\dg}{\dg}}:\qquad 1 &= \frac{2!}{2 \cdot 1}(\underbrace{1}_{\ytab{{\dg}{\dg}}}) \\
\ytab{{}{\dg},{}{\dg}}:\qquad 1 &= \frac{2!}{3 \cdot 2 \cdot 2 \cdot 1} (\underbrace{3 \cdot 2}_{\ytab{{\bd}{}{\dg},\none{\bd}{\dg}}})
\end{align*}
}
Observe that with vertical translations or conjugation, the terms change, meaning that Proposition \ref{prop:ootgeometric} does not extend to those actions. Also, note that in $32/3$ and $33/31$, there are factors of $1$ and $4 \cdot 3$ which appear in both the prefactor and every term, as per Proposition \ref{prop:emptycolumns}.
\end{example}

\section{Determinantal Formulas Enumerating Okounkov--Olshanski Tableaux} \label{sec: det formulas for oot}
This section begins with a quick proof of Theorem \ref{thm:ooterms}, aided by characterizations of positive Okounkov--Olshanski terms as row-flagged or column-flagged tableaux. We then proceed to relate these results through the framework of outside decompositions. While it is false that $\OOT(\lm)$ obeys the same Hamel--Goulden identities as Schur functions, it is true in the special case of Lascoux--Pragacz decompositions. Crucial to our analysis will be a recent generalization of Hamel--Goulden determinantal identities to a generalization of Schur functions by Macdonald called \vocab{ninth variation Schur functions}, which are useful for enumerative purposes.

\subsection{Proof of Theorem \ref{thm:ooterms}}
By Section \ref{sec:posterms}, $\OOT(\lm) = \abs{\cSF(\lm)}$.

For this section, it is convenient to refine the  notion of \vocab{flagged tableaux} from Section~\ref{subsec:prelim_tableaux}.

\begin{definition}
If $\lm$ is a skew shape with $d$ rows, and $\vec a = (a_1, \ldots, a_d)$ and $\vec b = (b_1, \ldots, b_d)$ are nondecreasing tuples of integers, then define $\SSYT(\lm, \vec a, \vec b)$ to be the set of semistandard tableaux in $\SSYT(\lm)$ so that every entry in row $i$ is at least $a_i$ and at most $b_i$. In addition, we relax the condition that entries are positive, allowing entries to be any integer provided that the order conditions are still met.
\end{definition}

For example, $\SSYT(\lm, k) = \SSYT(\lm, 1^d, k^d)$ and  $\SSYT(\lm,\vec f)=\SSYT(\lm,1^d,\vec f)$ by abuse of notation. Note that the non-decreasing condition of ${\bf a}$ and ${\bf b}$ in the definition is essential for the enumeration that follows (e.g. see \cite{proctor2017row}).

This definition is useful due to the following result.

\begin{theorem}[Wachs \cite{wachs85}, Gessel--Viennot \cite{gv89}] \label{thm:flagenum}
\[ \abs{\SSYT(\lm, \vec a, \vec b)} = \det \brac*{\binom{\lambda_i - \mu_j + j - i + b_i - a_j}{b_i - a_j}}_{i,j=1}^d.\]
\end{theorem}

By definition, $\cSF(\lm) = \SSYT(\lm, 1^d, (1,2,\ldots))$, so
\[\abs{\OOT(\lm)} = \det \brac*{\binom{\lambda_i - \mu_j + j - 1}{i - 1}}_{i,j=1}^d,\]
reproducing the first half of Theorem \ref{thm:ooterms}.

To obtain the second half, we need two easy results.

\begin{proposition} \label{prop:columnflag}
A tableau $T$ is in $\cSF(\lm)$ if and only if $T \in \SSYT(\lm)$ with each entry in column $i$ is at most $\lambda'_i$. 
\end{proposition}
\begin{proof}
The forward direction is straightforward, since the largest row number in the $i$th column of $[\lm]$ is $\lambda'_i$. For the converse, consider an arbitrary cell $(k, i) \in [\lm]$. Then given a tableau $T$ with $T(\lambda'_i, i) \leq \lambda'_i$, as required by the condition, column-strictness implies that $T(k, i) \leq \lambda'_i- (\lambda'_i - k) = k$, as desired.
\end{proof}

Now consider $T \in \cSF(\lm)$. One can associate a tableau $T'$ with shape $(\lm)'$ with $T'(i,j) = T(j,i)+i-j$.

\begin{proposition}
The map $T \mapsto T'$ is a bijection from $\cSF(\lm)$ to $\SSYT((\lm)', (1-\mu'_1, 2-\mu'_2, \ldots), (1,2,\ldots))$.
\end{proposition}
\begin{proof}
Note that $T'(i+1,j)-T'(i,j) = 1+T(j,i+1)-T(j,i) \geq 1$ and $T'(i,j+1)-T'(i,j) = T(j+1,i)-T(j,i)-1 \geq 0$, so $T'$ satisfies the order conditions of semistandard tableaux. By Proposition \ref{prop:columnflag},
\[i - \mu'_i = 1 + i - (\mu' + 1) \leq T'(i, \mu'_i+1) \leq T'(i,j) \leq T'(i, \lambda'_i) \leq \lambda'_i + i - \lambda'_i = i,\]
so $T' \in \SSYT((\lm)', (1-\mu'_1, 2-\mu'_2, \ldots), (1,2,\ldots))$.

Given $T' \in \SSYT((\lm)', (1-\mu'_1, 2-\mu'_2, \ldots), (1,2,\ldots))$, the unique $T$ that can map to it must satisfy $T(i,j) = T'(i,j) + i - j$. Since $T(i+1,j)-T(i,j) = 1+T'(j,i+1)-T'(j,i) \geq 1$ and $T(i,j+1)-T(i,j) = T'(j+1,i)-T'(j,i)-1 \geq 0$, a $T$ defined this way would satisfy the order conditions for semistandard tableaux. Finally,
\[1 = (j - \mu'_j) + \mu'_j + 1 - j \leq T(i,j) \leq j + \lambda'_j - j = \lambda_j',\]
so $T \in \cSF(\lm)$ by Proposition \ref{prop:columnflag}.
\end{proof}

As a result, from Theorem \ref{thm:flagenum} we obtain
\[\OOT(\lm) = \abs*{\SSYT((\lm)', (1-\mu'_1, 2-\mu'_2, \ldots), (1,2,\ldots))} = \det \brac*{\binom{\lambda'_i}{\mu'_j + j - i}}_{i,j=1}^{\lambda_1}
,\]
as desired.

\begin{remark}
The proof Theorem \ref{thm:flagenum} in \cite{gv89} applies the Gessel-Viennot lemma to path systems associated with flagged tableaux. In this case, the paths correspond to the blue and green paths in Figure \ref{fig:oopathsrgb}. Indeed, one can successfully prove Theorem \ref{thm:flagenum} directly by counting lozenge tilings, though we consider it easier to stay away from a geometric argument.
\end{remark}

\subsection{Evaluations of Macdonald's ninth variation  Schur functions} \label{subsec:ninth variation Schur functions}
For the rest of this section, it is useful to work with the content $c(i,j) = j-i$ to describe the position of a cell, as opposed to rows or columns.

\begin{definition}
Given a shape $\lambda$ and an integer $c \in \setz$, let $d^{(\lambda)}_c$ be the largest row number $i$ so that $(i, i+c)$ is a cell in $[\lambda]$, and is $0$ otherwise. If the shape in context is clear, we may drop the superscript.
\end{definition}

The following result is a natural analogue to Proposition \ref{prop:columnflag}.
\begin{proposition} \label{prop:contentflag}
A tableau $T$ is in $\cSF(\lm)$ if and only if $T \in  \SSYT(\lm)$ satisfies $T(u) \leq d^{(\lambda)}_{c(u)}$ for all $u \in  [\lm]$.
\end{proposition}
\begin{proof}
Again the forward direction is straightforward, as a cell with content $c$ can't have a row number greater than $d^{(\lambda)}_{c(u)}$, by definition.

For the converse, note that for all $(i,j) \in [\lm]$, there exists a $k \geq 0$ so that $(i+k,j+k) \in [\lambda]$ and $d^{(\lambda)}_{j-i} = i + k$. In this case, $(i+k,j+k) \in [\lm]$ as well, since $(i,j) \notin [\mu]$. Then, if $T(u) \leq d^{(\lambda)}_{c(u)}$ for all $u$,
\[T(i,j) \leq T(i+k,j+k) - k \leq i,\]
as desired.
\end{proof}

The ninth variation Schur functions $s^{\dagger}_{\lm}$ (see Section~\ref{sec: generalizations Schur functions}) are suitable to count tableaux with this condition, as we next show.

\begin{definition}
For a shape $\lambda$, let $\vd^{(\lambda)} = (d^{(\lambda)}_{i,j})_{i \geq 1, j \in \setz}$ where $d^{(\lambda)}_{i,j}$ is $1$ if $i \leq d^{(\lambda)}_j$ and $0$ otherwise.
\end{definition}

Given this notation, an immediate corollary of Proposition \ref{prop:contentflag} is the following result.
\begin{corollary} \label{cor:sdagger2OOT}
\[\OOT(\lm) = s^\dagger_\lm(\vd^\lambda).\]
\end{corollary}

\begin{proof}
This follows from combining the definition of $s_{\lm}^{\dagger}$ \eqref{eq:def ninth variation} and Proposition~\ref{prop:contentflag}. 
\end{proof}

In \cite{hamelgoulden1995} it was proved that a Schur function of a skew shape can be written as a determinant of certain border strips; specifically, given an outside decomposition $(\theta_1, \theta_2, \ldots, \theta_n)$ of a shape \lm,
\[ s_\lm = \det\brac*{s_{\theta_i \hash \theta_j}}_{i,j=1}^n.\]
This identity unifies other determinantal identities of Schur functions like the Jacobi--Trudi formula and the {\em Lascoux--Pragacz identity} \cite{lascoux1988ribbon}. This identity has been recently shown to hold for these ninth variation Schur functions as well:
\begin{lemma}[\cite{bc20}, see also \cite{fk20}] \label{thm:generalhamelgoulden}
If $(\theta_1, \theta_2, \ldots, \theta_n)$ is an outside decomposition of a skew shape \lm,
\[ s^\dagger_\lm = \det\brac*{s^\dagger_{\theta_i \hash \theta_j}}_{i,j=1}^n.\]
\end{lemma}

Combining Lemma \ref{thm:generalhamelgoulden} and Corollary \ref{cor:sdagger2OOT}, implies that
\begin{equation} \label{eq:hgoot}
    \OOT(\lm) = \det\brac*{s^\dagger_{\theta_i \hash \theta_j}(\vd^{(\lambda))}}_{i,j=1}^n.  
\end{equation} 
Unfortunately, for general outside decompositions, the quantities in the determinant are rather complicated, since there is no clear relation between $\lambda$ and $\theta_i \hash \theta_j$. However, for the Lascoux--Pragacz decomposition, it turns out that, while placing $\theta_i \hash \theta_j$ as a subset of the border strip of $\lambda$, there is a $\nu$ so that $[\theta_i \hash \theta_j] = [\lambda/\nu]$. This implies the following result.

\begin{theorem} \label{cor:oolascouxpragacz}
If $(\theta_1, \theta_2, \ldots, \theta_n)$ is the Lascoux--Pragacz decomposition of a skew shape \lm,
\[ \OOT(\lm) = \det \brac*{\OOT(\theta_i \hash \theta_j)}_{i,j=1}^n,\]
where $\theta_i\hash \theta_j$ is placed on the border strip of $\lambda$. Note that $\OOT(\theta_i \hash \theta_j)$ is well-defined due to Proposition \ref{prop:ootgeometric}.
\end{theorem}

\begin{proof}
As discussed above, it suffices to show there is a shape $\nu$ so that $[\theta_i \hash \theta_j] = [\lambda/\nu]$. Let $\theta$ be the border strip of $\lambda$.

Let $S = [\lambda] \setminus [\theta_i \hash \theta_j]$. We want to show that $S = [\nu]$ for some $\nu$, or that for all $(i,j) \in S$, $(i-1,j), (i, j-1) \in S$, provided that the respective coordinates are both positive. Equivalently, we want $(i-1,j), (i, j-1) \notin [\theta_i \hash \theta_j]$. We will show that $(i, j-1) \notin [\theta_i \hash\theta_j]$; the second statement follows similarly.

First of all, if $(i,j) \notin [\theta]$, it is clear that $(i,j-1)$ cannot be in $[\theta]$, a superset of $[\theta_i \hash \theta_j]$. So assume that $(i,j) \in [\theta]$. Then, if $(i, j-1) \in [\theta_i \hash \theta_j]$, the maximum content of $\theta_j$ must be $j-i+1$.

Suppose that the cell of maximum content in $\theta_j$ is $(i-k, j-k-1)$ for some $k \geq 0$. Since $(i,j) \in [\theta]$, by the definition of the Lascoux--Pragacz decomposition, $(i-k, j-k) \notin [\lm]$. However this is absurd, since $(i,j) \in [\lambda]$ implies that  $(i-k, j-k) \in [\lambda]$, and $(i-k, j-k-1) \notin [\mu]$ implies that  $(i-k,j-k) \notin[\mu]$.
\end{proof}

Finally, to show how this framework subsumes the earlier theory of flagged tableaux enumeration, we give a second proof of Theorem \ref{thm:ooterms} using \eqref{eq:hgoot}.

\begin{proof}[Second proof of Theorem \ref{thm:ooterms}]
First let the cutting strip be horizontal. Then $\theta_i$ is the $i$th row of $[\lm]$, from contents $\mu_i - i + 1$ to $\lambda_i - i$. Therefore $\theta_i \hash \theta_j$ is a horizontal strip ranging from contents $\mu_i - i + 1$ to $\lambda_j - j$.

The expression $s^\dagger_{\theta_j \hash \theta_i}(\vd^{(\lambda)})$ counts the number of sequences of positive integers $a_{\mu_j - j + 1} \leq \cdots \leq a_{\lambda_i - i}$ so that $a_c \leq d^{(\lambda)}_c$ for all $c$. However, since $d^{(\lambda)}_c$ decreasing in $c$, this is equivalent to choosing a sequence of $\lambda_i - \mu_j - i + j$ nondecreasing positive integers at most $d^{(\lambda)}_{\lambda_i - i} = i$. It is well-known that this can be done in
\[\binom{(\lambda_i - \mu_j - i + j) + (i-1)}{i-1} = \binom{\lambda_i - \mu_j + j - 1}{i-1}.\]
ways, proving the first formula by \eqref{eq:hgoot}.

Now let the cutting strip be vertical. In this case $\theta_i$ is the $i$th column, with contents ranging from $i - \lambda'_i$ to $i - \mu'_i - 1$. Therefore $\theta_i \hash \theta_j$ is a vertical strip ranging from contents $i - \lambda'_i$ to $j - \mu'_j - 1$.

The expression $s^\dagger_{\theta_i \hash \theta_j}(\vd^{(\lambda)})$ now counts strictly decreasing sequences of positive integers $(a_c)$, with $i - \lambda'_i \leq c \leq j - \mu'_j - 1$, so that $a_c \leq d^{(\lambda)}_c$. Since $d^{(\lambda)}_{c+1} \geq d^{(\lambda)}_c - 1$, if the condition is satisfied for $c = i - \lambda'_i$, it is satisfied for all greater $c$ by the decreasing condition. Therefore
\[s^\dagger_{\theta_i \hash \theta_j}(\vd^{(\lambda)}) = \binom{d^{(\lambda)}_{i - \lambda'_i}}{\lambda'_i - \mu'_j - i + j} = \binom{\lambda'_i}{\lambda'_i - \mu'_j - i + j} = \binom{\lambda'_i}{\mu'_j + i - j},\]
implying the second half by \eqref{eq:hgoot}.
\end{proof}

\section{Okounkov--Olshanski terms for special families of skew shapes} \label{sec:special}
In this section we study the Okouknov-Olshanski formula for some special cases of skew shapes $\lm$, with a focus on enumeration.

\subsection{Zigzags}
For $n \geq 0$ define $\delta_n = (n-1,n-2,\ldots,1)$ and for $n > 1$ let $\sigma_n = \delta_{n+1}/\delta_{n-1}$ be a zigzag shape. For convenience define also $\delta_{0} = \varnothing$ so that $\sigma_1 = (1)$. Recall that $f^{\sigma_{n}}=E_{2n-1}$.

\begin{theorem}
\label{thm:genocchi}
\if\thesection1{For the zigzag $\sigma_n = (n,n-1,\ldots,1)/(n-2,n-3,\ldots,1)$, we have }\fi
$\OOT(\sigma_n) = G_{2n}$.
\end{theorem}

\begin{proof}
By Corollary~\ref{cor:directbij} it suffices to enumerate $\cSF(\sigma_n)$. By considering the reverse row word of such tableaux, i.e.\ the word created by reading the entries from the top right to the bottom left, $\OOT(\sigma_n)$ is the number of sequences $a_1, a_2, \ldots, a_{2n-1}$ so that $a_{2i-1} \geq a_{2i} < a_{2i+1}$ for $1 \leq i < n$, $a_{2i-1} \leq i$ for $1 \leq i \leq n$, and $a_{2i} \leq i$ for $1 \leq i < n$. These are exactly the strictly alternating pistols of length $2n - 1$ that are counted by $G_{2n}$.
\end{proof}

If we delete the bottom-most square from a zigzag, the number of nonzero Okounkov--Olshanski terms can be described using median Genocchi numbers:
\begin{theorem}
$\OOT((n+1,n,\ldots,2)/\delta_{n}) = H_{2n+1}$.
\end{theorem}
\begin{proof}
Similar to the proof of Theorem \ref{thm:genocchi}, $\OOT((n+1,n,\ldots,2)/\delta_n)$ is the number of strict alternating pistols of length $2n$, which is $H_{2n+1}$.
\end{proof}

We use the Okounkov--Olshanski formula to give an identity for odd Euler numbers.

\begin{corollary}
\begin{equation} \label{eq:oot_euler}
E_{2n-1} \,=\, \frac{(2n-1)!}{\prod_{i=1}^n (2n-2i+1)^i} \sum_{T \in \cSF(\sigma_n)} \prod_{(j,t) \in M(T)} (n-t+i(j,t)-j+2)
\end{equation}
\end{corollary}

\begin{proof}
The result follows by applying the flagged tableaux formulation, Corollary~\ref{cor:flag-form}, of the Okounkov--Olshanski formula to the shape $\sigma_n$ and using that $f^{\sigma_n}=E_{2n-1}$.
\end{proof}

\begin{figure}
    \centering
    \includegraphics[scale=0.6]{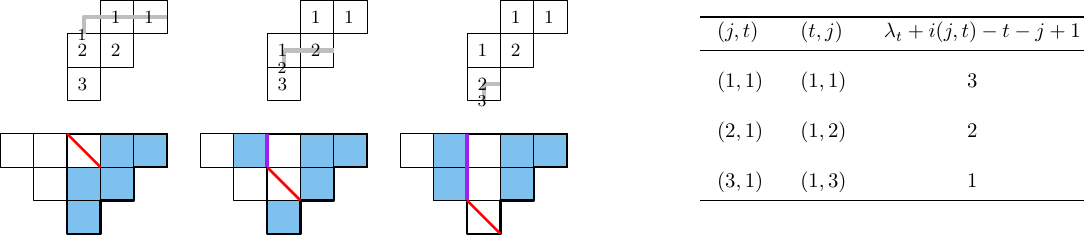}
    \caption{Flagged tableaux and reverse excited diagrams of the zigzag $\sigma_3=321/21$. Table of $(j,t) \in M(T)$ for each flagged tableaux $T$ and the restricted hook length of each $(t,j)$.}
    \label{fig: flagged tableaux reverse excited zigzag}
\end{figure}

\begin{example}
For the zigzag $\sigma_3$, we have that $\OOT(\sigma_3)=\abs{\cSF(\sigma_3)}=G_6 =3$. The flagged tableaux and corresponding reverse excited diagrams are in Figure~\ref{fig: flagged tableaux reverse excited zigzag} and \eqref{eq:oot_euler} in this case gives
\[
16 \,=\, E_5 \,=\, \frac{5!}{5\cdot 3^3} \left( 3 + 2 + 1 \right).
\]
\end{example}

Since $\OOT(\cdot)$ is not invariant under vertical translation (see \S \ref{sec: oof geometric?} and Example~\ref{ex: oof geometric?}), it is also interesting to consider shifted zigzags.
\begin{definition}
For $n > k \geq 0$, let $\sigma_{n}^{(k)}$ be $\sigma_{n-k}$, but shifted down $k$ units. This can be realized as the skew shape
\[\delta_{n+1}/(n+1, n, \ldots, n-k+2, n-k-1, n-k-2, \ldots, 1).\]
The notation stems from the fact that $\sigma_n^{(k)}$ is obtained from $\sigma_n$ by removing the first $k$ rows. Also, let $G_{2n}^{(k)} = \OOT(\sigma_n^{(k)})$.
\end{definition}
The quantity $G_{2n}^{(k)}$ is well defined due to  Proposition \ref{prop:ootgeometric}. In fact, it has a simple formula in terms of the Genocchi numbers.
\begin{proposition} \label{prop:shiftedgenocchi}
For nonnegative integers $n>k\geq 0$ we have that 
\[
G_{2n}^{(k)} = \frac{1}{2k+1}\binom{2n}{2k}G_{2(n-k)}.
\]
\end{proposition}

Let $g_{2n}^{(k)} \coloneqq \frac{1}{2k+1}\binom{2n}{2k}G_{2(n-k)}$.  We will show that both sequences of numbers obey the same recurrence relation given in the next two lemmas. 

\begin{lemma}[{Cigler \cite{cigler09}}] \label{rec: g}
\[
\binom{2n-k-1}{k} = \sum_{i=k}^{n-1} (-1)^{n-i-1} \binom{2i-k}{k} g_{2n}^{(i)}.
\]
\end{lemma}

\begin{proof}
The identity follows by taking the coefficient $t^k$ in the following identity of Cigler \cite[Cor. 1.3]{cigler09} for {\em Fibonacci polynomials} $F_n(t)\coloneqq\sum_{i=0}^n \binom{n-i-1}{i} t^i$: 
\[
F_{2n}(t) = \sum_{i=0}^{n-1} (-1)^{n-i-1} g_{2n}^{(i)} F_{2i+1}(t).
\]
\end{proof}

\begin{lemma} \label{rec:shifted G}
\[
\binom{2n-k-1}{k} = \sum_{i=k}^{n-1} (-1)^{n-i-1} \binom{2i-k}{k} G_{2n}^{(i)}.
\]
\end{lemma}

\begin{proof}
Observe that $G_{2n}^{(k)}$ is the number of sequences of positive integers $a_{2k+1}, a_{2k+2}, \ldots, a_{2n-1}$ with $a_i \leq \frac{i+1}{2}$ for all $2k+1 \leq i \leq 2n-1$, and $a_{2i-1} \geq a_{2i} < a_{2i+1}$ for all $k < i < n$. These are similar to the \vocab{truncated pistols} of \cite{zeng06}.

For integers $0 \leq i \leq j < n$, let $N_{i, j}$ be the the number of sequences of positive integers $a_{2i+1}, a_{2i+2}, \ldots, a_{2n-1}$ with $a_m \leq \frac{m+1}{2}$ for all $2i+1 \leq m \leq 2n-1$, and $a_{2m-1} \geq a_{2m} \geq a_{2m+1}$ for all $i < m \leq j$, and $a_{2m-1} \geq a_{2m} < a_{2m+1}$ for all $j < m < n$.

Now observe that for $k \leq j \leq n-2$, the only definition between the definitions $N_{k, j}$ and $N_{k, j+1}$ is the relative order of $a_{2j+2}$ and $a_{2j+3}$. Therefore $N_{k,j} + N_{k,j+1}$ is can be defined analogously to $N_{k,j}$ but with this condition removed. In this case, the integers $2(j-k+1)$ integers $a_{2k+1}, \ldots, a_{2j+2}$, which must be nonstrictly decreasing and thus all at most $k+1$, behave independently from $a_{2j+3}, \ldots, a_{2n-1}$. Namely,
\[N_{k,j} + N_{k, j+1} = \binom{2j-k+2}{k} G^{(j+1)}_{2n}.\]
The identity then follows by observing that $N_{k,k} = G^{(k)}_{2n}$ and $N_{k, n-1} = \binom{2n-k-1}{k}$.
\end{proof}

\begin{proof}[Proof of Proposition~\ref{prop:shiftedgenocchi}]
As $k$ varies from $0$ to $n-1$, the relation in Lemma~\ref{rec:shifted G} allows $G_{2n}^{(k)}$ to be calculated for all $k$. Since by Lemma~\ref{rec: g}, the numbers $g_{2n}^{(k)}$ satisfy the same recurrence the result follows.
\end{proof}

\begin{proof}[Proof of Theorem~\ref{thm:ootthickzigzag}]
The shape $\delta_{n+2k}/\delta_n$ has a Lascoux--Pragacz decomposition $(\theta_1, \ldots, \theta_k)$ where $\theta_i = \delta_{n+2i}/\delta_{n+2i-2}$. Therefore, $\theta_i \hash \theta_j$ is a zigzag running from contents $2-n-2i$ to $n+2j-2$. Placing this on the rim of $\delta_{n+2k}$, it is a zigzag with $2n+2j+2i-3$ cells, shifted down $\frac 12 ((n+2k-2)-(n+2j-2)) = k-j$ units, meaning that $\OOT(\theta_i \hash \theta_j) = G_{2(n+i+k-1)}^{(k-j)}$. By Theorem~\ref{cor:oolascouxpragacz}, we now have
\[\OOT(\delta_{n+2k}/\delta_n) = \det\brac*{G_{2(n+i+k-1)}^{(k-j)}}_{i,j=1}^{k}.\]
By Proposition \ref{prop:shiftedgenocchi}, this is therefore
\begin{align*}
\lefteqn{\det\brac*{\frac{1}{2(k-j)+1}\binom{2(n+i+k-1)}{2(k-j)}G_{2(n+i+j-1)}}_{i,j=1}^{k} =}\\
& \det\brac*{\frac{(2(n+i+k-1))!}{(2(k-j)+1)!}\hat{G}_{2(n+i+j-1)}}_{i,j=1}^{k} = \prod_{i=1}^k \frac{(2n+2i+2k-2)!}{(2i-1)!} \det \brac*{\hat{G}_{2(n+i+j-1)}}_{i,j=1}^{k},\end{align*}
as claimed.
\end{proof}

\begin{remark}
 Theorem~\ref{thm:ootthickzigzag} for $k=1$ gives $\OOT(\delta_{n+2}/\delta_n)=G_{n+2}$ and for $k=2$ gives  $\OOT(\delta_{n+4}/\delta_n)$  as a determinant of a  $2\times 2$ matrix. We have added the latter sequence to \cite[\href{http://oeis.org/A336674}{A336674}]{oeis}. 
\end{remark}

Interestingly, if $n=0,1$, the identity in Theorem~\ref{thm:ootthickzigzag} can be reversed to give formulas for Hankel determinants. For example, since the Okounkov--Olshanski tableaux for ordinary shapes are simply empty tableaux, $\OOT(\delta_{2k}/\delta_{0}) = \OOT(\delta_{2k+1}/\delta_1) = 1$. This yields the following result.

\begin{corollary} \label{cor:genocchihankel}
\[\det\brac*{\hat{G}_{2(i+j+1)}}_{i,j=0}^{k-1} = \prod_{i=0}^{k-1} \frac{(2i+1)!}{(2i+2k)!}, \quad 
\text{and} \quad
\det\brac*{\hat{G}_{2(i+j+2)}}_{i,j=0}^{k-1} = \prod_{i=0}^{k-1} \frac{(2i+1)!}{(2i+2k+2)!}.\]
\end{corollary}

Hankel determinants, such as these, are often evaluated through the theory of continued fractions and orthogonal polynomials (see \cite{adc2}). While we have been unable to find the exact statement of Corollary \ref{cor:genocchihankel} in the literature, a slightly modified version has been well-studied through the Euler numbers $E_{2n-1} = (-1)^{n-1}2^{2n-1}(4^n-1)B_{2n}/n$. Indeed, after defining $\hat E_{2n-1} = E_{2n-1}/(2n-1)!$, we have $\hat E_{2n-1} = 2^{2n-1} \hat G_{2n}$\footnote{We refer the reader to the recent article by Han \cite{han20} for a useful compilation of these identities, though some of them date back over a century.}, which means we rewrite Theorem \ref{thm:ootthickzigzag} as follows.

\begin{corollary}  \label{cor:ootthickzigzag-Euler}
\[\OOT(\delta_{n+2k}/\delta_n) = 2^{-k(2n+2k-1)}\prod_{i=1}^k \frac{(2n+2i+2k-2)!}{(2i-1)!} \det \brac*{\hat{E}_{2(n+i+j-1)-1}}_{i,j=1}^{k},\]
\end{corollary}

This leads to the extremely mysterious identity that states that for thick zigzags $\theta$, the number $f^{\theta}$ of SYT is proportional to the number $\OOT(\theta)$ of Okounkov--Olshanski terms.

\begin{corollary} \label{cor: rel SYTzigzag OOTzigzag}
\begin{equation} \label{eq: rel SYTzigzag OOTzigzag}  
f^{\delta_{n+2k}/\delta_n} \,=\, 2^{k(2n+2k-1)}\,(k(2n+2k-1))!\, \prod_{i=1}^k \frac{(2i-1)!}{(2n+2i+2k-2)!}\,\OOT(\delta_{n+2k}/\delta_n).
\end{equation}
\end{corollary}

\begin{proof}
The result follows by combining 
Corollary~\ref{cor:ootthickzigzag-Euler} with the identity in \cite[Cor. 8.9]{mpp2} that states 
\[
f^{\delta_{n+2k}/\delta_n}/\abs{\delta_{n+2k}/\delta_n}! \,=\, \det \brac*{\hat{E}_{2(n+i+j-1)-1}}_{i,j=1}^{k}.
\] 
\end{proof}

\subsection{An asymptotic application to thick zigzags} \label{sec:asymptotic thick zigzags}
In \cite{mpp_asymptotics}, the asymptotic behavior the quantity $f^{\delta_{2k}/\delta_{k}}$ is studied, where $n = \abs{\delta_{2k}/\delta_{k}}$. In this section we use Corollary~\ref{cor: rel SYTzigzag OOTzigzag} to give bounds for $\log f^{\delta_{2k}/\delta_{k}}$.

\begin{proof}[Proof of Theorem~\ref{thm: bounds thick zigzags}]
Since $f^{\delta_{2k}/\delta_k} \leq f^{\delta_{2k+1}/\delta_{k+1}} \leq f^{\delta_{2k+2}/\delta_{k+1}}$, we will bound $f^{\delta_{2k+1}/\delta_{k+1}}$ instead. Redefine $n = \abs{\delta_{2k+1}/\delta_{k+1}}$, which is asymptotically equivalent to the original $n$. Also, assume that $k$ is even.

With the notation of superfactorials in Section~\ref{sec:superfactorials}, \eqref{eq: rel SYTzigzag OOTzigzag} becomes
\[f^{\delta_{2k+1}/\delta_{k+1}} = n!\, 2^n \frac{\Psi(k/2)\,\Psi(2k)\,\Phi(3k)}{\Psi(3k/2)\,\Phi(4k)} \OOT(\delta_{2k+1}/\delta_{k+1}).\]
Using known asymptotics for $\Phi(\cdot)$ and $\Psi(\cdot)$ (e.g. see \cite[\SS 2.10]{mpp3}) and since $n  \sim \frac{3}{2} k^2$ 
, we have
\begin{multline*}
    \log f^{\delta_{2k+1}/\delta_{k+1}} - \log \OOT(\delta_{2k+1}/\delta_{k+1}) \sim  n \log n - n + n \log 2 + \frac{k^2}{4} \log(k/2) + \frac{k^2}{4} \paren*{\log 2 - \frac{3}{2}} \\- \frac{9k^2}{4} \log(3k/2) - \frac{9k^2}{4}\paren*{\log 2 - \frac{3}{2}} + 4k^2 \log(2k) + 4k^2 \paren*{\log 2 - \frac{3}{2}} \\+ \frac{9}{2}k^2 \log(3k) - \frac{27}{4} k^2 - 8k^2 \log(4k) + 12k^2 + o(n).
\end{multline*}

The terms involving $k$ simplifies to 
\[-\frac{3}{2}k^2 \log k + k^2 \paren*{\frac{9}{4} - 8\log 2 + \frac{9}{4} \log 3} = -\frac{1}{2}n \log n + n \paren*{\frac{3}{2} - \frac{35}{6} \log 2 + 2 \log 3} + o(n).\]
Combining everything, we have
\begin{equation} \label{eq: pf bound f thick zigzag}
\log f^{\delta_{2k+1}/\delta_{k+1}} - \log \OOT(\delta_{2k+1}/\delta_{k+1}) = \frac{1}{2} n\log n + n\paren*{\frac{1}{2} - \frac{29}{6} \log 2 + 2 \log 3} + o(n),
\end{equation}
where the constant is $\approx -0.65299$.

We now bound $\OOT(\delta_{2k+1}/\delta_{k+1})$. We describe tableaux in $\cOOT(\delta_{2k+1}/\delta_{k+1})$ as tableaux $T$ in $\SSYT(\delta_{k+1}, 2k)$ where $T(u) > c(u)$ for all $u$. This is automatically satisfied when if $T(u) > k$ for all $u$, so we can bound
\begin{equation} \label{eq: pf bound oot thick zigzag}
    \log \abs{\SSYT(\delta_{k+1}, k)} \leq \log \OOT(\delta_{2k+1}/\delta_{k+1}) \leq \log \abs{\SSYT(\delta_{k+1}, 2k)}.
\end{equation}
Next, by a direct calculation from the hook-content formula (e.g. \cite[Cor. 7.21.4]{EC2}) we have that $\abs{\SSYT(\delta_{k+1},k)} \,=\, 2^{\binom{k}{2}}$ and $\abs{\SSYT(\delta_{k+1},2k)} \,=\,\frac{\Psi(3k/2)\,\Phi(k-1)}{\Psi(k/2)\,\Phi(2k-1)\,\Lambda(k)}$  (see \cite[Cor. 5.10]{mpp3}).  Thus the LHS of \eqref{eq: pf bound oot thick zigzag} becomes
\[\log 2^{\binom{k}{2}} = \frac{n\,\log 2}{3} + o(n),\]
while the RHS of \eqref{eq: pf bound oot thick zigzag} simplifies to
\begin{align*}
\lefteqn{\log \frac{\Psi(3k/2)\,\Phi(k-1)}{\Psi(k/2)\,\Phi(2k-1)\,\Lambda(k)}  = }\\
&= k^2 \paren*{\frac{9}{4}\log 3 - \frac{27}{8} + \frac{3}{8} - \frac{3}{4} - 2\log 2 + 3 - \frac{1}{2}\log 2 + \frac{3}{4}} + o(k^2) \\
&= k^2 \paren*{\frac{5}{2}\log 2 + \frac{9}{4} \log 3} + o(k^2) \\
&= n \paren*{ - \frac{5}{3} \log 2 + \frac{3}{2} \log 3} + o(n).
\end{align*}
Thus \eqref{eq: pf bound oot thick zigzag} becomes
\[
  \frac{n\,\log 2}{3} + o(n)   \,\leq\, \log \OOT(\delta_{2k+1}/\delta_{k+1}) \,\leq\, n \paren*{ - \frac{5}{3} \log 2 + \frac{3}{2} \log 3} + o(n).
\]
Combining this with \eqref{eq: pf bound f thick zigzag} gives the stated bounds.
\end{proof}

\begin{remark}
The RHS in \eqref{eq:bounds thick zigzags} numerically evaluates to $\approx -0.4219$ and the LHS to $\approx -0.1603$. While the lower bound is weaker than previously known results, $\approx -0.3237$), the upper bound is sharper than $\approx -0.0621$ and quite close to the conjectured value of $\approx -0.1842$ \cite[Thm. 1.1, \S 13.7]{mpp_asymptotics} which is known to be a constant by the main result in \cite{MPT}. See also \cite{Sun} and \cite{Gord} for recent results on the existence and conjectured description of limit curves for general skew stable limit shapes (the diagram $[\lambda]$ scaled by $1/\sqrt{n}$ in both directions converges to a curve) like $\delta_{2k+1}/\delta_{k+1}$.
\end{remark}

\subsection{Slim shapes and rectangles} \label{subsec: slim shapes and rects}
Under certain circumstances, the Okounkov--Olshanski formula becomes ``degenerate'' and begins to resemble a product formula. This phenomenon first manifests itself with the number of terms.

\begin{proposition} \label{prop:oot_is_ssytmud}
For a shape \lm{} with $\mu_1 \leq \lambda_d$, $\cOOT(\lm) = \SSYT(\mu, d)$.
\end{proposition}
\begin{proof}
Note that for all $T \in \SSYT(\mu, d)$, and all $u \in [\mu]$,
\[\lambda_{d+1-T(u)} \geq \lambda_d \geq \mu_1 > \mu_1 - 1 \geq c(u). \qedhere\]
\end{proof}

Therefore, in this case Okounkov--Olshanski terms can be enumerated from Stanley's hook content formula \cite{EC2}:

\begin{corollary}
For a shape \lm{} with $\mu_1 \leq \lambda_d$,
\[ \OOT(\lm) = \prod_{u \in [\mu]} \frac{d + c(u)}{h(u)}.\]
\end{corollary}

As special cases of the result above, we obtain that for thick hooks $\lm = (b+c')^{a+c}/b^a$ and slim shapes $\lambda/\delta_{d+1}$, the number of terms $\OOT(\lm)$ is given by the {\em MacMahon box formula} for parameters $a,b,c$ and by a power of $2$. The proof follows from known special cases of the hook-content formula (e.g. see \cite[Thm. 7.21.7]{EC2} and \cite[Prop. 10.3]{mpp_asymptotics}). 

\begin{corollary}
For nonnegative integers $a,b,c,c'$ and $\lm = (b+c')^{a+c}/b^a$ we have that 
\begin{equation} \label{eq:oof_is_MacMahon}
     \OOT(\lm) = \prod_{i=1}^a\prod_{j=1}^b \prod_{k=1}^c \frac{i+j+k-1}{i+j+k-2}.
\end{equation}
\end{corollary}

\begin{corollary}
For $\lm = \lambda/\delta_d$ where $\lambda_d\geq d$, we have that 
\[
\OOT(\lambda/\delta_{d+1}) \,=\, 2^{\binom{d}{2}}.
\]
\end{corollary}

\begin{remark} \label{rmk: equal terms}
If we strengthen the condition of slim shapes to $\lambda_d \geq \mu_1 + d - 1$, as discussed in Section \ref{subsec:ooesum}, we have a coincidence of sets $\cale(\lm) = \cOOE(\lm) = \cale(\Lambda^d/\mu)$. This leads to the curious result that there exists a correspondence between $\cale(\lm)$ and $\cRE(\lm)$. As this problem has a strong geometric flavor, it would be interesting to see if there exists an intuitive reason for the existence of such a correspondence, since all bijections we have constructed either require broken diagonals (as in Lemma \ref{lem:redaffine}), or the column-by-column analysis of semistandard tableaux (as in Corollary \ref{cor:directbij})
\end{remark}

Finally, we turn to the case where $\lambda = a^d$ is a rectangle. Here, $\lambda_d = \lambda_1 \geq \mu_1$, so the number of term is given by a product formula. Moreover, the weight of every term becomes the same, namely $\prod_{u \in [\mu]} (a - c(u))$. As a result, one attains the following result:
\begin{proposition}
\[f^{a^d/\mu} = \frac{(ad-\abs{\mu})!}{\prod_{i=1}^d \prod_{j=1}^a (i+j-1)} \prod_{u \in [\mu]} \frac{(a-c(u))(d+c(u))}{h(u)},\]
where $h(u)$ is the hook length of cell $u$ with respect to $\mu$.
\end{proposition}

\begin{remark}
Of course, the existence of such a product formula is hardly surprising, given that $a^d/\mu$ must be an ordinary shape rotated 180 degrees, making $f^{a^d/\mu}$ a quantity that can be evaluated via the ordinary hook length formula. Nonetheless, the equality of these two expressions is not immediately obvious and therefore this result may be useful in some applications, especially when $\mu$ is relatively small and but also complicated.
\end{remark}

\section{A Reverse Plane Partition \texorpdfstring{$q$}{q}-Analogue} \label{sec:qoo}

The aim of this section is to give $q$-analogues of the Okounkov--Olshanski formula for reverse plane partition of skew shapes in a manner similar to the Okounkov--Olshanski formula. 

\begin{remark}
Recall that the generating function of reverse plane partitions and semistandard Young tableaux of straight shapes differ only by a power of $q$ (e.g. see \cite[Cor. 7.21.3, Thm. 7.22.1]{EC2}). This is not true for skew shapes and so our reverse plane partition $q$-analogue of \eqref{eq:oof} is different to the SSYT $q$-analogue of Chen--Stanley \cite{chenstanley2016qoo} (see Theorem~\ref{thm: Chen-Stanley q-analogue}). 
\end{remark}

We first state a version of such a result using the language of Grothendieck polynomials; we will connect it to the Okounkov--Olshanski formula in Section \ref{subsec: qoof}. To state the result we need the following notation, let $G_\lm(\vy) \coloneqq G_\mu\mpar{\ominus \vy_\lambda}{\vy}$. The following result can be viewed as a $q$-analogue of Theorem~\ref{thm:flmintermsofs}.

\begin{theorem} \label{thm:oorpp}
For a skew shape $\lm$ we have that 
\[\rppratio = G_\lm(1-q^{-1}, 1-q^{-2}, \ldots).\]
\end{theorem}

We present two proofs of this result: one involving equivariant $K$-theory of Grassmannians and the other involving determinantal manipulations.



\subsection{Equivariant \texorpdfstring{$K$}{K}-theory proof}
Let $X_\mu$ be the Schubert variety in equivariant $K$-theory corresponding to a partition $\mu$. Then let $[O_\mu]$ be the class of the structure sheaf of $X_\mu$ (eg. see \cite[\S 2.2]{grahamkreiman2015}). There are localization maps corresponding to permutations that can be applied $[O_\mu]$; let the image of $[O_\mu]$ under the map corresponding to a partition $\lambda$ be $\eval{[O_\mu]}_\lambda$, which is a Laurent polynomial of expressions of the form $e^{\epsilon_i}$ ($1\leq i \leq n$). Our usage of results regarding equivariant $K$-theory is limited to these two results, after which Theorem \ref{thm:oorpp} follows easily.

\begin{theorem}[{Graham-Kreiman \cite[Thm. 4.15]{grahamkreiman2015}}] \label{thm:grothendieckisk}
\[\eval{[O_\mu]}_\lambda = \eval{G_\lm(\vy)}_{y_i=1-e^{\epsilon_i}}\]
\end{theorem}

\begin{theorem}[{Naruse-Okada \cite[Thm. 5.2]{naruseokada2018}}] \label{thm:naruseokadarppk}
\[\rpplmq = \eval{
\frac{\eval{[O_\mu]}_\lambda}{\eval{[O_\lambda]}_\lambda}
}_{e^{\alpha_i} = q},\]
where $\alpha_i = \epsilon_{i+d} - \epsilon_{i+d+1}$.
\end{theorem}

\begin{proof}[Proof of Theorem \ref{thm:oorpp}]
Since the factorial Grothendieck polynomial $G_\varnothing$ is $1$,
\[\rppratio = \eval{\frac{\eval{[O_\mu]}_\lambda/\eval{[O_\lambda]}_\lambda}{\eval{[O_\varnothing]}_\lambda/\eval{[O_\lambda]}_\lambda}}_{e^{\alpha_i} = q} = \eval{\eval{[O_\mu]}_\lambda}_{e^{\alpha_i} = q}.\]
Now, if the $e^{\alpha_i}$ are set equal to $q$, we get that $e^{\epsilon_i} = cq^{-i}$ for some $c$. Therefore
\[\rppratio = G_\lm(1-cq^{-1}, 1-cq^{-2}, \ldots)\]
for some constant $c$.

We claim that $c$ is irrelevant in the above expression. Indeed,
\[\ominus(1-cq^i) \oplus (1-cq^j) = (1-c^{-1}q^{-i}) \oplus (1-cq^j) = 1-q^{j-i}.\]
Thus we can set $c=1$, completing the proof.
\end{proof}

\subsection{Determinantal algebra proof}
An alternate proof of Theorem \ref{thm:oorpp} follows from expanding the Grothendieck as a determinant and comparing to known expressions for $\rpplmq$. For this section only let $[a] \coloneqq 1-q^a$, $[n]! \coloneqq\prod_{i=1}^n [i]$, $\ell_i \coloneqq \lambda_i + d - i$, and $m_i \coloneqq \mu_i + d - i$.

\begin{theorem}[{Ikeda-Naruse \cite[Eq. (2.12)]{ikedanaruse2013grothendieck}}] \label{thm:bialternantGrothendieck}
If $\vx = (x_1, x_2, \ldots, x_n)$, then
\[G_\mu\mpar{\vx}{\va} = \frac{1}{\prod_{i<j} (x_i - x_j)} \det \brac*{[x_i\vert \va]^{\mu_j+n-j}(1-x_i)^{j-1}}^n_{i,j=1} \]
where $[x_i\vert \va]^k = (x_i \oplus a_1)(x_i \oplus a_2)\cdots(x_i \oplus a_k)$.
\end{theorem}
\begin{theorem}[{Krattenthaler \cite[Cor. 8]{krattenthaler89}}] \label{thm:kratt}
\[
\rpplmq = \det \brac*{\frac{q^{\binom{-\lambda_i + i - j}{2}-\binom{-\lambda_i}{2}}}{[\lambda_i - i - \mu_j + j]!}}_{i,j=1}^d
\]
\end{theorem}
\begin{remark}
Krattenthaler uses a definition of reverse plane partitions so that the minimum entry is $1$, so the result in \cite{krattenthaler89} has an extra factor of $q^{\abs{\lm}}$.
\end{remark}
We now write this in a more convenient form:
\begin{corollary} \label{cor:krattcor}
\[\rpplmq = \det \brac*{\frac{q^{\ell_i(j-i)}}{[\lambda_i - i - \mu_j + j]!}}_{i,j=1}^d\]
\end{corollary}
\begin{proof}
We have that
\[
\paren*{\binom{-\lambda_i + i - j}{2}-\binom{-\lambda_i}{2}} - \ell_i(j-i) = \paren*{(d-1)i - \binom{i}{2}} - \paren*{(d-1)j - \binom{j}{2}},\]
so the result follows from simple row manipulations on Theorem \ref{thm:kratt}.
\end{proof}

The next lemma is a $q$-analogue of the ordinary product formula for $f^{\lambda}$ \cite[Thm. 14.5.1]{adinroichman}.

\begin{lemma} \label{lem:ordproductq}
For a partition $\lambda$ we have that 
\[\rpplq = \frac{\prod_{i<j} [\ell_i - \ell_j]}{\prod_i [\ell_i]!}.\]
\end{lemma}
\begin{proof}
The result follows by Stanley's $q$-analogue of the hook length formula \cite{stanley1971theory} \cite[Cor. 7.21.3]{EC2}
\begin{equation}
    \rpplq = \prod_{u \in [\lambda]} \frac{1}{[h(u)]},
\end{equation}
and the fact from \cite[Lem. 7.21.1]{EC2} that
$\prod_{u\in \lambda} [h(u)] = \prod_{i\geq 1} [\ell_i]!/\prod_{i<j} [\ell_i-\ell_j]$.
\end{proof}

\begin{proof}[Proof of Theorem \ref{thm:oorpp}]
Since $[a] \oplus [b] = [a+b]$, applying Theorem~\ref{thm:bialternantGrothendieck} results in the following:
\begin{align*}
G_\lm([-1],[-2],\ldots) &= G_\mu\mpar{[\ell_1 + 1],[\ell_2 + 1],\ldots,[\ell_d + 1]}{[-1],[-2],\ldots} \\
&= \frac{\det \brac*{[\ell_i][\ell_i - 1]\cdots[\ell_i - m_j + 1] q^{(\ell_i+1)(j-1)}}^d_{i,j=1}}{\prod_{i<j} ([\ell_i+1]-[\ell_j+1])}.
\end{align*}
Since $[a] - [b] = q^a - q^b = q^b[a-b]$, this can be rewritten as
\[
\frac{1}{\prod_{i<j} ([\ell_i-\ell_j]q^{\ell_j+1})} \det \brac*{\frac{[\ell_i]!}{[\ell_i - m_j]!}q^{(\ell_i+1)(j-1)}}^d_{i,j=1}.
\]
We now pull out all the factors of $[\ell_i]!$ from the determinant. Also, in the product, observe that for all $j$ there are exactly $j - 1$ possibilities for $i$. So the expression simplifies to
\[\frac{\prod_i [\ell_i]!}{\prod_{i<j} [\ell_i - \ell_j]} \det \brac*{\frac{q^{(\ell_i+1)(j-1)}}{[\ell_i - m_j]!}}^d_{i,j=1}\cdot q^{-\sum_i (\ell_i+1)(i-1)} = \frac{\prod_i [\ell_i]!}{\prod_{i<j} [\ell_i - \ell_j]} \det \brac*{\frac{q^{(\ell_i+1)(j-i)}}{[\ell_i - m_j]!}}^d_{i,j=1}.\]

By Lemma \ref{lem:ordproductq} we have therefore proven that
\[G_\lm([-1], [-2], \ldots) \cdot \rpplq = \det \brac*{\frac{q^{(\ell_i+1)(j-i)}}{[\ell_i - m_j]!}}^d_{i,j=1}.\]
This, combined with Corollary \ref{cor:krattcor}, concludes the proof.
\end{proof}

\subsection{Evaluating the Grothendieck polynomial} \label{subsec: qoof}
Expanding out Theorem \ref{thm:oorpp} in the obvious way leads to a sum over set-valued semistandard tableaux. In order to manipulate this into an expression resembling \eqref{eq:oof}, we need a technical result that categorizes set-valued semistandard tableaux by the maximum entries in each cell.
\begin{definition} \label{def:mtu}
Given a semistandard tableau $T$ of shape $\mu$ and a cell $u \in [\mu]$, let $m_T(u)$ be the minimum $k \leq T(u)$ such that replacing $T(u)$ with $k$ still results in a semistandard tableau. Similarly, define $M_T(u)$ to be the maximum such $k \geq T(u)$, or $\infty$ if such a $k$ does not exist.
\end{definition}
\begin{proposition} \label{prop:sssytdecomp}
The set of all set-valued tableaux of shape $\mu$ admits parallel decompositions
\begin{align*}
\SVT(\mu)&= \bigsqcup_{T_0 \in \SSYT(\mu)}\: \set[\big]{T \mid T(u) \subseteq \set{m_{T_0}(u), m_{T_0}(u)+1, \ldots, T_0(u)}, T_0(u) \in T(u)} \\
&=\bigsqcup_{T_0 \in \SSYT(\mu)}\: \set[\big]{T \mid T(u) \subseteq \set{T_0(u), T_0(u)+1, \ldots, M_{T_0}(u)}, T_0(u) \in T(u)}.
\end{align*}
\end{proposition}
\begin{proof}
Given  $T$ in $\SVT(\mu)$, let $\max T$ be the tableau such that $(\max T)(u) = \max T(u)$. We claim that the set in the first decomposition corresponding to a particular $T_0$ is the set of tableaux $T$ such that $\max T = T_0$.

Consider a tableau $T$ so that $\max T = T_0$. Pick a cell $u$ and a value $t \in T(u)$. If we choose $t$ for the cell $u$ and $T_0(u')$ for every other cell $u'$, we obtain that $m_{T_0}(u) \leq t \leq T_0(u)$. Obviously $T_0(u) \in T(u)$, so $T$ is indeed in the desired set.

Now consider a tableau $T$ that lies in this set; clearly $\max T = T_0$, so we only need to check that $T$ is a valid set-valued tableau. To show this, note that if $u_1, u_2$ are cells so that $u_1$ is immediately to the left of $u_2$ then $t_1 \leq T_0(u_1) \leq m_{T_0}(u_2) \leq t_2$ for any $t_1 \in T(u_1)$ and $t_2 \in T(u_2)$. Similarly, if $u_1$ were immediately above $u_2$, we would have $t_1 \leq T_0(u_1) < m_{T_0}(u_2) \leq t_2$ for any such $t_1$, $t_2$. This proves the first decomposition. The proof of the second is almost exactly the same, so we omit it here.
\end{proof}

The following proof of Theorem \ref{thm:oorpptableaux} now follows analogously from the proof of Theorem \ref{thm:oof} given in Section \ref{sec:background}.

\begin{proof}[Proof of Theorem \ref{thm:oorpptableaux}]
Since factorial Grothendieck polynomials are symmetric \cite{mcnamara2006}, we have that
\begin{align*}
\rppratio &= G_\mu\mpar{1 - q^{\lambda_d + 1}, 1 - q^{\lambda_{d-1} + 2}, \ldots}{1 - q^{-1}, 1-q^{-2}, \ldots} \\
&= \sum_{T\in \SVT(\mu,d)} (-1)^{\abs{T} - \abs{\mu}} \prod_{\substack{u\in [\mu]\\r\in T(u)}} (1-q^{w(u,r)}) \\
&= \sum_{T\in \SVT(\mu,d)} (-1)^{\abs{\mu}} \prod_{u\in [\mu]} \prod_{r\in T(u)} (q^{w(u,r)}-1).
\end{align*}

By using the first part of Proposition \ref{prop:sssytdecomp}, this is
\begin{align*}
\MoveEqLeft \sum_{T \in \SSYT(\mu, d)} \prod_{u \in [\mu]} \paren*{ (1 - q^{w(u,T(u))}) \sum_{S \subseteq \set{m_T(u), m_T(u)+1, \ldots, T(u)-1}} \prod_{k \in S} (q^{w(u,k)} - 1)} \\
&= \sum_{T \in \SSYT(\mu, d)} \prod_{u \in [\mu]} \paren*{(1 - q^{w(u,T(u))}) \prod_{m_T(u) \leq k < T(u)} q^{w(u,k)}} \\
&= \sum_{T \in \SSYT(\mu, d)}q^{p(T)}\prod_{u \in [\mu]} (1 - q^{w(u,T(u))}). \qedhere
\end{align*}
\end{proof}

Using extremely similar methods, one can also deduce Theorem~\ref{thm:oorpptalt}.

\begin{proof}[Proof of Theorem~\ref{thm:oorpptalt}]
The proof is completely analogous to that of Theorem~\ref{thm:oorpptableaux} where we use the second part of Proposition~\ref{prop:sssytdecomp} instead.
\end{proof}



\begin{example} \label{ex:rpp q-oof}
For the shape $\lm=2221/11$, we have
\[
\rppratio = 1 - q^5 - q^6 - q^7 + q^8 + q^{10}.
\]
Applying Theorems \ref{thm:oorpptableaux} and \ref{thm:oorpptalt} results in two different expressions for this polynomial as a sum over terms corresponding to the six tableaux in $\SSYT(\lambda,4)$ in \eqref{eq:oo-ex-ssyt}:
\begin{align} 
\label{eq: 1st qanalogue rpp} \lefteqn{\rppratio  =} \notag \\
&=\begin{multlined}[t]
q^3(1-q^2)(1-q^3) + q^4(1-q^2)(1-q^3) + q^1(1-q^2)(1-q^3) + q^6(1-q^1)(1-q^3) \\
+q^3(1-q^1)(1-q^3)+q^0(1-q^1)(1-q^3)
\end{multlined} \\
\label{eq: 2nd qanalogue rpp} &= \begin{multlined}[t]
q^0(1-q^2)(1-q^3) + q^2(1-q^2)(1-q^3) + q^3(1-q^2)(1-q^3) + q^4(1-q^1)(1-q^3)  \\
+q^5(1-q^1)(1-q^3)+q^6(1-q^1)(1-q^3).
\end{multlined}
\end{align}
\end{example}

We conclude this section by  showing that Theorem \ref{thm:oorpptableaux} is a $q$-analogue of the Okounkov--Olshanski formula by using a result from Stanley's theory of \vocab{$P$-partitions} and taking the limit $q \to 1$.

\begin{theorem}[Stanley \cite{stanley1972}] \label{thm:ppartitions}
\[\rpplmq = \frac{P(q)}{(1-q)(1-q^2) \cdots (1-q^{\abs{\lm}})}\]
where $P(q)$ is a polynomial with $P(1) = f^\lm$.
\end{theorem}

\begin{proposition} \label{prop: q-analogue imply oof}
Theorems~\ref{thm:oorpptableaux} and \ref{thm:oorpptalt} each imply the Okounkov--Olshanski formula \eqref{eq:oof}.
\end{proposition}

\begin{proof}
By Theorem \ref{thm:ppartitions}, we have
\[f^\lm = \abs{\lm}! \lim_{q \to 1} (1-q)^{\abs{\lm}} \cdot  \rpplmq.\]
Therefore by either the formula for $\rpplmq$ in Theorem~\ref{thm:oorpptableaux} or Theorem~\ref{thm:oorpptalt} we have that, 
\begin{multline*}
\frac{f^\lm}{f^\lambda} = \frac{\abs{\lm}!}{\abs{\lambda}!} \lim_{q \to 1} (1-q)^{-\abs{\mu}} \cdot \rppratio  = \frac{\abs{\lm}!}{\abs{\lambda}!} \sum_{T \in \SSYT(\mu, d)} \prod_{u \in [\mu]} w(u,T(u)) \\ = \frac{\abs{\lm}!}{\abs{\lambda}!} \sum_{T \in \SSYT(\mu, d)} \prod_{u \in [\mu]} (\lambda_{d+1-T(u)} - c(u)).
\end{multline*}
Using the hook length formula to rewrite $f^\lambda$ indeed yields the Okounkov--Olshanski formula.
\end{proof}

\subsection{An excited diagram reformulation} \label{subsec: q-analogue oof excited diagram}
Observe that the nonzero terms in Theorems \ref{thm:oorpptableaux} and \ref{thm:oorpptalt} are those with $w(u, T(u)) \neq 0$ for all $u$, i.e.\ $\cOOT(\lm)$. For the remainder of this section we will apply the correspondences of Section \ref{subsec:weighttrans} in order to reformulate Theorems \ref{thm:oorpptableaux} and \ref{thm:oorpptalt}. Instead of applying every correspondence to both theorems, which would yield at least eight reformulations, we will only highlight those which relate to existing constructions in the literature or those which are especially elegant.

First, we rewrite Theorem~\ref{thm:oorpptableaux} in terms of the excited diagrams $\cOOE(\lm)$. The following definition is useful:
\begin{definition}[\cite{ikedanaruse2013grothendieck,mpp1}]
Given an excited diagram $D \in \cale(\lm)$, we define its \vocab{excited peaks} $\EP(D)$ recursively by defining $\EP([\mu]) = \varnothing$. Furthermore, if an excited move takes $D$ to $D'$ by moving a cell of $D$ out of cell $u$, $\EP(D') = (\EP(D) \cup \set{u}) \setminus \set{u_\rightarrow, u_\downarrow}$, where $u_\rightarrow$ and $u_\downarrow$ are the cells immediately to the right and immediately below $u$, respectively. This is represented in the following diagram, where a gray triangle represents an excited peak:
\[\ytab[centertableaux]{{\dg}{},{}{}} \; \longrightarrow\; \ytab[centertableaux]{{\peak}{\del},{\del}{\dg}} \]
\end{definition}

Although there may be multiple ways to reach an excited diagram using excited moves, $\EP(D)$ is well-defined \cite{ikedanaruse2013grothendieck}.

\begin{theorem} \label{thm:oorppexcited}
	\[\rppratio =  \sum_{D \in \cOOE(\Lambda^d/\mu)} q^{p(D)} \prod_{u \in D} (1 - q^{w(u)}),\]
	where $w(i,j)=\lambda_{d+1-i}+i-j$ and $p(D) = \sum_{u \in \EP(D)} w(u)$.
\end{theorem}

\begin{example}
Figure~\ref{fig:ex_oo_excited} has the excited diagrams $\cOOE(2221/11)$ with the cells with excited peaks marked with gray. Calculating $p(D)$ for each of these diagrams gives the numerators in the RHS of \eqref{eq: 1st qanalogue rpp}.
\end{example}

The rest of this subsection is devoted to the proof of Theorem~\ref{thm:oorppexcited}.

\begin{proposition} \label{prop:excitedpeaks}
Given a tableau $T \in \cOOT(\lm)$, the excited peaks of the associated diagram $D$ in $\cOOE(\lm)$ are
\[\bigsqcup_{u \in [\mu]}\set{(s, s + c(u)) \mid m_T(u) \leq s < T(u)}.\]
\end{proposition}
\begin{proof}
Denote the set by $S(T)$. Say a cell $u$ is a \vocab{witness} of an element in $S(T)$ if it is $(s, s+c(u))$ for some $m_T(u) \leq s < T(u)$. If $u_1, u_2, \ldots$ are the cells of a given content going down the diagonal, we have $m_T(u_1) \leq T(u_1) < m_T(u_2) \leq T(u_2) < \cdots$, so each element in $S(T)$ has a unique witness, meaning that taking the disjoint union is valid. Also observe that $S(T)$ is disjoint from $D$.

We use induction on the excited diagram. The diagram $[\mu]$ corresponds to the tableau $T$ with $T(i,j) = i$; it is easy to show that in this case $m_T(i, j) = i = T(i, j)$, so $S(T)$ is empty, as desired.

Say an excited cell corresponding to $u \in [\mu]$ is excited from row $i$ to row $i+1$. Let $u_\uparrow$, $u_\rightarrow$, and $u_\downarrow$ be the cells adjacent to $u$ in the respective directions, if they exist. The corresponding tableau has $T(u)$ increase from $i$ to $i+1$. Since $m_T(u')$ depends only on entries the cells directly to the left and above $u'$ (and, in particular, not on $T(u')$), the only changes to $m_T(u')$ are that
\begin{itemize}
    \item $m_T(u_\rightarrow)$ increases to $i+1$ if it is originally $i$;
    \item $m_T(u_\downarrow)$ increases to $i+2$ if it is originally $i+1$.
\end{itemize}
Since $T(u_\uparrow) < i$ and $T(u_\rightarrow) \leq T(u)$ for both the initial and final tableau, if $(i, i + c(u) + 1)$ was in $S(T)$, it must have been witnessed by $u_\rightarrow$ and therefore must disappear. Similarly, $(i+1, i+c(u))$ must also be removed from $S(T)$, if it existed. Finally, since $T(u)$ increases by one, $(i+1, i+1+c(u))$ must be added to $S(T)$ as well, if it was not already there. (It is in fact very easy to prove that it could not have originally been there, but we will not use this fact.) Since the only changes to the pairs $(m_T(u), T(u))$ occur at $u$, $u_\rightarrow$, and $u_\downarrow$, there can be no more changes to $S(T)$.

Therefore, we have shown that $S(T)$ undergoes the same changes as $\EP(D)$, concluding the proof.
\end{proof}

\begin{proof}[Proof of Theorem~\ref{thm:oorppexcited}]
We rewrite Theorem \ref{thm:oorpptableaux} as the desired expression by Proposition \ref{prop:excitedpeaks}.
\end{proof}

\subsection{A reverse excited diagram formulation} \label{subsec: q-analogue oof reverse excited diagram}

Second, we rewrite Theorem~\ref{thm:oorpptalt} in terms of the reverse excited diagrams $\cRE(\lm)$.

\begin{theorem} \label{thm: 2nd q-analogue rev excited}
\[\rppratio = \sum_{D \in \cRE(\lm)} q^{p^*(D)}\prod_{u \in B(D)} (1 - q^{w(u)}) ,\]
where $w(i, j) = \lambda_i + d - j$ and $p^*(D)= \sum_{u \in F(D)} q^{w(u)}$ for $F(D) \subseteq [\lambda^*] \setminus D$ is the set of cells immediately to the right of an excited cell of $D$.
\end{theorem}

\begin{example}
Figure~\ref{fig:ex_reverse_excited} has the reverse excited diagrams $\cRE(2221/11)$ with the cells  in $F(D)$ marked with a purple left border. Calculating $p^*(D)$ for each of these diagrams gives the numerators in the RHS of \eqref{eq: 2nd qanalogue rpp}.
\end{example}

\begin{proof}
Fix a tableau $T \in \cOOT(\lm)$ with associated $D \in \cRE(\lm)$. We claim that
\[F(D) = \bigsqcup_{u \in [\mu]}\set{(d+1-s, c(u)+d) \mid T(u) < s \leq \widetilde M(u)},\]
where for this proof we let $\widetilde M(u) = \min(d, M_T(u))$. This claim would prove that
\[p^*(T) = \sum_{(i,j) \in F(D)} (\lambda_i +d- j),\]
meaning that the term corresponding to $T$ transforms to the term involving $D$, as we have previously shown that the normal weights transform accordingly as well (see Corollary \ref{cor:redform}). Thus it suffices to prove the claim.

We begin by noting that by similar reasoning to the proof of Proposition \ref{prop:excitedpeaks}, the sets in the disjoint union are indeed all disjoint.

Take a cell $u = (i,j)$ and some integer $T(u) < s \leq \widetilde M(u)$; we aim to show that the cell $(d+1-s, c(u) + d)$ is in $F(D)$. By Proposition \ref{prop:wherecellsgo}, $(d + 1 - T(u), c(u) + d) \in \mathsf{d}_j(D)$. Moreover, if $u_\downarrow$ is the cell below $u$, $(d + 1 - T(u_\downarrow), c(u) - 1 + d) \in \mathsf{d}_j(D)$ as well. If $u_\downarrow$ does not exist, then $(d+1-T(u), c(u)+d)$ is the highest member of $\mathsf{d}_j(D)$. In either case, we have $(d+1-t, c(u)+d-1) \in D$ for all $T(u) < t < T(u_\downarrow)$ if $u_\downarrow$ exists or for all $T(u) < t \leq d$ otherwise. Therefore, $(d+1-s, c(u)+d) \in D$.

Since $(d+1-T(u), c(u)+d) \in B(D) \subseteq [\lambda^*]$, $(d+1-s, c(u)+d) \in [\lambda^*]$ as well. Now suppose $(d+1-s, c(u)+d) \in D$. Then, the two points $(d+1-s, c(u)+d-1)$ and $(d+1-s, c(u)+d)$ correspond to entries of $d+1-s$ in the corresponding flagged tableau with contents $c(u)+s-d-1$ and $c(u)+s-d$. The only way for this to occur is for the cells to be horizontally adjacent, meaning that $(d+1-s, c(u)+d)$ must have been excited from column $j+1$. Now, since $(d+1-T(u), c(u)+d)$ is located somewhere below $(d+1-s, c(u)+d)$, the cells excited from column $j+1$ cannot extend to the bottom of $[\lambda^*]$, meaning that there must be some $(d+1-s', c(u)+d+1) \in \mathsf{d}_{j+1}(D)$ with $T(u) \leq s' < s$. This translates to an entry of $s'$ in $T$ with content $c(u)+1$, which contradicts the inequality $s \leq \widetilde M(u)$.

Conversely, consider $(d+1-s, c+d) \in F(D)$ for some $1\leq s \leq d$ and $c \in \setz$ and assume $(d+1-s, c+d-1)$ was excited from a cell $(i,j) \in [\lm]$. If it has not moved at all, then no cell to the right of it could have moved either, so $(d+1-s, c+d)$ is either in $D$ or not in $[\lambda^*]$, a contradiction. Therefore, there must exist some element $(d+1-t, c+d) \in \mathsf{d}_j(D)$ with $t < s$, which corresponds to an entry of $t$ in a cell $u$ with $c = c(u)$. Let $u_\rightarrow$ and $u_\downarrow$ be the cells immediately to the right or below $u$, if they exist.

We now want to show that $t < s \leq \widetilde M(u)$. The first inequality is clear from our discussion above, so we only need to show that $s \leq d$, $s \leq T(u_\rightarrow)$, $s \leq T(u_\downarrow)$, if these quantities exist. The first statement is obvious. To prove the second, note that $(d+1-T(u_\rightarrow), c(u)+d+1) \in \mathsf{d}_{j+1}(D)$. Therefore, if $T(u_\rightarrow) < s$, there would either have to be a cell $(d+1-s, c(u)+d+1)$ excited from column $j+1$, or a cell $(d+1-s', c(u)+d) \in \mathsf{d}_{j+1}(D)$ for some $T(u_\rightarrow) < s' \leq s$. Both cases are easily checked to be impossible.

Finally, $u_\downarrow$ induces an element $(d+1-T(u_\downarrow), c(u)+d-1) \in \mathsf{d}_j(D)$. In order for $(d+1-s, c(u)+d-1)$ to be in $D$, by Proposition \ref{prop:diagglobal} we must have $s < T(u_\downarrow)$ as desired.
\end{proof}

\section{Final remarks} \label{sec: final remarks}

\subsection{Tableaux from lozenge tilings} \label{sec: symmetries lozenge tilings} 
Section~\ref{sec:posterms} introduced the main bijections and several tableaux objects  related to the terms in the Okounkov--Olshanski formula coming from lozenge tilings. Given a lozenge tiling $\tiling$ in $L(\nabla^*_\lm)$, there are twelve different ways to construct a corresponding semistandard tableau: given one of the three path systems:
$(\rho'_1,\ldots,\rho'_d)$, $(\beta'_1,\ldots,\beta'_d)$, $(\gamma'_1,\ldots,\gamma'_d)$; there are two choices regarding how to orient the paths and two choices regarding which of the step directions should be interpreted as representing an entry. In practice, reversing the orientation essentially rotates the tableau by a half turn. Also, as suggested above, a tableau generated this way has its column data encoded in a different path system.

As a result, there are only three ``conjugate'' tableaux that can be generated with nontrivial differences from each other, of which our analysis of Okounkov--Olshanski tableaux in Sections~\ref{subsec:OOtoLT},\ref{subsec:LTtoFT} has mentioned two (the third appears to be highly technical and thus not particularly useful). This notion of conjugacy appears promising for not only the analysis of tableaux, but also of lattice paths.

A more familiar context for this idea is that of {\em plane partitions}. Given the standard correspondence between plane partitions and sets of unit cubes in $3$-space, rotation of the axes induces two conjugates to each plane partition . While these operations have been a classical topic of study (e.g. see \cite{stanley-symm-pp}), it does not appear that the analogue for semistandard tableaux is as established. 

\subsection{Flagged skew tableaux and dual stable Grothendieck polynomials} \label{sec:dual stable}

Damir Yeliussizov (private communication) pointed out that the flagged tableaux in $\cSF(\lm)$ appear in \cite{LP} as \emph{elegant fillings} of skew shapes $\lambda/\mu$ and their number $\OOT(\lm)$ gives the Schur expansion of \emph{dual stable Grothendieck polynomials} denoted by $g_{\lambda}({\bf x})$. That is  \cite[Theorem 9.8]{LP} states  that 
\[
g_{\lambda}({\bf x}) = \sum_{\mu} \OOT(\lm) \cdot s_{\mu}.
\]
It would be interesting to understand this connection further.

\subsection{Similarities with the Naruse formula} \label{sec:finalnaruse}
In 2014, Naruse announced a new formula for $f^\lm$ involving a sum indexed by excited diagrams $\cale(\lm)$. 

\begin{theorem}[Naruse \cite{naruse2014}]
\begin{equation}
\tag{NHLF}
\label{eq:nhlf}
f^\lm = \frac{\abs{\lm}!}{\prod_{u \in [\lambda]} h(u)} \sum_{D \in \cale(\lm)} \prod_{u \in D} h(u).
\end{equation}
\end{theorem}

It is interesting to compare these two formulas \eqref{eq:nhlf} and \eqref{eq:oof} to compute $f^{\lm}$. For instance, both formulas come from evaluations of factorial Schur functions \cite[\S 9.4]{mpp1}. Also, with the appropriate reformulation of \eqref{eq:oof} (Corollary~\ref{cor:edform} or \ref{cor:redform}) both their terms are indexed by some sort of excited diagrams. In addition, Morales, Pak, and Panova \cite{mpp1} gave a $q$-analogue enumerating reverse plane partitions that was later extended to \emph{$d$-complete posets} by Okada and Naruse \cite{naruseokada2018}.

\begin{theorem}[Morales-Pak-Panova \cite{mpp1}]
\[\rppratio = \sum_{D \in \cale(\lambda/\mu)} q^{a'(D)} \prod_{u \in D}(1-q^{h(u)})\]
where $a'(D) = \sum_{u \in \EP(D)} h(u)$.
\end{theorem}
This $q$-analogue of the Naruse hook length formula bears a strong resemblance to one of the $q$-analogues of \eqref{eq:oof}, Theorem~\ref{thm:oorppexcited}. However, unlike the Okounkov--Olshanski formula, the Naruse formula is strongly geometric in the sense of Section \ref{sec: oof geometric?}.

\subsection{Proportionality for zigzags and reverse hooks}
It is also interesting to compare \eqref{eq:oof} and \eqref{eq:nhlf} on certain families of skew shapes. For slim shapes and rectangles the number of terms of both formulas are equal but the term-by-term weights can differ (Section~\ref{subsec: slim shapes and rects}). Also for \vocab{thick reverse hooks} $\lm=(b+c)^{a+c}/b^a$ we have that $f^{\lm}$ is proportional to $\OOT(\lm)=|\cale(\lm)|=M(a,b,c)$, where $M(a,b,c)$ denotes the MacMahon box formula. Namely by \cite[Eq. (1.1)]{mpp3} and \eqref{eq:oof_is_MacMahon} we have that \[
f^{(b+c)^{a+c}/b^a} \,=\,  n! \cdot  M(a,b,c) \cdot \frac{\Phi(c+1)\Phi(a+b+c+1)}{\Phi(a+b+2c+1)},
\]
where $n = |\lm|$. Surprisingly, by Corollary~\ref{cor: rel SYTzigzag OOTzigzag} this proportionality phenomenon also holds for thick zigzags $\lm = \delta_{n+2k}/\delta_n$ and $\OOT(\lm)$. For the case of thick reverse hooks this proportionality can be proved by a \vocab{hidden symmetry} of a multivariate extension of factorial Schur functions \cite[\S 3]{mpp3}. It would be interesting to understand the formula \eqref{eq:oof} for thick zigzag shapes in a similar way.

\subsection{Combinatorial proofs}
The hook length formula for $f^{\lambda}$ has a celebrated bijective proof by Novelli-Pak-Stoyanovskii \cite{NPS}. Also, a recurrence equivalent to Naruse's formula \eqref{eq:nhlf} has bijective proofs by Konvalinka \cite{konvalinka1,konvanlinka2}. It would be interesting to find bijective proofs of \eqref{eq:oof}.

In \cite[\S 6,7]{mpp1} there are bijective arguments  for the SSYT and RPP $q$-analogues of the Naruse hook length formula using the {\em Hillman--Grassl correspondence} (eg. see \cite[\S 7.22]{EC2}). It would be of interest to study the relation between the latter (or its variants \cite{Sulz}) and the SSYT and RPP $q$-analogues of \eqref{eq:oof}.

\subsection{An Okounkov--Olshanski formula for skew \texorpdfstring{$d$}{d}-complete posets}

In \cite{naruse2014} Naruse also announced positive formulas for the number standard tableaux of {\em skew shifted shapes} (see \cite[\S 8.5]{mpp2}). Both this formula and \eqref{eq:nhlf} were generalized by Naruse and Okada \cite{naruseokada2018} to {\em skew $d$-complete posets}. These posets are obtained from Proctor's $d$-complete posets (e.g. see \cite{d-complete_survey}) by removing a certain order filter. It would be interesting to find analogues of \eqref{eq:oof} for skew shifted shapes and for skew $d$-complete posets.

\subsection{Asymptotics for \texorpdfstring{$f^\lm$}{f\^{}(λ/μ)}}
In Section~\ref{sec:asymptotic thick zigzags} we gave asymptotic applications of the Okounkov--Olshanski formula for thick zigzags $\delta_{2k}/\delta_{k}$. In \cite{stanley_asymptotic_skew} Stanley used \eqref{eq:oof} to give asymptotic results for $f^{\lambda/\mu}$ when both $\lambda,\mu$ have the \emph{Thoma–Vershik–Kerov} limit (see \cite[\S 5]{mpp_asymptotics}).
Another approach to give bounds on $f^\lm$ is using bounds on characters \cite{stanley_asymptotic_skew, CGS, DF}. 

Since \eqref{eq:oof} and \eqref{eq:nhlf} are nonnegative formulas for $f^\lm$, we can derive lower and upper bounds for $f^\lm$ from them. This was done for Naruse's formula \eqref{eq:nhlf} in \cite{mpp_asymptotics} and we do it for \eqref{eq:oof} next. For a skew shape $\lm$ let 
\[
G(\lm)\coloneqq \frac{n!}{\prod_{u \in [\lambda]} h(u)} \prod_{(i,j)\in B([\lm])} (\lambda_i-j+d),
\]
where $B([\lm])$ are the cells of the broken diagonals of $[\lm]$, i.e. the cells of the diagonals $\mathsf{d}_1,\ldots,\mathsf{d}_{\ell(\mu')}$ of associated to the reverse excited diagram $[\lm]$. 

\begin{corollary}
For a skew shape $\lm$ we have that 
\begin{equation}
G(\lm) \leq f^{\lm} \leq \OOT(\lm) \cdot G(\lm).
\end{equation}
\end{corollary}

\begin{proof}
We use the excited diagram formulation Corollary~\ref{cor:redform} of the Okounkov--Olshanski formula. The lower bound follows from $\eqref{eq:oof}$ since $G(\lm)$ is one of the terms of the sum. To show the upper bound it suffices to show that the product $\prod_{u \in B(D)} (\lambda_i-j+d)$ of arm lengths is maximal for the reverse excited diagram $D=[\lm]$. To see this note that with each reverse excited move this product decreases, since the broken diagonals go down by column. 
\end{proof}

It would be interesting to study these general bounds coming from \eqref{eq:oof}. Recently, the bounds from the Naruse formula in \cite{mpp_asymptotics} were used in \cite{CPP} to bound \vocab{sorting probabilities} for standard tableaux of skew shape.

\subsection{Equivariant semistandard tableaux} \label{subsec:eqssyt}
Thomas and Yong \cite{thomasyong2012} define an \vocab{equivariant semistandard tableau} of shape $\theta$ to be a filling of the cells of $[\theta]$ with positive integers and a labeling of the horizontal edges in $[\theta]$ with a possibly empty set of positive integers such that in each row the entries in the cells are weakly increasing and in each column the entries and the edge labels are strictly increasing.

As mentioned in Section \ref{subsec:weighttrans}, the members of $\cSF(\lm)$ are in bijection with equivariant semistandard tableaux of shape $\lm$ so that all entries from $1$ to $\lambda'_i$ appear in the $i$th column of $\lm$. At the same time, there exists a rule for $f^{\lm}$ indexed by equivariant \emph{standard} tableaux that rectify to the \vocab{row superstandard tableau} of shape $\lambda$ under the process of \vocab{equivariant jeu de taquin} defined in \cite{thomasyong2012} (see also \cite[\S 9.4.2]{mpp1}).

These terms are not exactly the same; in the case of $\lm = 22/1$, the tableaux in bijection with $\cSF(\lm)$ are
\[
\begin{tikzpicture}[x=0.45cm,y=0.45cm,scale=0.5,font=\footnotesize,baseline={(0,0)}]
\draw (0, 0) rectangle (4,2) (2,0) rectangle (4,4);
\node at (1, 1) {1};
\node at (3, 1) {2};
\node at (3, 3) {1};
\node at (1, 0) {2};
\end{tikzpicture}
\quad\text{and}\quad
\begin{tikzpicture}[x=0.45cm,y=0.45cm,scale=0.5,font=\footnotesize,baseline={(0,0)}]
\draw (0, 0) rectangle (4,2) (2,0) rectangle (4,4);
\node at (1, 1) {2};
\node at (3, 1) {2};
\node at (3, 3) {1};
\node at (1, 2) {1};
\end{tikzpicture},
\]
with weights $2$ and $2$ respectively while the relevant equivariant standard tableaux are
\[
\begin{tikzpicture}[x=0.45cm,y=0.45cm,scale=0.5,font=\footnotesize,baseline={(0,0)}]
\draw (0, 0) rectangle (4,2) (2,0) rectangle (4,4);
\node at (1, 1) {1};
\node at (3, 1) {4};
\node at (3, 3) {2};
\node at (1, 0) {3};
\end{tikzpicture},
\quad
\begin{tikzpicture}[x=0.45cm,y=0.45cm,scale=0.5,font=\footnotesize,baseline={(0,0)}]
\draw (0, 0) rectangle (4,2) (2,0) rectangle (4,4);
\node at (1, 1) {3};
\node at (3, 1) {4};
\node at (3, 3) {2};
\node at (1, 2) {1};
\end{tikzpicture},
\quad
\begin{tikzpicture}[x=0.45cm,y=0.45cm,scale=0.5,font=\footnotesize,baseline={(0,0)}]
\draw (0, 0) rectangle (4,2) (2,0) rectangle (4,4);
\node at (1, 1) {3};
\node at (3, 1) {4};
\node at (3, 3) {1};
\node at (3, 2) {2};
\end{tikzpicture},\quad
\text{and}\quad
\begin{tikzpicture}[x=0.45cm,y=0.45cm,scale=0.5,font=\footnotesize,baseline={(0,0)}]
\draw (0, 0) rectangle (4,2) (2,0) rectangle (4,4);
\node at (1, 1) {1};
\node at (3, 1) {3};
\node at (3, 3) {2};
\node at (3, 0) {4};
\end{tikzpicture}.
\]
with weights $2,2,0$, and $0$ respectively.

However, the flagged tableaux formulation of the Okounkov-Olshanski formula and the equivariant tableaux formula are term-by-term equivalent once trivial terms are removed. The bijection between terms proceeds as follows:
\begin{itemize}
\item Given an equivariant tableau $T$ of shape $\lm$ that rectifies to a tableau of shape $\lambda$, we claim that we may ignore the cases where a label ever moves horizontally. To see this, suppose that a label moves from column $i + 1$ to column $i$. This implies that the number of labels in $T$ strictly to the right of column $i$ exceeds the sum $\sum_{j = i+1}^{\ell(\lambda')} \lambda'_j$. However, this implies that after rectifying the columns to the right of column $i$, there is at least one edge label to the right of column $i$. This means that the weight associated with the tableau is zero, so we may safely ignore this case.
\item As a result, there are exactly $\lambda'_i$ labels in column $i$ of $T$, which are equal to the entries in column $i$ of the superstandard tableau of shape $\lambda$. For each $i$, relabel the entries in column $i$ with the numbers $1, 2, \ldots, \lambda'_i$, preserving their relative order.
\item Finally, remove all the edge labels. Since $T$ is a standard equivariant tableau, the remaining tableau $T'$ is a semistandard tableau of shape $\lm$. Moreover, it is in $\cSF(\lm)$, since no entry in column $i$ exceeds $\lambda'_i$.
\end{itemize}
It is straightforward to check that this bijection is reversible: given an element of $\cSF(\lm)$, edge labels of the ``missing entries'' in each column can be added, and reversing the renumbering only requires knowledge of $\lambda$, not of the original tableau. The resulting tableau can be verified to be standard and rectifying to the superstandard tableau. Therefore, this bijection $T \mapsto T'$ exists.

The weight of a valid equivariant tableau \cite[\S 9.4.2]{mpp1} is $\prod_{\text{edge labels }\mathfrak{l}} (\lambda_{i_{\mathfrak{l}}} - i_{\mathfrak{l}} + 1 - c(u_{\mathfrak{l}}))$, where, each label $\mathfrak{l}$ is originally immediately under cell $u_{\mathfrak{l}}$ and rectifies to row $i_{\mathfrak{l}}$. Under our bijection, edge labels correspond to elements of $M(T')$. Specifically, an element $(j, t) \in M(T')$ corresponds to a label $\mathfrak{l}$ with $u_{\mathfrak{l}} = (i(j, t), j)$ and $i_{\mathfrak{l}} = t$. This proves the term-by-term equivalence.

\subsection{On the use of shifted shapes}
In this paper, we have used shifted Young diagrams in order to define reverse excited diagrams. However, in some sense, this usage is somewhat superficial, since we do not make use of any deep result from the theory of shifted Young tableaux. In other words, the definition of $\cSF(\lm)$ would work equally well if the shape $\lambda^*$ were replaced with a band of height $d$ extending infinitely far to the left. It would be interesting to see if there is a deeper connection between reverse excited diagrams and shifted shapes, beyond the geometric restriction on where the excited cells can move.



\section*{Acknowledgements}
The formulations of Theorem~\ref{thm:ooiskt}, and Theorem~\ref{thm:genocchi} were first conjectured by Igor Pak, Greta Panova and the first author in the course of work on \cite{mpp1,mpp2}. They were subsequently proved in this paper.  We are very grateful to them for graciously allowing us to study and settle these conjectures, for their encouragement, and helpful suggestions throughout this project. 

This research was made possible by MIT-PRIMES program of 2018. We thank Christian Krattenthaler and an anonymous referee for the idea of the proof of Lemma~\ref{lem:ordproductq} that simplified a previous argument; and another anonymous referee for clarifications and improvements to Section~\ref{subsec:eqssyt}, as well as thank Pavel Etingof, Claude Eicher, Angèle Foley, Igor Pak, Damir Yeliussizov, Alex Yong, and the anonymous referees for helpful comments and suggestions.


%

\bibliographystyle{alphaurl}

\end{document}